\documentclass[reqno,nosumlimits,12pt]{amsart}
\usepackage{amssymb}
\usepackage{color,float}
\usepackage{graphicx,hyperref}
\usepackage[all]{xy}

\newtheorem{theorem}{\sc Theorem}[section]
\newtheorem{proposition}[theorem]{\sc Proposition}
\newtheorem{lemma}[theorem]{\sc Lemma}
\newtheorem{corollary}[theorem]{\sc Corollary}
\theoremstyle{definition}
\newtheorem{definition}[theorem]{\sc Definition}

\theoremstyle{remark}
\newtheorem{remark}[theorem]{\sc Remark}

\newtheorem{noname}[theorem]{\sc\hspace*{-1ex}}
\newtheoremstyle{mystyle}{}{}{\slshape}{}{\scshape}{}{ }{\thmname{#1}\ \textnormal{\thmnote{(#3).}}}
\theoremstyle{mystyle}

\newtheoremstyle{mystyle1}{}{}{\slshape}{}{\scshape}{}{ }{\thmname{#1}.}
\theoremstyle{mystyle1}

\begin{document}
\title[Coactions on Hochschild Homology of Hopf-Galois Extensions]{Coactions on Hochschild
Homology of Hopf-Galois Extensions and Their Coinvariants}
\author{A. Makhlouf and D. \c{S}tefan}
\thanks{2000 \textit{Mathematics Subject Classification}. Primary 16E40;
Secondary 16W30}
\thanks{D. \c Stefan was financially supported by CNCSIS, Contract
560/2009 (CNCSIS code ID\_69).}
\date{}
\keywords{Hopf-Galois extensions, Hochschild homology and cyclic homology}

\begin{abstract}
Let ${\mathcal{B}\subseteq \mathcal{A}}$ be an $\mathcal{H}$-Galois
extension, where $\mathcal{H}$ is a Hopf algebra over a field $\mathbb{K}$.
If $M$ is a Hopf bimodule then $\mathrm{HH}_{\ast }(\mathcal{A},M)$, the
Hochschild homology of $\mathcal{A}$ with coefficients in $M$, is a right
comodule over the coalgebra $\mathcal{C}_{\mathcal{H}}=\mathcal{H}/[\mathcal{%
H},\mathcal{H}]$. Given an injective left $\mathcal{C}_{\mathcal{H}}$%
-comodule $V$, our aim is to understand the relationship between $\mathrm{HH}%
_{\ast }(\mathcal{A},M)\square_{\mathcal{C_H}}V$ and $\mathrm{HH}_{\ast }(%
\mathcal{B},M\square_{\mathcal{C_H}}V)$. The roots of this problem can be
found in \cite{Lo2}, where $\mathrm{HH}_{\ast }(\mathcal{A},\mathcal{A})^G$
and $\mathrm{HH}_{\ast }(\mathcal{B},\mathcal{B})$ are shown to be
isomorphic for any centrally $G$-Galois extension. To approach the above
mentioned problem, in the case when $\mathcal{A}$ is a faithfully flat $%
\mathcal{B}$-module and $\mathcal{H}$ satisfies some technical conditions,
we construct a spectral sequence
\begin{equation*}
\mathrm{Tor}_p^{\mathcal{R_H}}\left(\mathbb{K},\mathrm{HH}_q(\mathcal{B}%
,M\square_{\mathcal{C_H}}V)\right)\Longrightarrow \mathrm{HH}_{p+q}(\mathcal{%
A},M)\square_{\mathcal{C_H}}V,
\end{equation*}
where $\mathcal{R_H}$ denotes the subalgebra of cocommutative elements in $%
\mathcal{H}$. We also find conditions on $\mathcal{H}$ such that the edge
maps of the above spectral sequence yield isomorphisms
\begin{equation*}
\mathbb{K}\otimes_{\mathcal{R_H}}\mathrm{HH}_\ast(\mathcal{B},M\square_{%
\mathcal{C_H}}V)\cong \mathrm{HH}_{\ast}(\mathcal{A},M)\square_{\mathcal{C_H}%
}V.
\end{equation*}
In the last part of the paper we define centrally Hopf-Galois extensions
and we show that for such an extension $\mathcal{B}\subseteq\mathcal{A}$,
the $\mathcal{R_H}$-action on $\mathrm{HH}_\ast(\mathcal{B},M\square_{
\mathcal{C_H}}V)$ is trivial. As an application, we compute the subspace of $%
\mathcal{H}$-coinvariant elements in $\text{HH}_{\ast }(\mathcal{A},M)$. A
similar result is derived for $\text{HC}_{\ast }(\mathcal{A})$, the cyclic
homology of $\mathcal{A}$.
\end{abstract}

\maketitle

\section*{Introduction}

Chase and Sweedler \cite{CS}, more than 35 years ago, defined a special case
of Hopf-Galois extensions, similar to the theory of Galois group actions on
commutative rings that had been developed by Chase, Harrison and Rosenberg
\cite{CHR}. The general definition, for arbitrary Hopf algebras, is due to
Takeuchi and Kreimer \cite{KT}. Besides Galois group actions that we
already mentioned, strongly graded algebras are Hopf-Galois extensions.
Other examples come from the theory of affine quotients and of enveloping
algebras of Lie algebras.

By definition, an $\mathcal{H}$-Galois extension is given by an algebra map $%
\rho _{\mathcal{A}}:\mathcal{A}\rightarrow \mathcal{A}\otimes \mathcal{H}$
that defines a coaction of a Hopf algebra $\mathcal{H}$ on the algebra $%
\mathcal{A}$ such that the map $\ $
\begin{equation*}
\;\beta _{\mathcal{A}}:\mathcal{A}\otimes _{\mathcal{B}}\mathcal{A}%
\rightarrow \mathcal{A}\otimes \mathcal{H},\quad \beta _{A}=\left( m_{%
\mathcal{A}}\otimes \mathcal{H}\right) \circ \left( \mathcal{A}\otimes _{%
\mathcal{B}}\rho _{\mathcal{A}}\right)
\end{equation*}%
is bijective. Here $m_{\mathcal{A}}:\mathcal{A}\otimes \mathcal{A}%
\rightarrow \mathcal{A}$ denotes the multiplication map in
$\mathcal{A}$ and $\mathcal{B}$ is the subalgebra of coinvariant
elements, i.e. of all $a\in A\ $ such that $\rho
_{\mathcal{A}}(a)=a\otimes 1.$  For the remaining part of
the introduction we fix an $\mathcal{H}$-Galois extension $\mathcal{B}%
\subseteq \mathcal{A}$. If there is no danger of confusion, we
shall say that $\mathcal{B}\subseteq\mathcal{A}$ is a Hopf-Galois
extension.

For an $\mathcal{A}$-bimodule $M,$ let $\mathrm{HH}_{\ast }(\mathcal{A},M)$
denote Hochschild homology of $\mathcal{A}$ with coefficients in $M.$ Since $%
\mathcal{B}$ is a subalgebra of $\mathcal{A},$ the Hochschild homology of $%
\mathcal{B}$ with coefficients in $M$ also makes sense. A striking
feature of $\mathrm{HH}_{\ast }(\mathcal{B},M)$ is that
$\mathcal{H}$ acts on these linear spaces via the Ulbrich-Miyashita
action, cf. \cite{St1}. By
taking $M$ to be a Hopf bimodule, more structure can be defined not only on $\mathrm{%
HH}_{\ast }(\mathcal{B},M)$ but on $\mathrm{HH}_{\ast
}(\mathcal{A},M)$ too. Recall that $M$ is a Hopf bimodule if $M$
is an $\mathcal{A}$-bimodule and an $\mathcal{H}$-comodule such
that the maps that define the module
structures are $\mathcal{H}$-colinear. By definition of Hopf modules, $%
\mathrm{HH}_{0}(\mathcal{B},M)$ is a quotient $\mathcal{H}$-comodule of $M$.
This structure can be extended to an $\mathcal{H}$-coaction on $\mathrm{HH}%
_{n}(\mathcal{B},M),$ for every $n.$ On the other hand, $\mathrm{HH}_{0}(%
\mathcal{A},M)$ is not an $\mathcal{H}$-comodule, in general. Nevertheless,
the quotient coalgebra $\mathcal{C}_{\mathcal{H}}:=\mathcal{H}/[\mathcal{H},%
\mathcal{H}]$ coacts on $\mathrm{HH}_{\ast }(\mathcal{A},M),$ where $[%
\mathcal{H},\mathcal{H}]$ denotes the subspace spanned by all
commutators in $\mathcal{H}$, see \cite{St2}.

These actions and coactions played  an important role in the study
of Hochschild (co)homology of Hopf-Galois extensions, having deep
application in the field. Let us briefly discuss some of them,
that are related to the present work. First, in degree zero, the
coinvariants of Ulbrich-Miyashita action
\begin{equation*}
\mathrm{HH}_{0}(\mathcal{B},M)_{\mathcal{H}}:=\mathbb{K}\otimes _{\mathcal{H}%
}\mathrm{HH}_{0}(\mathcal{B},M)
\end{equation*}%
equal $\mathrm{HH}_{0}(\mathcal{A},M)$. This identification is one of the
main ingredients that are used in \cite{St1} to prove the existence of the
spectral sequence%
\begin{equation}
\mathrm{Tor}_{p}^{\mathcal{H}}(\mathbb{K},\mathrm{HH}_{q}(\mathcal{B}
,M))\implies \mathrm{HH}_{\ast }(\mathcal{A},M).  \label{ec:ss1}
\end{equation}%
It generalizes, in an unifying way, Lyndon-Hochschild-Serre spectral
sequence for group homology \cite[p. 195]{We}, Hochschild-Serre
spectral sequence for Lie algebra homology \cite[p. 232]{We} and
Lorenz spectral sequence for strongly graded algebras \cite{Lo1}. If
$\mathcal{H}$ is semisimple then (\ref{ec:ss1}) collapses and yields
the isomorphisms
\begin{equation*}
\mathrm{HH}_{n}(\mathcal{B},M)^{\mathcal{H}}\cong \mathrm{HH}_{n}(\mathcal{B}
,M)_{\mathcal{H}}\cong \mathrm{HH}_{n}(\mathcal{A},M),
\end{equation*}%
where $\mathrm{HH}_{n}(\mathcal{B},M)^{\mathcal{H}}$ denotes the space of $%
\mathcal{H}$-invariant Hochschild homology classes. Similar
isomorphisms in Hochschild cohomology were proved in the same way
and used, for example, to investigate algebraic deformations
arising from orbifolds with discrete torsion \cite{CGW}, to
characterize deformations of certain bialgebras \cite{MW} and to
study the $G$-structure on the cohomology of a Hopf algebra
\cite{FS}.

In \cite{St2}, for a subcoalgebra $C$ of $\mathcal{C}_{\mathcal{H}},$ a new
homology theory $\mathrm{HH}_{\ast }^{C}(\mathcal{A},-)$ with coefficients
in the category of Hopf bimodules is defined. In the case when $C$ is
injective as a left $\mathcal{C}_{\mathcal{H}}$-comodule we have
\begin{equation*}
\mathrm{HH}_{\ast }^{C}(\mathcal{A},M)\cong \mathrm{HH}_{\ast }(\mathcal{A}%
,M)\square _{\mathcal{C}_{\mathcal{H}}}C,
\end{equation*}%
for every Hopf bimodule $M$ (for the definition of the cotensor product $%
\square _{\mathcal{C}_{\mathcal{H}}}$ see the preliminaries of this paper).
Thus, $\mathrm{HH}_{\ast }^{C}(\mathcal{A},M)$ may be regarded as a sort of $%
C$-coinvariant part of $\mathrm{HH}_{\ast }(\mathcal{A},M).$ The main result
in loc. cit. is the spectral sequence%
\begin{equation}
\mathrm{Tor}_{p}^{\mathcal{H}}(\mathbb{K},\mathrm{HH}_{q}(\mathcal{B}%
,M\square _{\mathcal{C}_{\mathcal{H}}}C))\implies \mathrm{HH}_{\ast }^{C}(%
\mathcal{A},M),  \label{ec:ss2}
\end{equation}%
that exists for every Hopf bimodule $M,$ provided that $\mathcal{H}$ is
cocommutative and $C$ is injective as a left $\mathcal{C}_{\mathcal{H}}$%
-comodule.

Let $G$ be a group and let $\mathbb{K}$ be a field. If
$\mathcal{H}$ is the group algebra $\mathbb{K}G$, then
$\mathcal{C}_{\mathbb{K} G}:=\bigoplus_{\sigma \in T(G)}C_{\sigma
}$, where $T(G)$ denotes the set of conjugacy classes in $G$ and
$C_{\sigma }$ is a  subcoalgebra of dimension one, for every
$\sigma \in T(G).$ Hence, in this particular case,
the homogeneous components $\mathrm{HH}_{\ast }^{\sigma }(\mathcal{A},M):=%
\mathrm{HH}_{\ast }^{C_{\sigma }}(\mathcal{A},M)\ $completely
determine Hochschild homology of $\mathcal{A}$ with coefficients
in $M,$ cf. \cite{Lo1, St2}. For strongly graded algebras (i.e.
$\mathbb{K}G$-Galois extensions) and $C=C_{\sigma}$, the spectral
sequence (\ref{ec:ss2}) is due to Lorenz \cite{Lo1}. On the other
hand, Burghelea and Nistor defined and studied similar homogeneous
components of Hochschild and cyclic cohomology of group algebras
and crossed products in \cite{Bu} and \cite{Ni}, respectively.

We have already remarked that, for an arbitrary
$\mathcal{H}$-Galois extension $\mathcal{B}\subseteq \mathcal{A}$
and every Hopf bimodule $M$, Hochschild homology
$\mathrm{HH}_{\ast }(\mathcal{B},M)$ is a right
$\mathcal{H}$-comodule and a left $\mathcal{H}$-module.
In particular, $\mathcal{H}$ acts and coacts on $%
\mathcal{A}_{\mathcal{B}}:=\mathrm{HH}_{0}(\mathcal{B},\mathcal{A})$.
Notably, with respect to these structures,
$\mathcal{A}_{\mathcal{B}}$ is a stable-anti-Yetter-Drinfeld
$\mathcal{H}$-module (SAYD $\mathcal{H}$-module, for short). These
modules were independently discovered in \cite{JS} and \cite{HKRS},
and they can be thought of as coefficients for Hopf-cyclic homology.
In \cite{JS} the authors showed that Hopf-cyclic homology of
$\mathcal{H}$ with coefficients in $\mathcal{A}_{
\mathcal{B}}$ equals relative cyclic homology $\mathrm{HC}_{\ast }(\mathcal{%
A }/\mathcal{B}).$ This identification is then used to compute cyclic
homology of a strongly $G$-graded algebra with separable component of degree $%
1$ (e.g. group algebras and quantum tori). It is worthwhile to
mention that cyclic homology of a groupoid was computed in
\cite{BS}, using the theory of (generalized) SAYD modules.\bigskip

Let $G$ be a finite group of\ automorphisms of an algebra $\mathcal{A}$ over
a field $\mathbb{K}\ $and let ${\mathcal{A}}^{G}$ denote the ring of $G$%
-invariants in $\mathcal{A}$. Since $G$ is finite, the dual vector space $(%
\mathbb{K}G)^{\ast }$ has a canonical structure of Hopf algebra and $%
\mathcal{A}$ is a $(\mathbb{K}G)^{\ast }$-comodule algebra. Clearly, the
coinvariant subalgebra with respect to this coaction equals $\mathcal{A}^G$.
It is well-known that $\mathcal{A}^G\subseteq \mathcal{A}$ is $(\mathbb{K}
G)^{\ast }$-Galois if and only if this extension is Galois in the sense of
\cite{CHR}. The center $\mathcal{Z}$ of $\mathcal{A}$ is $G$-invariant.
Following \cite{Lo2}, we say that ${\mathcal{A}}^{G}\subseteq\mathcal{A}$ is
centrally Galois if ${\mathcal{Z}}^{G}\subseteq\mathcal{Z}$ is $(\mathbb{K}
G)^{\ast }$-Galois.

The Galois group $G$ acts, of course both on Hochschild homology $\mathrm{HH%
}_{\ast }(\mathcal{A},\mathcal{A})$ and cyclic homology $\mathrm{HC}_{\ast }(%
\mathcal{A})$. To simplify the notation, we shall write $\mathrm{HH}_{\ast }(%
\mathcal{A})$ for $\mathrm{HH}_{\ast }(\mathcal{A},\mathcal{A})$.
By \cite[\S 6]{Lo2}, for a centrally Galois extension
$\mathcal{A}^G\subseteq \mathcal{A}$,
\begin{equation}
\mathrm{HH}_{\ast }(\mathcal{A})^{G}\cong \mathrm{HH}_{\ast
}(\mathcal{A}^G) \label{izo_intro}
\end{equation}%
and a similar isomorphism exists in cyclic homology, provided that
the order of $G$ is invertible in $\mathbb{K}$.

Since $(\mathbb{K}G)^{\ast }$ is commutative, the coalgebras $\mathcal{C}_{(%
\mathbb{K}G)^{\ast }}\ $and $(\mathbb{K}G)^{\ast }$ are equal. The category
of left $\mathbb{K}G$-modules is isomorphic to the category of right $(%
\mathbb{K}G)^{\ast }$-comodules and through this identification
$X^{G}\cong X\square _{(\mathbb{K}G)^{\ast }}\mathbb{K}.$  In
particular,
\begin{equation*}
\mathrm{HH}_{\ast }(\mathcal{A})^{G}\cong \mathrm{HH}_{\ast }(\mathcal{A}%
)\square _{(\mathbb{K}G)^{\ast }}\mathbb{K}.
\end{equation*}%
This isomorphism suggests that the main result in \cite{Lo2} might be
approached in the spirit of \cite{St2}, i.e. using the theory of Hopf-Galois
extensions and an appropriate spectral sequence that converges to $\mathrm{HH}%
_{\ast }(\mathcal{A})\square _{(\mathbb{K}G)^{\ast }}\mathbb{K}.$ Since in
general $(\mathbb{K}G)^{\ast }$ is not cocommutative, the spectral sequence
in (\ref{ec:ss2}) cannot be used directly. The main obstruction to extend it
for a not necessarily cocommutative Hopf algebra $\mathcal{H}$, is the fact
that Ulbrich-Miyashita action does not induce an $\mathcal{H}$-action on $%
\mathrm{HH}_{\ast }(\mathcal{B},M\square _{\mathcal{C}_{\mathcal{H}}}C).$ To
overcome this difficulty we define
\begin{equation*}
\mathcal{R}_{\mathcal{H}}:=\ker (\Delta -\tau \circ \Delta ),
\end{equation*}%
where $\tau :\mathcal{H}\otimes \mathcal{H}\longrightarrow \mathcal{H}%
\otimes \mathcal{H}$ denotes the usual flip map. Since $\Delta $ and $\tau
\circ \Delta $ are morphisms of algebras, $\mathcal{R}_{\mathcal{H}}$ is a
subalgebra in $\mathcal{H}$. Moreover, if $\mathcal{A}$ is faithfully flat
as a left (or right) $\mathcal{B}$-module and the antipode of $\mathcal{H}$
is an involution, then we prove that $\mathrm{HH}_{\ast }(\mathcal{B}%
,M\square _{\mathcal{C}_{\mathcal{H}}}V)$ are left $\mathcal{R}_{\mathcal{H}}
$-modules, for every injective left $\mathcal{C}_{\mathcal{H}}$-comodule $V$%
; see Proposition \ref{pr:G_V}. Thus, under the above assumptions, for every
pair $(p,q)$ of natural numbers, it makes sense to define the vector spaces
\begin{equation}
E_{p,q}^{2}:=\mathrm{Tor}_{p}^{\mathcal{R_{H}}}\left( \mathbb{K},\mathrm{HH}%
_{q}(\mathcal{B},M\square _{\mathcal{C}_{\mathcal{H}}}V)\right) .
\label{ss1}
\end{equation}%
Furthermore, in Theorem \ref{te:sir} we prove that there is a
spectral sequence that has $E_{p,q}^{2}$ in the $(p,q)$-spot of
the second page and converges to $\mathrm{HH}_{\ast
}(\mathcal{A},M)\square _{\mathcal{H}}V$. This result follows as
an application of Proposition \ref{pr:sir spectral}, where we
indicate a new form of Grothendieck's spectral sequence. It also
relies on several properties of the Ulbrich-Miyashita action that
are proved in the first part of paper. Here, we just mention
equation (\ref{ec:SAYD}) that plays a key role, as it explains the
relationship between the module and comodule structures on
$\mathrm{HH}_{\ast }(\mathcal{B},M)$. More precisely, using the
terminology from \cite{HKRS, JS}, relation (\ref{ec:SAYD}) means
that Hochschild homology of $\mathcal{B}$ with coefficients in $M$
is an SAYD $\mathcal{H}$-module.

The most restrictive conditions that we impose in Theorem \ref{te:sir} are
the relations
\begin{equation}
\mathcal{R}_{\mathcal{H}}^{+}\mathcal{H}=\mathcal{H}^{+}\ \ \ \ \text{and~\
\ \ \ }\mathrm{Tor}_{n}^{\mathcal{R}_{\mathcal{H}}}(\mathbb{K},\mathcal{H)}%
=0,  \label{ec:cond}
\end{equation}%
for every $n>0,$ where $\mathcal{H}^{+}$ is the kernel of the counit and $%
\mathcal{R}_{\mathcal{H}}^{+}:=\mathcal{R}_{\mathcal{H}}\bigcap \mathcal{H}%
^{+}.$ The second relation in (\ref{ec:cond}) is easier to handle.
For example, if $\mathcal{H}$ is semisimple and cosemisimple over
a field of characteristic zero we show that
$\mathcal{R}_{\mathcal{H}}$ is semisimple and
$\mathcal{C}_{\mathcal{H}} $ is cosemisimple, cf. Proposition
\ref{pr:R_H-semisimple}. The proof of this result is based on the
identification $\mathcal{R}_{\mathcal{H}^{\ast
}}\cong \mathbb{K\otimes }_{\mathbb{Q}}C_{\mathbb{Q}}(\mathcal{H}),$ where $%
C_{\mathbb{Q}}(\mathcal{H})$ denotes the character algebra of $\mathcal{H},\
$and on the fact that $C_{\mathbb{Q}}(\mathcal{H})$ is semisimple if $%
\mathcal{H}$ is so. Thus in this case the second relation in
(\ref{ec:cond})
holds true. Obviously both relations in (\ref{ec:cond}) are satisfied if $%
\mathcal{H}$ is cocommutative. Notably, by Proposition \ref{pr:ss}, they are
also verified if $\mathcal{H}$ is semisimple and commutative. The other two
assumptions in Theorem \ref{te:sir} are not very strong. The antipode of $%
\mathcal{H}$ is involutive for commutative and cocommutative Hopf
algebras. In characteristic zero, by a result of Larson and
Radford, the antipode of a finite-dimensional Hopf algebra is an
involution if and only if the Hopf algebra is semisimple and
cosemisimple. Faithfully flat Hopf-Galois extensions are
characterized in \cite[Theorem 4.10]{SS}. In view of this
theorem, if the antipode of $\mathcal{H}$ is bijective, then an $\mathcal{H}$%
-Galois extension $\mathcal{B}\subseteq \mathcal{A}$ is faithfully
flat if and only if $\mathcal{A}$ is injective as an
$\mathcal{H}$-comodule. Hence, every $\mathcal{H}$-Galois
extension is faithfully flat if $\mathcal{H}$ is cosemisimple and
its antipode is bijective.

If the algebra $\mathcal{R}_{\mathcal{H}}$ is semisimple then the spectral
sequence in Theorem \ref{te:sir} collapses. We have already noticed that $%
\mathcal{R}_{\mathcal{H}}$ is semisimple if $\mathcal{H}$ is either
semisimple and cosemisimple over a field of characteristic zero, or
commutative and semisimple. In these situations, the edge maps induce
isomorphisms
\begin{equation}
\mathbb{K}\otimes _{\mathcal{R}_{\mathcal{H}}}\mathrm{HH}_{n}(\mathcal{B}%
,M\square _{\mathcal{H}}V)\cong \mathrm{HH}_{n}(\mathcal{A},M)\square _{%
\mathcal{H}}V.  \label{iso_general}
\end{equation}%
By specializing the isomorphism in (\ref{iso_general}) to the case $\mathcal{%
H}:=(\mathbb{K}G)^{\ast }$, we deduce the isomorphism in Corollary \ref%
{co:iso} that can be thought of as a generalization of
(\ref{izo_intro}).

Our result also explains why the isomorphism (\ref{izo_intro}) works for
centrally Galois extensions but not for arbitrary ones. Namely, the action
of $\mathcal{R}_{(\mathbb{K}G)^{\ast }}$ on $\mathrm{HH}_{n}(\mathcal{B)}$
is trivial for centrally Galois extensions, but not in general. In the last
part of the paper we show that a similar result holds for Hopf-Galois
extensions. Let $\mathcal{B}\subseteq \mathcal{A}$ be an $\mathcal{H}$%
-comodule algebra, where $\mathcal{H}$ is a commutative Hopf algebra. If the
center $\mathcal{Z}$ of $\mathcal{A}$ is an $\mathcal{H}$-subcomodule of $%
\mathcal{A}$ and $\mathcal{Z}\bigcap \mathcal{B}\subseteq \mathcal{Z}$ is a
faithfully flat $\mathcal{H}$-Galois extension, then we say that ${\mathcal{B%
}\subseteq \mathcal{A}}$ is centrally $\mathcal{H}$-Galois. We fix such an extension $%
{\mathcal{B}\subseteq \mathcal{A}}$. In view of Proposition
\ref{pr: trivial}, the $\mathcal{R}_{\mathcal{H}}$-action on
$\mathrm{HH}{}_{\ast }(\mathcal{B},M\square _{\mathcal{H}}V)$ is
trivial. If $\mathcal{H}$ is finite-dimensional and $\dim
\mathcal{H}$ is not zero in $\mathbb{K}$, then by Theorem
\ref{izo1}
\begin{equation*}
\mathrm{HH}{}_{\ast }(\mathcal{A},M)\square _{\mathcal{H}}V\simeq \mathrm{HH}%
{}_{\ast }(\mathcal{B},M\square _{\mathcal{H}}V),
\end{equation*}%
for every Hopf bimodule $M$ which is symmetric as a $\mathcal{Z}$-bimodule
and every left $\mathcal{H}$-comodule $V.$ Assuming that $\mathcal{H}:=(%
\mathbb{K}G)^{\ast }$ and that the order of $G$ is not zero in $\mathbb{K}$,
and taking $M:=\mathcal{A}$ and $V:=\mathbb{K}$ in the above
isomorphism, we obtain the isomorphism in (\ref{izo_intro}), cf. Corollary %
\ref{iso3}. Another application of Theorem \ref{izo1}  is given in Corollary %
\ref{iso2}.

Our approach has also the advantage that one can easily recover $\mathrm{HH}%
_{\ast }(\mathcal{A},M)$ from $\mathrm{HH}_{\ast }(\mathcal{B},M^{\mathrm{co}%
\mathcal{H}}).$ More precisely, in Theorem \ref{izo2}, we prove the
following isomorphism of $\mathcal{Z}$-modules and $\mathcal{H}$-comodules%
\begin{equation*}
\mathrm{HH}_{\ast }(\mathcal{A},M)\simeq \mathcal{Z}\otimes _{\mathcal{Z}%
\bigcap \mathcal{B}}\mathrm{HH}_{\ast }(\mathcal{B},M^{\mathrm{co}\mathcal{H}%
}),
\end{equation*}%
for any centrally $\mathcal{H}$-Galois extension ${\mathcal{B}\subseteq
\mathcal{A}}$ and any Hopf bimodule  $M$  which is symmetric as a $\mathcal{Z%
}$-bimodule, provided that $\mathcal{H}$ is finite-dimensional and that $%
\dim \mathcal{H}$\textrm{\ }is not zero in $\,\mathbb{K}$.  Under the same
assumptions, we also show that $\mathrm{HC}_{\ast }(\mathcal{A})^{\mathrm{co}%
\mathcal{H}}$ and $\mathrm{HC}_{\ast }(\mathcal{B})$ are isomorphic, cf.
Theorem \ref{izo3}. We conclude the paper by indicating a method to produce
examples of centrally Hopf-Galois extensions of non-commutative algebras.

\section{Preliminaries}

In order to state and prove our main result we need several basic facts
concerning Hochschild homology of Hopf-Galois extensions. Those that are
well-known will be only stated, for details the reader being referred to
\cite{JS,SS,St1,St2}.

\begin{noname}
\label{Hopf_module}Let $\mathcal{H}$ be a Hopf algebra with comultiplication
$\Delta _{\mathcal{H}}$ and counit $\varepsilon _{\mathcal{H}}.$ To denote
the element $\Delta _{\mathcal{H}}(h)$ we shall use the $\Sigma$-notation
\begin{equation*}
\textstyle \Delta _{\mathcal{H}}(h)={\sum }h_{(1)}\otimes h_{(2)}.
\end{equation*}%
Similarly, for a left $\mathcal{H}$-comodule $(N,\rho _{N})$ and a right $%
\mathcal{H}$-comodule $(M,\rho _{M})$ we shall write
\begin{equation*}
\textstyle \rho_N(n):={\sum }n_{\left\langle -1\right\rangle }\otimes
n_{\left\langle 0\right\rangle }\qquad\text{and}\qquad \rho _{M}(m)={\sum }%
m_{\left\langle 0\right\rangle }\otimes m_{\left\langle 1\right\rangle }.
\end{equation*}
For $(M,\rho _{M})$ as above we define the set of coinvariant elements in $M$
by
\begin{equation*}
M^{\mathrm{co}\mathcal{H}}:=\{m\in M\mid \rho _{M}(m)=m\otimes 1\}.
\end{equation*}
Recall that a comodule algebra is an algebra $\mathcal{A}$ which is a right $%
\mathcal{H}$-comodule via a morphism of \emph{algebras} $\rho _{\mathcal{A}}:%
\mathcal{A}\rightarrow \mathcal{A}\otimes \mathcal{H}$. Equivalently, $(%
\mathcal{A},\rho _{\mathcal{A}})$ is an $\mathcal{H}$-comodule algebra if and only if $\rho _{\mathcal{A}}(1)=1\otimes 1$ and
\begin{equation*}
\textstyle \rho _{A}(ab)={\sum }a_{\left\langle 0\right\rangle
}b_{\left\langle 0\right\rangle }\otimes a_{\left\langle 1\right\rangle
}b_{\left\langle 1\right\rangle },\;
\end{equation*}%
for any $a,$ $b$ in $\mathcal{A}.$ The set $\mathcal{A}^{\mathrm{co}\mathcal{%
H}}$ is a subalgebra in $\mathcal{A}$. If there is no danger of confusion we
shall also denote this subalgebra by $\mathcal{B}$ and we shall say that ${%
\mathcal{B} \subseteq \mathcal{A}}$ is an $\mathcal{H}$-comodule algebra.

For an $\mathcal{H}$-comodule algebra $\mathcal{A}$, the category $\mathfrak{%
M}_{\mathcal{A}}^{\mathcal{H}}$ of right Hopf modules is defined as follows.
An object in $\mathfrak{M}_{\mathcal{A}}^{\mathcal{H}}$ is a right $\mathcal{%
A}$-module $M$ together with a right $\mathcal{H}$-coaction $\rho
_{M}:M\rightarrow M\otimes \mathcal{H}$ such that, for any $m\in M\;$and $%
a\in \mathcal{A}$, the following compatibility relation is verified
\begin{equation}
\textstyle \rho _{M}(ma)={\sum }m_{\left\langle 0\right\rangle
}a_{\left\langle 0\right\rangle }\otimes m_{\left\langle 1\right\rangle
}a_{\left\langle 1\right\rangle }.  \label{Hopf-mod}
\end{equation}
Obviously, a morphism in $\mathfrak{M}_{\mathcal{A}}^{\mathcal{H}}$ is a map
which is both $\mathcal{A}$-linear and $\mathcal{H}$-colinear.

The category $_{\mathcal{A}}\mathfrak{M}^{\mathcal{H}}$ is defined
similarly. A left $\mathcal{A}$-module and right $\mathcal{H}$-comodule $%
(M,\rho _{M})$ is a left Hopf module if, for any $m\in M$ and $a\in \mathcal{%
A}$,
\begin{equation}
\textstyle \rho _{M}(am)={\sum }a_{\left\langle 0\right\rangle
}m_{\left\langle 0\right\rangle }\otimes a_{\left\langle 1\right\rangle
}m_{\left\langle 1\right\rangle }.  \label{Hopf-modl}
\end{equation}%
By definition, a Hopf bimodule is an $\mathcal{A}$-bimodule $M$ together
with a right $\mathcal{H}$-coaction $\rho _{M}$ such that relations (\ref%
{Hopf-mod}) and (\ref{Hopf-modl}) are satisfied for all $m\in M$
and $a\in \mathcal{A}$. A morphism between two Hopf bimodules is,
by definition, a map of $\mathcal{A}$-bimodules and
$\mathcal{H}$-comodules. The category of Hopf
bimodules will be denoted by$_{\mathcal{A}}\mathfrak{M}_{\mathcal{A}}^{%
\mathcal{H}}$. For example, $\mathcal{A}$ is a Hopf bimodule.
\end{noname}

\begin{noname}
\label{ff_Galois-ext}Let ${\mathcal{B}\subseteq \mathcal{A}}$ be an $%
\mathcal{H}$-comodule algebra. Recall that ${\mathcal{B}\subseteq \mathcal{A}%
}$ is an $\mathcal{H}$-Galois extension if the canonical $\mathbb{K}$-linear
map
\begin{equation*}
\textstyle \beta :\mathcal{A}\otimes _{\mathcal{B}}\mathcal{A}\rightarrow
\mathcal{A}\otimes \mathcal{H},\;\beta (a\otimes x)={\sum }ax_{\left\langle
0\right\rangle }\otimes x_{\left\langle 1\right\rangle }
\end{equation*}%
is bijective. Note that $\mathcal{A}\otimes _{\mathcal{B}}\mathcal{A}$ is an
object in $_{\mathcal{A}}\mathfrak{M}_{\mathcal{A}}^{\mathcal{H}}$ with
respect to the canonical bimodule structure and the $\mathcal{H}$-coaction
defined by $\mathcal{A}\otimes _{\mathcal{B}}\rho _{\mathcal{A}}.$ One can
also regard $\mathcal{A}\otimes \mathcal{H}$ as an object in $_{\mathcal{A}}%
\mathfrak{M}_{\mathcal{A}}^{\mathcal{H}}$ with the $\mathcal{A}$-bimodule
structure
\begin{equation*}
\textstyle a\cdot (x\otimes h)\cdot a^{\prime}=\sum axa_{\left\langle
0\right\rangle }^{\prime}\otimes ha_{\left\langle 1\right\rangle }^{\prime},
\end{equation*}%
and the $\mathcal{H}$-coaction defined by $\mathcal{A}\otimes \Delta _{%
\mathcal{H}}.$ With respect to these Hopf bimodule structures $\beta $ is a
morphism in $_{\mathcal{A}}\mathfrak{M}_{\mathcal{A}}^{\mathcal{H}}.$
\end{noname}

\begin{definition}
\label{nn:X_R}For a $\mathbb{K}$-algebra $\mathcal{R}$ and an $\mathcal{R}$%
-bimodule $X$ we define
\begin{equation*}
X_{\mathcal{R}}:=X/[\mathcal{R},X]\qquad\text{and}\qquad X^{\mathcal{R}%
}:=\{x\in X\mid rx=xr,\forall r\in\mathcal{R}\},
\end{equation*}%
where $[\mathcal{R},X]$ is the $\mathbb{K}$-subspace of $X$ generated by all
commutators $rx-xr,$ with $r\in \mathcal{R}$ and $x\in X$. The class of $%
x\in X$ in $X_{\mathcal{R}}$ will be denoted by $[x]_{\mathcal{R}}.$
\end{definition}

\begin{remark}
If $X$ is an $\mathcal{R}$-bimodule then $X_{\mathcal{R}}\cong \mathcal{R}%
\otimes _{\mathcal{R}^{e}}X$, where $\mathcal{R}^{e}:=\mathcal{R}\otimes
\mathcal{R}^{op}$ denotes the enveloping algebra of $\mathcal{R}$. The
isomorphism is given by $[x]_{\mathcal{R}}\mapsto 1\otimes_{\mathcal{R}^e}x$.
\end{remark}

\begin{noname}
\label{pa:Galois} Let now ${\mathcal{B}\subseteq \mathcal{A}}$ be
an arbitrary extension of algebras. By \cite[p. 145]{JS} it
follows that $(\mathcal{A} \otimes
_{\mathcal{B}}\mathcal{A})^{\mathcal{B}}$ is an associative
algebra with the multiplication given by
\begin{equation}
\textstyle zz^{\prime }=\sum_{i=1}^{n}\sum_{j=1}^{m}a_{i}a_{j}^{\prime
}\otimes _{\mathcal{B}}b_{j}^{\prime }b_{i},  \label{ec:enveloping}
\end{equation}%
where $z=\sum_{i=1}^{n}a_{i}\otimes _{\mathcal{B}}b_{i}${\ and }$z^{\prime
}=\sum_{j=1}^{m}a_{j}^{\prime }\otimes _{\mathcal{B}}b_{j}^{\prime }$ are
arbitrary elements in $(\mathcal{A}\otimes _{\mathcal{B}}\mathcal{A})^{
\mathcal{B}}.$ Moreover, if $M$ is an $\mathcal{A}$-bimodule then $M_{%
\mathcal{B}}$ is a right $(\mathcal{A}\otimes _{ \mathcal{B}}\mathcal{A})^{%
\mathcal{B}}$-module with respect to the action that, for $m$ in $M$ and $%
z=\sum_{i=1}^{n}a_{i}\otimes _{\mathcal{B}}b_{i}$ in $(\mathcal{A}\otimes _{%
\mathcal{B}}\mathcal{A})^{\mathcal{B}}$, is defined by%
\begin{equation}
\textstyle \lbrack m]_{\mathcal{B}}\cdot z=\sum_{i=1}^{n}[b_{i}ma_{i}]_{%
\mathcal{B}}.  \label{ec:M_B=modul}
\end{equation}
\end{noname}

\begin{noname}
Suppose now that ${\mathcal{B}\subseteq \mathcal{A}}$ is an $\mathcal{H}$%
-Galois extension and let $M$ be an $\mathcal{A}$-bimodule. Let $i:\mathcal{%
H}\longrightarrow {\mathcal{A}\otimes {\mathcal{H}}}$ denote the canonical
map $i(h)=1\otimes {h}$. Following \cite[p. 146]{JS} we define
\begin{equation*}
\kappa :\mathcal{H}\longrightarrow (\mathcal{A} \otimes _{\mathcal{B}}%
\mathcal{A})^{\mathcal{B}},\ \ \ \kappa :=\beta ^{-1}\circ {i.}
\end{equation*}
For $h\in {\mathcal{H}}$ we shall use the notation $\kappa (h)=\sum \kappa
^{1}(h) \otimes _{\mathcal{B}}\kappa ^{2}(h)$. Thus, by definition,
\begin{equation}
\textstyle \sum \kappa ^{1}(h)\kappa ^{2}(h)_{\langle {0}\rangle }\otimes
\kappa ^{2}(h)_{\langle {1}\rangle }=1\otimes {h}.  \label{eq:k}
\end{equation}%
By \cite[p. 146]{JS}, $\kappa $ is an anti--morphism of algebras. Hence
\begin{equation}  \label{eq:R-module}
\textstyle h\cdot \lbrack m]_{\mathcal{B}}=\sum [\kappa ^{2}(h)m\kappa
^{1}(h)]_{\mathcal{B}}
\end{equation}%
defines a left $\mathcal{H}$-action on $M_{\mathcal{B}}$. Obviously this
structure is functorial in $M,$ so we get a functor $(-)_{\mathcal{B}}:{}_{%
\mathcal{A}}\mathfrak{M}_{\mathcal{A}}\rightarrow {}_{\mathcal{H}}\mathfrak{M%
}$.
\end{noname}

\begin{noname}
Let $\mathcal{H}$ be a Hopf algebra with multiplication $m$ and
comultiplication $\Delta $. Let $\tau :\mathcal{H}\otimes \mathcal{H}%
\rightarrow \mathcal{H}\otimes \mathcal{H}$ denote the usual flip map $%
x\otimes y\mapsto y\otimes x.$ One can prove that
\begin{equation*}
\mathcal{C}_{\mathcal{H}}:=\mathrm{coker}\,(m-m\circ \tau)
\end{equation*}
is a quotient coalgebra of $\mathcal{H},$ as the linear space generated by
all commutators in $\mathcal{H}$ is a coideal of $\mathcal{H}$. The
canonical projection onto $\mathcal{C}_{\mathcal{H}}$ will be denoted by $%
\pi _{\mathcal{H}}.$ Note that $\pi _{\mathcal{H}}\ $is a trace map, that is
$\pi _{\mathcal{H}}(hk)=\pi _{\mathcal{H}}(kh)$ for all $h$ and $k\ $in $%
\mathcal{H} .$

Dually,
\begin{equation*}
\mathcal{R}_{\mathcal{H}}:=\mathrm{ker}(\Delta -\tau \circ \Delta) =\{r\in
\mathcal{H}\mid \sum r_{(1)}\otimes r_{(2)} =\sum r_{(2)}\otimes r_{(1)}\}
\end{equation*}
is a subalgebra of $\mathcal{H}$. It is clear that
\begin{equation}
\textstyle \mathcal{R}_{\mathcal{H}} =\{r\in \mathcal{H}\mid \sum
r_{(1)}\otimes r_{(2)}\otimes r_{(3)}=\sum r_{(2)}\otimes r_{(3)}\otimes
r_{(1)}\}.  \label{eq:R_H}
\end{equation}
\end{noname}

%\begin{noname}
%Let $\mathcal{B}\subseteq \mathcal{A}$ be an $\mathcal{H}$-comodule algebra.
%For an arbitrary $(\mathcal{A},\mathcal{H})$-Hopf bimodule $(M,\rho )$ we
%set:
%\begin{equation*}
%M^{\mathrm{co}{}\mathcal{H}}:=\{m\in M\mid \rho (m)=m\otimes 1_{\mathcal{H}%
%}\}.
%\end{equation*}%
%On the other hand, since any right $\mathcal{H}$-comodule can be regarded as
%a right $\mathcal{C}_{\mathcal{H}}$ comodule with the coaction $\rho
%^{\prime }:=(M\otimes \pi _{\mathcal{H}})\circ \rho $ one can define%
%\begin{equation*}
%M^{\mathrm{co}{}\mathcal{C}_{\mathcal{H}}}:=\{m\in M\mid \rho ^{\prime
%}(m)=m\otimes \pi _{\mathcal{H}}(1_{\mathcal{H}})\}.
%\end{equation*}%
%Obviously $M^{\mathrm{co}{}\mathcal{H}}\subseteq M^{\mathrm{co}{}\mathcal{C}%
%_{\mathcal{H}}}$. Note that, in general, this inclusion is strict.
%\end{noname}

\begin{noname}
\label{nn:cotensor}If $\mathcal{C}$ is a coalgebra and $\mathcal{R}$ is an
algebra we define the category $_{\mathcal{R}}\mathfrak{M}^{\mathcal{C}}$ as
follows. The objects in $_{\mathcal{R}}\mathfrak{M}^{\mathcal{C}}$ are left $%
\mathcal{R}$-modules and right $\mathcal{C}$-comodules such that the map $%
\rho_M$ that defines the coaction on $M$ is $\mathcal{R}$-linear,
that is
\begin{equation*}
\textstyle \rho_M (rm)=\sum rm_{\left\langle 0\right\rangle }\otimes
m_{\left\langle 1\right\rangle }.
\end{equation*}%
A map $f:M\rightarrow N$ is a morphism in $_{\mathcal{R}}\mathfrak{M}^{%
\mathcal{C}}$ if it is $\mathcal{R}$-linear and $\mathcal{C}$-colinear.

For a right $\mathcal{C}$-comodule $(M,\rho_M)$ and a left $\mathcal{C}$%
-comodule $(V,\rho_V)$ we define their cotensor product by%
\begin{equation*}
M\square _{\mathcal{C}}V:=\ker \left( \rho _{M}\otimes N-M\otimes \rho
_{V}\right) ,
\end{equation*}
Recall that $V$ is said to be \emph{coflat} if the functor $(-)\square _{%
\mathcal{C}}V : \mathcal{M}^{\mathcal{C}}\to{}_\mathbb{K}\mathcal{M}$ is
exact. By \cite[Theorem 2.4.17]{DNR} $V$ is coflat if and only if $V$ is
an injective object in the category of left $\mathcal{C}$-comodules.

Note that, if $M\in {}_{\mathcal{R}}\mathfrak{M}^{\mathcal{C}}$ and $V\in
{}^{\mathcal{C}}\mathfrak{M,}$ then $M\square _{\mathcal{C}}V$ is an $%
\mathcal{R}$-submodule of $M\otimes V,$ as $M\square _{\mathcal{C}}V$ is the
kernel of an $\mathcal{R}$-linear map. Dually, for a right $\mathcal{R}$%
-module $X$ and an object $M$ in $_{\mathcal{R}}\mathfrak{M}^{\mathcal{C}},$
the tensor product $X\otimes _{\mathcal{R}}M$ is a quotient $\mathcal{C}$%
-comodule of $X\otimes M.$

In some special cases the cotensor product and the tensor product
\textquotedblleft commute\textquotedblright . For instance, if $X$ is a
right $\mathcal{R}$-module, $V$ is a left $\mathcal{C}$-comodule and $M\in
{}_{\mathcal{R}}\mathfrak{M}^{\mathcal{C}}$ then
\begin{equation}
(X\otimes _{\mathcal{R}}M)\square _{\mathcal{C}}V\cong X\otimes _{\mathcal{R}%
}(M\square _{\mathcal{C}}V),  \label{eq:tens-cotens}
\end{equation}%
provided that either $X$ is flat or $V$ is injective.
\end{noname}

\begin{noname}
\label{delta_functor}Let $\mathfrak{A}$ and $\mathfrak{B}$ be two abelian
categories and assume that, for each $n\in \mathbb{N},$ a functor $T_{n}:%
\mathfrak{A\longrightarrow B\ \ }$is given. We say that $T_{\ast }$ is a
\emph{homological }$\delta $\emph{-functor} if, for every short exact
sequence
\begin{equation}
0\longrightarrow X^{\prime }\longrightarrow X\longrightarrow X^{\prime
\prime }\longrightarrow 0  \label{eq:short_exact}
\end{equation}%
in $\mathfrak{A}$ and $n>0$, there are ``connecting'' morphism $\delta
_{n}:T_{n}(X^{\prime \prime })\rightarrow T_{n-1}(X^{\prime })$ such that
\begin{equation*}
\cdots \longrightarrow T_{n}(X^{\prime })\longrightarrow
T_{n}(X)\longrightarrow T_{n}(X^{\prime \prime })\overset{\delta _{n}}{%
\longrightarrow }T_{n-1}(X^{\prime })\longrightarrow \cdots
\end{equation*}%
is exact and functorial in the sequence in (\ref{eq:short_exact}).
Furthermore, $T_{\ast }$ is said to be \emph{effaceable} if, for each object
$X$ in $\mathfrak{A,}$ there is an object $P$ in $\mathfrak{A}$ together
with an epimorphism from $P$ to $X$ such that $T_{n}(P)=0$ for any $n>0.$ A
morphism of $\delta $-functors $T_{\ast }$ and $S_{\ast }$ with connecting
homomorphisms $\delta _{\ast }$ and $\partial _{\ast }$, respectively, is a
sequence of natural transformations $\phi _{\ast }:T_{\ast }\rightarrow
S_{\ast }$ such that, for $n>0,$%
\begin{equation*}
\phi _{n-1}\circ \delta _{n}=\partial _{n}\circ \phi _{n}.
\end{equation*}
By Theorem 7.5 in \cite[Chapter III]{Br}, homological and effaceable $\delta
$-functors have the following universal property. If $T_{\ast }$ and $%
S_{\ast }$ are homological and effaceable $\delta$-functors and $\phi
_{0}:T_{0}\rightarrow S_{0}$ is a natural transformation, then there is a
unique morphism of $\delta$-functors $\phi _{\ast }:T_{\ast }\rightarrow
S_{\ast }$ that lifts $\phi _{0}$.
\end{noname}

\begin{proposition}
\label{le:PropHopfBimod}Let $\mathcal{H}$ be a Hopf algebra with antipode $%
S_{\mathcal{H}}$. Suppose that ${\mathcal{B}\subseteq \mathcal{A}}$ is an $%
\mathcal{H}$-Galois extension, $(M,\rho_M )$ is a Hopf bimodule and $V$ is a
left $\mathcal{H}$-comodule.

\begin{enumerate}
\item There is a $\mathbb{K}$-linear map $\rho _0(M):M_{\mathcal{B}%
}\rightarrow M_{\mathcal{B}}\otimes \mathcal{H}$ such that $(M_{\mathcal{B}%
},\rho _0(M))$ is a quotient $\mathcal{H}$-comodule of $(M,\rho_M).$ Moreover, for $%
m\in M$ and $h\in \mathcal{H}$ we have
\begin{equation}
\textstyle \rho _0(M)(h\cdot \lbrack m]_{\mathcal{B}})=\sum h_{(2)}\cdot
\lbrack m_{\left\langle 0\right\rangle }]_{\mathcal{B}}\otimes
h_{(3)}m_{\left\langle 1\right\rangle }S_{\mathcal{H}}h_{(1)}.
\label{ec:M_B=CrossedModule}
\end{equation}

\item If $S_\mathcal{H}$ is an involution, then $(M_{\mathcal{B}}\otimes \pi
_{\mathcal{H}})\circ \rho _0(M):M_{\mathcal{B}}\rightarrow M_{\mathcal{B}%
}\otimes \mathcal{C}_{\mathcal{H}}$ is a morphism of left $\mathcal{R}_{%
\mathcal{H}}$-modules. Hence, with respect to the above $\mathcal{C}_{%
\mathcal{H}}$-comodule structure, $M_{\mathcal{B}}\ $is an object in $_{%
\mathcal{R}_{\mathcal{H}}}\mathfrak{M}^{\mathcal{C}_{\mathcal{H}}}$ and, for
any right $\mathcal{R}_\mathcal{H}$-module $X$, the coalgebra $\mathcal{C}_{%
\mathcal{H}}$ coacts canonically on $X\otimes _{\mathcal{R}_\mathcal{H}}M_{%
\mathcal{B}}$. The $\mathcal{C}_{\mathcal{H}}$-coaction on $M_{\mathcal{B}}$
will be denoted by $\rho _0(M)$ too.

\item If $V$ is a left $\mathcal{C}_{\mathcal{H}}$-comodule then $M\square _{%
\mathcal{C}_{\mathcal{H}}}V$ is a $\mathcal{B}^{e}$-submodule of $M\otimes V$%
. Under the additional assumption that $V$ is injective, $(M\square _{%
\mathcal{C}_{\mathcal{H}}}V)_{\mathcal{B}}$ and $M_{\mathcal{B}}\square _{%
\mathcal{C}_{\mathcal{H}}}V$ are isomorphic linear spaces. In particular,
the action of $\mathcal{R}_{\mathcal{H}}$ on the latter vector space can be
transported to $(M\square _{\mathcal{C}_{\mathcal{H}}}V)_{\mathcal{B}}$.
\end{enumerate}
\end{proposition}

\begin{proof}
(1) Obviously, $(M_{\mathcal{B}},\rho _0(M))$ is a quotient $\mathcal{H}$%
-comodule of $M$, as $[\mathcal{B},M]$ is a subcomodule of $M$. For $M=%
\mathcal{A},$ identity (\ref{ec:M_B=CrossedModule}) is proven in
\cite[ Proposition 2.6]{JS}. The general case can be handled in a
similar manner,
replacing $a\in \mathcal{A}$ by $m\in M$ everywhere in the proof of \cite[%
Relation (6)]{JS}.

(2) Recall that $S_\mathcal{H}$ is an involution, i.e. $S_{\mathcal{H}}^{2}=%
\mathrm{Id}_{\mathcal{H}},$ if and only if
\begin{equation}
\textstyle \sum r_{(2)}S_{\mathcal{H}}r_{(1)}=\sum S_{\mathcal{H}%
}r_{(2)}r_{(1)}=\varepsilon (r)1_{\mathcal{H}}.  \label{ec:involutie}
\end{equation}
Clearly $(M_{\mathcal{B}}\otimes \pi_\mathcal{H})\circ\rho_0(M) $ defines a $%
\mathcal{C}_{\mathcal{H}}$-comodule structure on $M_{\mathcal{B}}$. For
brevity we shall denote this map by $\rho_0(M)$ too. For $r\in \mathcal{R}_{%
\mathcal{H}}$ and $m\in M$ we get
\begin{align*}
\rho _0(M)(r\cdot \lbrack m]_0(M))& =\textstyle\sum r_{(2)}\cdot \lbrack
m_{\left\langle 0\right\rangle }]_{\mathcal{B}}\otimes \pi _{\mathcal{H}%
}(r_{(3)}m_{\left\langle 1\right\rangle }S_{\mathcal{H}}r_{(1)}) \\
& =\textstyle\sum r_{(2)}\cdot \lbrack m_{\left\langle 0\right\rangle }]_{%
\mathcal{B}}\otimes \pi _{\mathcal{H}}(m_{\left\langle 1\right\rangle }S_{%
\mathcal{H}}r_{(1)}r_{(3)}) \\
& =\textstyle\sum r_{(3)}\cdot \lbrack m_{\left\langle 0\right\rangle }]_{%
\mathcal{B}}\otimes \pi _{\mathcal{H}}(m_{\left\langle 1\right\rangle
}r_{(2)}S_{\mathcal{H}}r_{(1)}) \\
& =\textstyle\sum r\cdot \lbrack m_{\left\langle 0\right\rangle }]_{\mathcal{%
B}}\otimes \pi _{\mathcal{H}}(m_{\left\langle 1\right\rangle }) \\
& =r\cdot \rho _0(M)([m]_{\mathcal{B}}).
\end{align*}%
Note that the second and the third equalities are consequences of the fact
that $\pi _{\mathcal{H}}$ is a trace map and respectively of relation (\ref%
{eq:R_H}). To deduce the penultimate identity we use
(\ref{ec:involutie}).
In conclusion, $\rho _0(M)$ is a morphism of $\mathcal{R}_{\mathcal{H}}$%
-modules. Hence $M$ is an object in $_{\mathcal{R}_{\mathcal{H}}}\mathfrak{M}%
^{\mathcal{C}_{\mathcal{H}}}\ $and, in view of \S \ref{nn:cotensor}, one can
regard $M_{\mathcal{B}}\otimes _{\mathcal{R}}X$ as a quotient comodule of $%
M_{\mathcal{B}}\otimes X.$

(3) Obviously, $\rho^{\prime }:=(M\otimes \pi _{\mathcal{H}})\circ \rho$
defines a $\mathcal{C}_{\mathcal{H}}$-coaction on $M$ and it is a morphism
of $\mathcal{B}$-bimodules, as $\rho$ is so. Thus $(M,\rho ^{\prime })$ is
an object in $_{B^{e}}\mathfrak{M}^{\mathcal{C}_{\mathcal{H}}}$ and $%
M\square _{\mathcal{C}_{\mathcal{H}}}V$is a $\mathcal{B}^{e}$-submodule of $%
M\otimes V$, cf. \S \ref{nn:cotensor}. If $V$ is an injective $\mathcal{C}_{%
\mathcal{H}}$-comodule, then
\begin{equation}
(M\square _{\mathcal{C}_{\mathcal{H}}}V)_{\mathcal{B}}\cong B\otimes _{%
\mathcal{B}^{e}}(M\square _{\mathcal{C}_{\mathcal{H}}}V)\cong (B\otimes _{%
\mathcal{B}^{e}}M)\square _{\mathcal{C}_{\mathcal{H}}}V\cong M_{\mathcal{B}%
}\square _{\mathcal{C}_{\mathcal{H}}}V.  \label{eq:iso_B}
\end{equation}%
Since $M_{\mathcal{B}}$ is an object in $_{\mathcal{R}_{\mathcal{H}}}%
\mathfrak{M}^{\mathcal{C}_{\mathcal{H}}},$ it follows that $M_{\mathcal{B}%
}\square _{\mathcal{C}_{\mathcal{H}}}V$ is an $\mathcal{R}_{\mathcal{H}}$%
-submodule of $M_{\mathcal{B}}\otimes V.$ In particular, $M_{\mathcal{B}%
}\square _{\mathcal{C}_{\mathcal{H}}}V$ is an $\mathcal{R}_{\mathcal{H}}$%
-module$.$ To conclude the proof, we take on $(M\square _{\mathcal{C}_{%
\mathcal{H}}}V)_{\mathcal{B}}$ the unique $\mathcal{R}_{\mathcal{H}}$-action
that makes the composition of the isomorphisms in (\ref{eq:iso_B}) an $%
\mathcal{R}_{\mathcal{H}}$-linear map.
\end{proof}

\begin{remark}
We keep the assumptions in the third part of Proposition \ref%
{le:PropHopfBimod}. Let $z:=\sum_{i=1}^{n}m_{i}\otimes v_{i}$ be an element
in $M\square _{\mathcal{C}_{\mathcal{H}}}V$. The composition of the $%
\mathcal{R}_{\mathcal{H}}$-linear isomorphisms in (\ref{eq:iso_B}) maps $%
[z]_{\mathcal{B}}\ $ to $\sum_{i=1}^{n}[m_{i}]_{\mathcal{B}}\otimes v_{i}.$
Obviously, this isomorphism is natural in $M\in {}_{\mathcal{A}}\mathfrak{M}%
_{\mathcal{A}}^{\mathcal{H}}.$ It is not hard to see that, for $h\in\mathcal{%
R}_\mathcal{H}$,
\begin{equation*}
h\cdot [z]_{\mathcal{B}}=\textstyle\sum_{i=1}^{n}[\kappa ^{2}(h)m_i\kappa
^{1}(h)\otimes v_i]_{\mathcal{B}}.
\end{equation*}
\end{remark}

\section{The spectral sequence}

In this section, given an $\mathcal{H}$-Galois extension
$\mathcal{B}\subseteq
\mathcal{A}$, a Hopf bimodule $M$ and an injective left $\mathcal{C}_{%
\mathcal{H}}$-comodule $V$, we construct a spectral sequence that converges
to $\mathrm{H}_{\ast }(\mathcal{A},M)\square _{\mathcal{C}_{\mathcal{H}}}V$.
Our result, Theorem \ref{te:sir}, will be obtained as a direct application
of a variant of Grothendieck's spectral sequence, which will be deduced from
the following two lemmas and \cite[Corollary 5.8.4]{We}.

Recall that a category $\mathfrak{A}$ is cocomplete if and only if
any set of objects in $\mathfrak{A}$ has a direct sum. If $X$ is
an object in a category, then we shall also write $X$ for the
identity map of $X$.

\begin{lemma}
\label{le:phi=izo}Let $\mathfrak{A}$ and $\mathfrak{B}$ be
cocomplete abelian categories and let $H,\;H^{\prime
}:\mathfrak{A}\rightarrow \mathfrak{B}$ be two right exact
functors that commute with direct sums. If $U$ is a generator in
$\mathfrak{A}$ and there is a natural morphism $\phi
:H\rightarrow H^{\prime }$ such that $\phi (U)$ is an isomorphism, then $%
\phi (X)$ is an isomorphism, for every object $X$ in $\mathfrak{A.}$
\end{lemma}

\begin{proof}
Let $X$ be an object in $\mathfrak{A}.$ Since $U$ is a generator in $%
\mathfrak{A},$ there is an exact sequence%
\begin{equation*}
U^{(J)}\overset{u}{\longrightarrow }U^{(I)}\overset{v}{\longrightarrow }%
X\longrightarrow 0,
\end{equation*}%
where $I$ and $J$ are certain sets. Hence in the following diagram%
\begin{equation*}
\xymatrix{ H(U^{(J)}) \ar[d]_{\phi (U^{(J)})} \ar[r]^{H(u)} & H(U^{(I)})
\ar[d]_{\phi (U^{(I)})}\ar[r]^{H(v)} & H(X)\ar[d]^{\phi (X)} \ar[r] & 0\\
H'(U^{(J)})\ar[r]_{H'(u)} & H'(U^{(I)})\ar[r]_{H'(v)} & H'(X)\ar[r]&0}
\end{equation*}
the squares are commutative and the lines are exact. Recall that $H$
commutes with direct sums if the canonical map $\alpha :\oplus _{i\in
I}H(X_{i})\rightarrow H(\oplus _{i\in I}X_{i})$ is an isomorphism for each
family of objects $(X_{i})_{i\in I}$ in $\mathfrak{A}$. Now one can see
easily that $\phi (U^{(I)})$ and $\phi (U^{(J)})$ are isomorphisms, as $H$
and $H^{\prime }$ commute with direct sums and $\phi (U)$ is an isomorphism.
Thus $\phi (X)$ is an isomorphism too.
\end{proof}

\begin{lemma}
\label{le:L_nF_1(F_2(A)=0}Let $\mathfrak{A},$ $\mathfrak{B}$ and $\mathfrak{%
\mathcal{C}}$ be cocomplete abelian categories with enough projective objects. Let $%
F:\mathfrak{\mathcal{B}}\rightarrow \mathfrak{C}$ and $G:\mathfrak{A}%
\rightarrow \mathfrak{B}$ be right exact functors that commute with direct
sums. If $U$ is a generator in $\mathfrak{A}$ such that $G(U)$ is $F$%
-acyclic, then $G(P)$ is $F$-acyclic for any projective object $P$ in $%
\mathfrak{A}.$
\end{lemma}

\begin{proof}
Recall that $G(U)$ is $F$-acyclic if $\mathrm{L}_{n}F(G(U))=0$ for any $%
n>0. $ Let $P$ be a projective object in $\mathfrak{A.}$ There is a set $I$
such that $P$ is a direct summand of $U^{(I)}.$ Hence $G(P)$ is a direct
summand of $G(U^{(I)}).$ On the other hand, the proof of \cite[Corollary
2.6.11]{We} works for any functor that commutes with direct sums. Thus $%
\mathrm{L}_{n}F$ commutes with direct sums, so $G(U^{(I)})\cong G(U)^{(I)}$
is $F$-acyclic. Then $G(P)$ is also $F$-acyclic.
\end{proof}

\begin{proposition}
\label{pr:sir spectral}Let $G:\mathfrak{A}\rightarrow \mathfrak{B}\mathfrak{,%
}$ $F:\mathfrak{B}\rightarrow \mathfrak{C}$ and $H:\mathfrak{A}\rightarrow
\mathfrak{C}$ be right exact functors that commute with direct sums, where $%
\mathfrak{A},$ $\mathfrak{B}$ and $\mathfrak{\mathcal{C}}$ are
cocomplete abelian categories with enough projective objects.
Assume that $U$ is a generator in $\mathfrak{A}$ and that $\phi
:F\circ G\rightarrow H$ is a natural transformation. If $\phi (U)$
is an isomorphism and $G(U)$ is $F$-acyclic then, for every object
$X$ in $\mathfrak{A},$ there is a functorial spectral sequence
\begin{equation}
\mathrm{L}_{p}F(\mathrm{L}_{q}G(X))\Longrightarrow \mathrm{L}_{p+q}H(X).
\label{ec:SirSpectralGeneral}
\end{equation}
\end{proposition}

\begin{proof}
By Lemma \ref{le:phi=izo}, for every object $X$ the morphism $\phi (X)$ is
an isomorphism. On the other hand, if $X$ is projective in $\mathfrak{A}$
then $\mathrm{L}_{n}F(G(X))=0$ for every $n\in \mathbb{N}^{\ast }.$ Hence,
we obtain (\ref{ec:SirSpectralGeneral}) as a particular case of \cite[%
Corollary 5.8.4]{We}.
\end{proof}

We take ${\mathcal{B}\subseteq \mathcal{A}}$ to be a faithfully flat Galois
extension. For proving Theorem \ref{te:sir}, one of our main results, we
shall apply Proposition \ref{pr:sir spectral}. In order to do that we need
some properties of the category $_\mathcal{A}\mathcal{M}_\mathcal{A}^%
\mathcal{H}$. We start with the following.

\begin{proposition}
\label{pr:progenerator}Let $\mathcal{H}$ be a Hopf algebra with bijective
antipode. Let ${\mathcal{B}\subseteq \mathcal{A}}$ be a faithfully flat $%
\mathcal{H}$-Galois extension. Then $\mathcal{A}\otimes \mathcal{A}$ is a
projective generator in the category of Hopf bimodules. It is also
projective as a $\mathcal{B}$-bimodule.
\end{proposition}

\begin{proof}
By \cite[Theorem 4.10]{SS} the induction functor $(-)\otimes _{\mathcal{B}}%
\mathcal{A}:{}\mathfrak{M}_{\mathcal{B}}{}\rightarrow \mathfrak{M}_{\mathcal{%
A}}^{\mathcal{H}}$ is an equivalence of categories and its inverse is $(-)^{%
\mathrm{co}\mathcal{H}}.$ We deduce that, for an arbitrary right Hopf module
$X$, the canonical map $X^{\mathrm{co}{}\mathcal{H}}\otimes _{\mathcal{B}}%
\mathcal{A}\rightarrow X$ induced by the module structure of $X$ is an
isomorphism of right Hopf modules. Let $M$ be a Hopf bimodule. Hence, the $%
\mathcal{A}$-bimodule structure on $M$ defines an epimorphism $\mathcal{A}%
\otimes M^{\mathrm{co}{}\mathcal{H}}\otimes \mathcal{A}\rightarrow M$ of
Hopf bimodules. Thus $\mathcal{A}\otimes \mathcal{A}$ is a generator in the
category ${}_{\mathcal{A}}\mathfrak{M}_{\mathcal{A}}^{\mathcal{H}}$.

Let $p:X\rightarrow Y$ be an epimorphism of Hopf bimodules and  $f:%
\mathcal{A}\otimes \mathcal{A}\rightarrow Y$ be an arbitrary morphism in $%
{}_{\mathcal{A}}\mathfrak{M}_{\mathcal{A}}^{\mathcal{H}}$. We want to show
that there is a morphism $g:\mathcal{A}\otimes \mathcal{A}\rightarrow X$ of
Hopf bimodules such that $p\circ g=f.$ Indeed, if $y:=f(1_{\mathcal{A}%
}\otimes 1_{\mathcal{A}})$ then $\emph{y}\in Y^{\mathrm{co}\mathcal{H}}.$
Since $(-)^{\mathrm{co}\mathcal{H}}:\mathfrak{M}_{\mathcal{A}}^{\mathcal{H}%
}\rightarrow \mathfrak{M}_{\mathcal{B}}$ is an equivalence of categories it
follows that $(-)^{\mathrm{co}\mathcal{H}}$ is exact. Hence $p(X^{\mathrm{co}%
\mathcal{H}})=$ $Y^{\mathrm{co}\mathcal{H}}.$ Let $x\in $ $X^{\mathrm{co}%
\mathcal{H}}$ be an element such that $p(x)=y.$ There is a unique morphism
of $\mathcal{A}$-bimodules $g:\mathcal{A}\otimes \mathcal{A}\rightarrow X$
such that $g(a^{\prime }\otimes a^{\prime \prime })=a^{\prime }xa^{\prime
\prime }.$ Since $x$ is an $\mathcal{H}$-coinvariant element in $X,$ one can
check easily that $g$ is a map of $\mathcal{H}$-comodules too. Obviously, $%
p\circ g=f.$

By \cite[Theorems 4.9 and 4.10]{SS} $\mathcal{A}$ is projective as a left
and right $\mathcal{B}$-module. Thus $\mathcal{A}^{e}$ is projective as a
left $\mathcal{B}^{e}$-module, that is $\mathcal{A}\otimes \mathcal{A}$ is a
projective $\mathcal{B}$-bimodule.
\end{proof}

\begin{corollary}
Let $\mathcal{H}$ be a Hopf algebra with bijective antipode. If ${\mathcal{B}%
\subseteq \mathcal{A}}$ is a faithfully flat $\mathcal{H}$-Galois extension
then ${}_{\mathcal{A}}\mathfrak{M}_{\mathcal{A}}^{\mathcal{H}}$ has enough
projective objects.
\end{corollary}

\begin{proof}
Every category with a projective generator has enough projective objects.
\end{proof}

\begin{definition}
\label{nn:action of the center}For a $\mathbb{K}$-algebra $\mathcal{R}$ and
an $\mathcal{R}$-bimodule $X,$ let $\left( C_{\ast }(\mathcal{R},X),b_{\ast
}\right) $ be the chain complex given by $C_{n}(\mathcal{R},X)=X\otimes
\mathcal{R}^{\otimes n}$ and
\begin{align*}
b_{n}(x\otimes r_{1}\otimes \cdots \otimes r_{n}) =xr_{1}\otimes \cdots
\otimes r_{n}&+\textstyle\sum_{i=1}^{n-1}(-1)^{i}x\otimes r_{1}\otimes
\cdots \otimes r_{i}r_{i+1}\otimes \cdots \otimes r_{n} \\
&+(-1)^{n}r_{n}x\otimes r_{1}\otimes \cdots \otimes r_{n-1}.
\end{align*}%
\emph{Hochschild homology} of $\mathcal{R}$ with coefficients in $X$ is, by
definition, the homology of $\left( C_{\ast }(\mathcal{R},X),b_{\ast
}\right) $. It will be denoted by $\mathrm{HH}{}_{\ast }(\mathcal{R},X).$
\end{definition}

\begin{noname}
Let $\mathcal{R}$ and $X$ be as in the above definition. It is well-known
that Hochschild homology of $\mathcal{R}$ with coefficients in $X$ may be
defined in an equivalent way by
\begin{equation*}
\mathrm{HH}{}_{\ast }(\mathcal{R},X)=\mathrm{Tor}{}_{\ast }^{\mathcal{R}%
^{e}}(\mathcal{R},X).
\end{equation*}%
Since $X_{\mathcal{R}}\cong \mathcal{R}\otimes _{\mathcal{R}^{e}}X$, it also
follows that $\mathrm{HH}{}_{\ast }(\mathcal{R},-)$ are the left derived
functors of $(-)_{\mathcal{R}}:{}_{\mathcal{R}}\mathfrak{M}_{\mathcal{R}%
}\rightarrow \mathfrak{M}_{\mathbb{K}}.$
\end{noname}

\begin{noname}
\label{nn:H-comodule}Let ${\mathcal{B}\subseteq \mathcal{A}}$ be an $%
\mathcal{H}$-comodule algebra and  $M$ be a Hopf bimodule. Obviously, $%
\rho _{n}(M):C{}_{n}(\mathcal{B},M)\rightarrow C{}_{n}(\mathcal{B},M)\otimes
\mathcal{H}$ given by
\begin{equation*}
\textstyle\rho _{n}(M)\left( m\otimes b_{1}\otimes \cdots \otimes
b_{n}\right) =\sum \left( m_{\left\langle 0\right\rangle }\otimes
b_{1}\otimes \cdots \otimes b_{n}\right) \otimes m_{\left\langle
1\right\rangle }.
\end{equation*}%
defines a comodule structure on $C{}_{n}(\mathcal{B},M)$ such that $C_{\ast
}(\mathcal{B},M)$ is a complex of right $\mathcal{H}$-comodules. Note that,
if $\mathcal{Z}$ is the center of $\mathcal{A}$ then $C_{\ast }(\mathcal{B}%
,M)$ is a complex of left $\mathcal{Z}_{0}$-modules, where $\mathcal{Z}_{0}:=%
\mathcal{Z}\bigcap \mathcal{B}$. Indeed, $\mathcal{Z}_{0}$-acts on $M\otimes
\mathcal{B}^{\otimes n}$ by
\begin{equation*}
z\cdot \left( m\otimes b_{1}\otimes \cdots \otimes b_{n}\right) =(z\cdot
m)\otimes b_{1}\otimes \cdots \otimes b_{n},
\end{equation*}%
and the differential maps $b_{\ast }$ are morphisms of $\mathcal{Z}_{0}$%
-modules. Clearly, $\rho _{n}(M)$ is a morphism of $\mathcal{Z}_{0}$%
-modules, so $C{}_{n}(\mathcal{B},-)$ can be seen as a functor from $_{%
\mathcal{A}}\mathfrak{M}_{\mathcal{A}}^{\mathcal{H}}$ to the category of
chain complexes in $_{\mathcal{Z}_{0}}\mathfrak{M}^{\mathcal{H}}.$
Therefore, a fortiori, the functors $\mathrm{HH}{}_{\ast }(\mathcal{B},-)%
\mathcal{\ }$map a Hopf bimodule to an object in $_{\mathcal{Z}_{0}}%
\mathfrak{M}^{\mathcal{H}}.$ The $\mathcal{H}$-coaction on $\mathrm{HH}%
{}_{\ast }(\mathcal{B},M)$ will still be denoted by $\rho _{\ast }(M)$.
\end{noname}

\begin{remark}
By definition, Hochschild homology of $\mathcal{B}$ with coefficients in $M$
in degree zero equals $M_\mathcal{B}$. Thus in Proposition \ref%
{le:PropHopfBimod} (1) and \S \ref{nn:H-comodule}, we constructed two $%
\mathcal{H}$-coactions on $M_\mathcal{B}$, both of them being denoted by $%
\rho_0(M)$. The notation we have used is consistent, as these coactions are
identical.
\end{remark}

\begin{lemma}
\label{le:HH_B}Let $\mathcal{H}$ be a Hopf algebra with bijective antipode.
If ${\mathcal{B}\subseteq \mathcal{A}}$ is a faithfully flat $\mathcal{H}$%
-Galois extension then $\mathrm{HH}{}_{\ast }(\mathcal{B},-):{}_{\mathcal{A}%
}\mathfrak{M}_{\mathcal{A}}^{\mathcal{H}}\rightarrow {}_{\mathcal{Z}_{0}}%
\mathfrak{M}^{\mathcal{H}}$ is a homological and effaceable $\delta $%
-functor, where $\mathcal{Z}$ is the center of $\mathcal{A}$ and $\mathcal{Z}%
_{0}:=\mathcal{Z}\bigcap \mathcal{B}$.
\end{lemma}

\begin{proof}
We take a short exact sequence$\ $of Hopf bimodules
\begin{equation}
0\longrightarrow M^{\prime }\longrightarrow M\longrightarrow M^{\prime
\prime }\longrightarrow 0.  \label{eq:S}
\end{equation}%
We have to prove that there are connecting maps $\delta _{n}:\mathrm{HH}%
{}_{n}(\mathcal{B},M^{\prime \prime })\rightarrow \mathrm{HH}{}_{n-1}(%
\mathcal{B},M^{\prime }),$ which are homomorphisms of $\mathcal{Z}_{0}$%
-modules and $\mathcal{H}$-comodules, making
\begin{equation*}
\cdots \longrightarrow \mathrm{HH}{}_{n}(\mathcal{B},M^{\prime
})\longrightarrow \mathrm{HH}{}_{n}(\mathcal{B},M)\longrightarrow \mathrm{HH}%
{}_{n}(\mathcal{B},M^{\prime \prime })\overset{\delta _{n}}{\longrightarrow }%
\mathrm{HH}{}_{n-1}(\mathcal{B},M^{\prime })\longrightarrow \cdots
\end{equation*}%
a functorial exact sequence. In our setting, the connecting maps are
obtained by applying the long exact sequence in homology to the following
short exact sequence of complexes in $_{\mathcal{Z}_{0}}\mathfrak{M}^{%
\mathcal{H}}$
\begin{equation*}
0\longrightarrow C_{\ast }(\mathcal{B},M^{\prime })\longrightarrow C_{\ast }(%
\mathcal{B},M)\longrightarrow C_{\ast }(\mathcal{B},M^{\prime \prime
})\longrightarrow 0.
\end{equation*}%
Let us prove that $\mathrm{HH}{}_{\ast }(\mathcal{B},-)$ is effaceable too.
Let $M$ be a given Hopf bimodule. By Proposition \ref{pr:progenerator},
there exists a certain set $I$ such that $M$ is the quotient of $P:=(%
\mathcal{A}\otimes \mathcal{A})^{(I)}$ as a Hopf bimodule. In view of the
same proposition, $\mathcal{A}\otimes \mathcal{A}$ is projective as a $%
\mathcal{B}$-bimodule. Thus, for $n>0$,
\begin{equation*}
\mathrm{HH}{}_{n}(\mathcal{B},P)\cong \mathrm{Tor}_{n}^{\mathcal{B}^{e}}(%
\mathcal{B},P)=0.
\end{equation*}%
Hence the lemma is completely proven.
\end{proof}

\begin{remark}
Both $\mathrm{HH}{}_{\ast }(\mathcal{B},-)$ and $\mathrm{HH}{}_{\ast }(%
\mathcal{B},-)\otimes \mathcal{H}$ can be seen as homological and effaceable
$\delta $-functors that map a Hopf bimodule to an object in $_{\mathcal{Z}%
_{0}}\mathfrak{M}^\mathcal{H}$. The natural transformations $\rho _{\ast
}(-) $ in \S \ref{nn:H-comodule} define a morphisms of $\delta $-functors
that lifts%
\begin{equation*}
\rho _{0}(-):\mathrm{HH}_0(\mathcal{B},-)\longrightarrow \mathrm{HH}_0(%
\mathcal{B},-) \otimes \mathcal{H}.
\end{equation*}
\end{remark}

\begin{proposition}
\label{pr:HH_n(B,M)}Let $\mathcal{H}$ be a Hopf algebra with bijective
antipode. We assume that ${\mathcal{B}\subseteq \mathcal{A}}$ is a
faithfully flat $\mathcal{H}$-Galois extension and that $M$ is a Hopf
bimodule.

\begin{enumerate}
\item There is an $\mathcal{H}$-action on $\mathrm{HH}{}_{n}(\mathcal{B},M)$
that extends the module structure defined in (\ref{eq:R-module}). Moreover,
for any $h\in \mathcal{H}$ and $\omega \in \mathrm{HH}{}_{n}(\mathcal{B},M),$
\begin{equation}
\rho _{n}(M)(h\cdot \omega )=\textstyle\sum h_{(2)}\cdot \omega
_{\left\langle 0\right\rangle }\otimes h_{(3)}\omega _{\left\langle
1\right\rangle }S_{\mathcal{H}}h_{(1)}.  \label{ec:SAYD}
\end{equation}

\item If the antipode of $\mathcal{H}$ is involutive then $\mathrm{HH}%
{}_{\ast }(\mathcal{B},-)$ is a homological and effaceable $\delta $-functor
that takes values in $_{\mathcal{R}_{\mathcal{H}}\otimes \mathcal{Z}_{0}}%
\mathfrak{M}^{\mathcal{C}_{\mathcal{H}}}.$
\end{enumerate}
\end{proposition}

\begin{proof}
(1) We fix $h\in \mathcal{H}.$ The module structure constructed in formula (%
\ref{eq:R-module}) defines a natural map
\begin{equation*}
\mu _{0}^{h}(M):\mathrm{HH}{}_{0}(\mathcal{B},M)\rightarrow \mathrm{HH}%
{}_{0}(\mathcal{B},M),\quad \mu _{0}^{h}(M)([m]_{\mathcal{B}})=h\cdot
\lbrack m]_{\mathcal{B}}.
\end{equation*}%
In view of Lemma \ref{le:HH_B} the $\delta $-functor $\mathrm{HH}{}_{\ast }(%
\mathcal{B},-):{}_{\mathcal{A}}\mathfrak{M}_{\mathcal{A}}^{\mathcal{H}%
}\rightarrow {}_{\mathcal{Z}_{0}}\mathfrak{M}^{\mathcal{H}}$ is homological
and effaceable. Hence, by the universal property of these functors (see \S %
\ref{delta_functor}) there is a unique morphism of $\delta $-functors
\begin{equation*}
\mu _{\ast }^{h}(-):\mathrm{HH}{}_{\ast }(\mathcal{B},-)\rightarrow \mathrm{%
HH}{}_{\ast }(\mathcal{B},-)
\end{equation*}%
that lifts $\mu _{0}^{h}(-).$ Note that, by definition, $\mu _{\ast }^{h}(-)$
and the connecting morphisms $\delta _{\ast }$ are morphisms of $Z_0$%
-modules and $\mathcal{H}$-comodules. For $\omega \in \mathrm{HH}{}_{n}(%
\mathcal{B},M), $ we set%
\begin{equation*}
h\cdot \omega :=\mu _{n}^{h}(M)(\omega ).
\end{equation*}%
Proceeding as in the proof of \cite[Proposition 2.4]{St1}, one can easily
see that the above formula defines a natural action of $\mathcal{H}$ on $%
\mathrm{HH}{}_{n}(\mathcal{B},M)$. By construction, it lifts the action in (%
\ref{eq:R-module}). Note that, for any $n,$ the connecting maps $\delta _{n}$
are morphisms of $\mathcal{H}$-modules, since $\mu _{\ast }^{h}(-)$ is a
morphism of $\delta $-functors.

To conclude the proof of this part, it remains to prove relation (\ref%
{ec:SAYD}). We proceed by induction. In degree zero the required identity
holds by (\ref{ec:M_B=CrossedModule}). Let us assume that (\ref{ec:SAYD})
holds in degree $n$ for any Hopf bimodule. Let $M$ be a given Hopf bimodule.
We take an exact sequence%
\begin{equation*}
0\longrightarrow K\longrightarrow P\longrightarrow M\longrightarrow 0
\end{equation*}%
of Hopf bimodules such that $P:=(\mathcal{A}\otimes \mathcal{A})^{(I)}.$
Since $\delta _{n+1}$ is a homomorphism of $\mathcal{H}$-modules and $%
\mathcal{H}$-comodules and using the induction hypothesis, for $h\in
\mathcal{H}$ and $\omega \in $ $\mathrm{HH}{}_{n+1}(\mathcal{B},M),$ we get
\begin{align*}
(\delta _{n+1}\otimes \mathcal{H})\left( \rho _{n+1}(M)(h\cdot \omega
)\right) & =\rho _{n}(K)\left( \delta _{n+1}(h\cdot \omega )\right) \\
& =\rho _{n}(K)\left( h\cdot \delta _{n+1}(\omega )\right) \\
& =\textstyle\sum h_{(2)}\cdot \delta _{n+1}(\omega _{\left\langle
0\right\rangle })\otimes h_{(3)}\omega _{\left\langle 1\right\rangle }S_{%
\mathcal{H}}h_{(1)} \\
& =(\delta _{n+1}\otimes \mathcal{H})\left( \textstyle\sum h_{(2)}\cdot
\omega _{\left\langle 0\right\rangle }\otimes h_{(3)}\omega _{\left\langle
1\right\rangle }S_{\mathcal{H}}h_{(1)}\right) .
\end{align*}%
As $\mathrm{HH}{}_{n+1}(\mathcal{B},P)=0$ it follows that $\delta _{n+1}$ is
injective. Consequently, $\delta _{n+1}\otimes \mathcal{H}$ is also
injective. Thus the foregoing computation implies relation (\ref{ec:SAYD}).

(2) Let $M$ be a given Hopf bimodule. For $z\in \mathcal{Z}_{0}$ we define
\begin{equation*}
\nu _{n}^{z}(M):\mathrm{HH}{}_{n}(\mathcal{B},M)\rightarrow \mathrm{HH}%
{}_{n}(\mathcal{B},M),\quad \nu _{n}^{z}(M)\left( \omega \right) =z\cdot
\omega .
\end{equation*}%
We claim that $\mu _{\ast }^{h}(-)$ and $\nu _{\ast }^{z}(-)$ commute for
all $h\in \mathcal{H}$, i.e.
\begin{equation}
\mu _{\ast }^{h}(-)\circ \nu _{\ast }^{z}(-)=\nu _{\ast }^{z}(-)\circ \mu
_{\ast }^{h}(-).  \label{eq:commuting}
\end{equation}%
In degree zero this identity follows from the computation below, where for
brevity we write $\mu _{0}^{h}$ and $\nu _{0}^{z}$ instead of $\mu
_{0}^{h}(M)$ and $\nu _{0}^{z}(M).$ Indeed,
\begin{equation*}
\textstyle\left( \mu _{0}^{h}\circ \nu _{0}^{z}\right) \left( \left[ m\right]
_{\mathcal{B}}\right) =\sum \left[ \kappa ^{2}(h)zm\kappa ^{1}(h)\right] _{%
\mathcal{B}}=z\cdot \sum \left[ \kappa ^{2}(h)m\kappa ^{1}(h)\right] _{%
\mathcal{B}}=\left( \nu _{0}^{z}\circ \mu _{0}^{h}\right) \left( \left[ m%
\right] _{\mathcal{B}}\right) ,
\end{equation*}%
where for the second equality we used that $z$ is a central element.
Furthermore, the natural transformations that appear in the left and right
hand sides of (\ref{eq:commuting}) are morphisms of $\delta $-functors that
lift respectively $\mu _{0}^{h}(-)\circ \nu _{0}^{z}(-)$ and $\nu
_{0}^{z}(-)\circ \mu _{0}^{h}(-)$. Hence (\ref{eq:commuting}) follows by the
universal property of homological and $\delta $-functors. In view of the
relation (\ref{eq:commuting}) it follows that $\mathrm{HH}{}_{n}(\mathcal{B}%
,M)$ is an $\mathcal{H}\otimes \mathcal{Z}_{0}$-module with respect to
\begin{equation*}
\left( h\otimes z\right) \cdot \omega =\left[ \mu _{n}^{h}(M)\circ \nu
_{n}^{z}(M)\right] (\omega ).
\end{equation*}%
We can now prove that $\mathrm{HH}{}_{n}(\mathcal{B},M)$ is an object in $_{%
\mathcal{R}_{\mathcal{H}}\otimes \mathcal{Z}_{0}}\mathfrak{M}^{\mathcal{C}_{%
\mathcal{H}}}.$ As $\mathcal{R}_{\mathcal{H}}\otimes $ $\mathcal{Z}_{0}$ is
subalgebra of $\mathcal{H}\otimes \mathcal{Z}_{0}$, it acts on $\mathrm{HH}%
{}_{n}(\mathcal{B},M)$. The coalgebra $\mathcal{C}_{\mathcal{H}}$ coacts on
the Hochschild homology of $\mathcal{B}$ with coefficients in $M$ via
\begin{equation*}
\bar{\rho}_\ast(M):=\big(\mathrm{HH}{}_{\ast}(\mathcal{B},M)\otimes \pi_%
\mathcal{H}\big)\otimes \rho(M).
\end{equation*}
To simplify the notation, we shall write $\rho _{\ast}(M)$ instead of $\bar{%
\rho}_ \ast(M)$. We have to show that $\rho _{n}(M)$ is a morphism of $%
\mathcal{R}_{\mathcal{H}}\otimes \mathcal{Z}_{0}$-modules. By Lemma \ref%
{le:HH_B} we already know that $\rho _{n}(M)$ is a morphism of $\mathcal{Z}%
_{0}$-modules. Thus, it remains to check that $\rho _{n}(M)$ is a morphism
of $\mathcal{R}_{\mathcal{H}}$-modules too. For the case $n=0$ see the proof
of Proposition \ref{le:PropHopfBimod} (2). In fact, the same proof works for
an arbitrary $n,$ just replacing $[m]_{\mathcal{B}}$ by an element $\omega
\in \mathrm{HH}{}_{n}(\mathcal{B},M)$ and using (\ref{ec:SAYD}) instead of (%
\ref{ec:M_B=CrossedModule}).

We still have to prove that $\mathrm{HH}{}_{\ast }(\mathcal{B},-):{}_{%
\mathcal{A}}\mathfrak{M}_{\mathcal{A}}^{\mathcal{H}}\rightarrow {}_{\mathcal{%
R}_{\mathcal{H}}\otimes \mathcal{Z}_{0}}\mathfrak{M}^{\mathcal{C}_{\mathcal{H%
}}}$ is a homological and effaceable $\delta $-functor, i.e. for every short
exact of Hopf bimodules the corresponding connecting maps $\delta _{\ast}$
are morphisms of $\mathcal{R}_{\mathcal{H}}\otimes \mathcal{Z}_{0}$-modules
and $\mathcal{H}$-comodules. By Lemma \ref{le:HH_B}, it follows that $\delta
_{\ast}$ are morphisms of $\mathcal{Z}_{0}$-modules and $\mathcal{H}$%
-comodules. By the proof of the first part of the proposition, $\delta
_{\ast}$ are also morphisms of $\mathcal{H}$-modules. Hence, a fortiori,
they are morphisms of $\mathcal{R}_{\mathcal{H}}$-modules.
\end{proof}

\begin{noname}
\label{nn:FunctorG} The natural transformations that define the $\mathcal{R}%
_{\mathcal{H}}$-module and the $\mathcal{C}_{\mathcal{H}}$-comodule
structures of $\mathrm{HH}{}_{\ast }(\mathcal{B},-),$ as in the above
proposition, will be denoted by $\mu _{\ast }(-)$ and $\rho _{\ast }(-),$
respectively.

Let us take an injective left $\mathcal{C}_{\mathcal{H}}$-comodule $V$. By
\S \ref{nn:cotensor}, for a Hopf bimodule $M,$ the cotensor product $\mathrm{%
HH}{}_{\ast }(\mathcal{B},M)\square _{\mathcal{C}_{\mathcal{H}}}V$ is a left
$\mathcal{R}_{\mathcal{H}}\otimes \mathcal{Z}_{0}$-module. It follows that $%
\mathrm{HH}{}_{\ast }(\mathcal{B},-)\square _{\mathcal{C}_{\mathcal{H}}}V%
\mathrm{\ }$is a homological and effaceable functor from $\ $the category of
Hopf bimodules to the category of left $\mathcal{R}_{\mathcal{H}}\otimes
\mathcal{Z}_{0}$-modules. Of course, its connecting maps are $\delta _{\ast
}\square _{\mathcal{C}_{\mathcal{H}}}V,$ where $\delta _{\ast }$ are the
connecting homomorphisms of the $\delta $-functor $\mathrm{HH}{}_{\ast }(%
\mathcal{B},-).$

To simplify the notation, we shall denote $\mathrm{HH}{}_{0}(\mathcal{B}%
,-)\square _{\mathcal{C}_{\mathcal{H}}}V$ by $G_{V}.$ By the foregoing
observations $G_{V}$ maps a Hopf bimodule to a left $\mathcal{R}_{\mathcal{H}%
}\otimes \mathcal{Z}_{0}$-module. Our aim now is to describe the left
derived functors of $G_{V}$.
\end{noname}

\begin{proposition}
\label{pr:G_V}Let $\mathcal{H}$ be a Hopf algebra such that $S_\mathcal{H}^2=%
\mathrm{Id}_\mathcal{H}$. If ${\mathcal{B}\subseteq \mathcal{A}}$ is a
faithfully flat $\mathcal{H}$ -Galois extension and $V$ is an injective $%
\mathcal{C}_{\mathcal{H}}$-comodule, then
\begin{equation*}
\mathrm{HH}{}_{\ast }(\mathcal{B},-\square _{\mathcal{C}_{\mathcal{H}%
}}V):{}_{\mathcal{A}}\mathfrak{M}_{\mathcal{A}}^{\mathcal{H}}\longrightarrow
{}_{\mathcal{R}_{\mathcal{H}}\otimes \mathcal{Z}_{0}}\mathfrak{M}
\end{equation*}%
is a homological and effaceable $\delta $-functor. As $\delta$-functors from
$_{\mathcal{A}}\mathfrak{M}_{\mathcal{A}}^{\mathcal{H}}$ to $_{\mathcal{R}_{%
\mathcal{H}}\otimes \mathcal{Z}_{0}}\mathfrak{M}$,
\begin{equation}
\mathrm{L}_{\ast}G_{V}\cong \mathrm{HH}{}_{\ast}(\mathcal{B},-)\square _{%
\mathcal{C}_{\mathcal{H}}}V\cong \mathrm{HH}{}_{\ast }\left( \mathcal{B}%
,-\square _{\mathcal{C}_{\mathcal{H}}}V\right).  \label{ec:izomorphisme}
\end{equation}
\end{proposition}

\begin{proof}
By Proposition \ref{le:PropHopfBimod} (3) the cotensor product $M\square _{%
\mathcal{C}_{\mathcal{\ H}}}V$ is a $\mathcal{B}$-bimodule, for every Hopf
bimodule $M$. Hence Hochschild homology of $\mathcal{B}$ with coefficients
in $M\square _{\mathcal{C}_{\mathcal{H}}}V$ makes sense. We set $T_{\ast }:=%
\mathrm{HH}{}_{\ast }(\mathcal{B},-\square _{\mathcal{C}_{\mathcal{H}}}V)$
and take a short exact sequence of Hopf bimodules as in (\ref{eq:S}). Since $%
V$ is an injective comodule,
\begin{equation*}
0\longrightarrow C_{\ast }(\mathcal{B},M^{\prime }\square _{\mathcal{C}_{%
\mathcal{H}}}V)\longrightarrow C_{\ast }(\mathcal{B},M\square _{\mathcal{C}_{%
\mathcal{H}}}V)\longrightarrow C_{\ast }(\mathcal{B},M^{\prime
\prime}\square _{\mathcal{C}_{\mathcal{H}}}V)\longrightarrow 0
\end{equation*}%
is exact. By the definition of $T_{\ast }$ and the long exact sequence in
homology we deduce that $T_{\ast }$ is homological, regarded as a $\delta $%
-functor to the category of vector spaces. Recall that $U:=\mathcal{A}%
\otimes \mathcal{A}$ is a generator in the category of Hopf bimodules.
Therefore, to prove that $T_{\ast }$ is effaceable, it is enough to show
that $T_{n}\left( X\right) =0,$ where $X$ is an arbitrary direct sum of
copies of $U$ and $n>0.$ In fact, as Hochschild homology and the cotensor
product commute with direct sums, we may assume that $X=U$. We claim that $%
U\square _{\mathcal{C}_{\mathcal{H}}}V$ is flat as a $\mathcal{B}$-bimodule.
By (\ref{eq:tens-cotens}), for an arbitrary $\mathcal{B}$-bimodule $N,$
\begin{equation*}
N\otimes _{\mathcal{B}^{e}}\left( U\square _{\mathcal{C}_{\mathcal{H}%
}}V\right) \cong \left( N\otimes _{\mathcal{B}^{e}}\mathcal{A}^{e}\right)
\square _{\mathcal{C}_{\mathcal{H}}}V\cong \left( \mathcal{A}\otimes _{%
\mathcal{B}}N\otimes _{\mathcal{B}}\mathcal{A}\right) \square _{\mathcal{C}_{%
\mathcal{H}}}V.
\end{equation*}%
Hence, the functors $(-)\otimes _{\mathcal{B}^{e}}\left( U\square _{\mathcal{C%
}_{\mathcal{H}}}V\right) $ and $\left( -\square _{\mathcal{C}_{\mathcal{H}%
}}V\right) \circ \left( \mathcal{A}\otimes _{\mathcal{B}}-\otimes _{\mathcal{%
B}}\mathcal{A}\right) $ are isomorphic. Since the antipode of $\mathcal{H}$
is bijective, by \cite[Theorems 4.9 and 4.10]{SS}, $\mathcal{A}$ is
faithfully flat as a left and a right $\mathcal{B}$-module. Therefore, $%
\mathcal{A}\otimes _{\mathcal{B}}-\otimes _{\mathcal{B}}\mathcal{A}$ is an
exact functor. As $V$ is injective, the functor $(-)\square _{\mathcal{C}_{%
\mathcal{H}}}V$ is also exact, so $U\square _{\mathcal{C}_{\mathcal{H}}}V$
is flat. Thus
\begin{equation*}
T_{n}\left( U\right) \cong \mathrm{Tor}{}_{n}^{\mathcal{B}^{e}}(\mathcal{B}%
,U\square _{\mathcal{C}_{\mathcal{H}}}V)=0.
\end{equation*}%
Summarizing, $T_{\ast }:{}_\mathcal{A}\mathfrak{M}_\mathcal{A}^\mathcal{H}
\longrightarrow{}_\mathbb{K}\mathfrak{M}$ is a homological and effaceable $%
\delta $-functor.

For each Hopf bimodule $M$, our aim now is to endow $T{}_{n}\left( M\right) $
with a left module structure over $\mathcal{R}_{\mathcal{H}}\otimes \mathcal{%
Z}_{0}$. Let us first consider the case $n=0.$ By Proposition \ref%
{le:PropHopfBimod} (3), there is a canonical left $\mathcal{R}_{\mathcal{H}}$%
-action on $T{}_{0}\left( M\right) $ such that
\begin{equation}
T{}_{0}\left( M\right) \cong \mathrm{HH}{}_{0}(\mathcal{B},M)\square _{%
\mathcal{C}_{\mathcal{H}}}V,  \label{ec:T_0}
\end{equation}%
the natural $\mathbb{K}$-linear isomorphism constructed in (\ref{eq:iso_B}),
is a homomorphism of $\mathcal{R}_{\mathcal{H}}$-modules. As $M$ is an
object in $_{\mathcal{Z}_{0}}\mathfrak{M}^{\mathcal{C}_{\mathcal{H}}},$ by
\S \ref{nn:cotensor}, it follows that $M\square _{\mathcal{C}_{\mathcal{H}%
}}V $ is a left $\mathcal{Z}_{0}$-submodule of $M\otimes V.$ In particular
this cotensor product is a left $\mathcal{Z}_{0}$-module. Furthermore, $%
T{}_{0}\left( M\right) $ is a quotient $\mathcal{Z}_{0}$-module of $M\square
_{\mathcal{C}_{\mathcal{H}}}V,$ as the commutator space $\left[ \mathcal{B}%
,M\square _{\mathcal{C}_{\mathcal{H}}}V\right] $ is a $\mathcal{Z}_{0}$%
-submodule of $M\square _{\mathcal{C}_{\mathcal{H}}}V.$ Obviously, with
respect to this module structure, $T{}_{0}\left( M\right) $ becomes a module
over $\mathcal{R}_{\mathcal{H}}\otimes \mathcal{Z}_{0}$ and the isomorphism
in (\ref{ec:T_0}) is a map of $\mathcal{R}_{\mathcal{H}}\otimes \mathcal{Z}%
_{0}$-modules. Since $T_{\ast }$ is a homological and effaceable functor,
one can proceed as in the proof of Proposition \ref{pr:HH_n(B,M)} to lift
the $\mathcal{R}_{\mathcal{H}}\otimes \mathcal{Z}_{0}$-action on $%
T{}_{0}\left( M\right) $ to a natural $\mathcal{R}_{\mathcal{H}}\otimes
\mathcal{Z}_{0}$-module structure on $T_{\ast }(M),$ for every Hopf bimodule
$M.\ \ $Again as in the proof of the above mentioned result, we can show
that $T_{\ast }:{}_{\mathcal{A}}\mathfrak{M}_{\mathcal{A}}^{\mathcal{H}%
}\rightarrow {}_{\mathcal{R}_{\mathcal{H}}\otimes \mathcal{Z}_{0}}\mathfrak{M%
}$ is homological and effaceable.

It remains to prove the isomorphisms in (\ref{ec:izomorphisme}). Note that
the left derived functors of a right exact functor define a homological and
effaceable $\delta $-functor. Thus $\mathrm{L}_{\ast }G_{V}$ is a
homological and effaceable $\delta $-functor. Clearly, $\mathrm{L}_{0}G_{V}=%
\mathrm{HH}_{0}(\mathcal{B,-)}\square _{\mathcal{C}_{\mathcal{H}}}V$. Hence,
by the universal property of homological and effaceable $\delta $-functors,
this identity may be lifted to give the first isomorphism in (\ref%
{ec:izomorphisme}). The second isomorphism is obtained in a similar manner,
by lifting the natural transformation in (\ref{ec:T_0}).
\end{proof}

\begin{noname}
\label{nn:FunctorH}Let $\mathcal{B\subseteq }{\mathcal{A}}$ be an $\mathcal{H%
}$-comodule algebra and let $M$ be a Hopf bimodule. Following \cite[Theorem
1.3]{St2} we regard $C_{\ast }(\mathcal{A},M)$ as a complex in the category $%
\mathfrak{M}^{\mathcal{C}_{\mathcal{H}}}$ with respect to the coaction that
in degree $n$ is given by
\begin{equation*}
\varrho_n(M) (m\otimes a^{1}\otimes \cdots \otimes a^{n})=\sum
m_{\left\langle 0\right\rangle }\otimes a_{\left\langle 0\right\rangle
}^{1}\otimes \cdots \otimes a_{\left\langle 0\right\rangle }^{n}\otimes \pi
_{\mathcal{H}}\left( m_{\left\langle 1\right\rangle }a_{\left\langle
1\right\rangle }^{1}\cdots a_{\left\langle 1\right\rangle }^{n}\right) .
\end{equation*}%
Recall that $\pi _{\mathcal{H}}$ denotes the projection of $\mathcal{H}$
onto $\mathcal{C}_{\mathcal{H}}$ and that $\mathcal{Z}_{0}=\mathcal{Z}%
\bigcap \mathcal{B}.$ It is not difficult to see that $C_{\ast }(\mathcal{A}%
,M)$ is a complex of left $\mathcal{Z}_{0}$-modules with respect to the
action that in degree $n$ is defined by%
\begin{equation*}
z\cdot (m\otimes a^{1}\otimes \cdots \otimes a^{n})=zm\otimes a^{1}\otimes
\cdots \otimes a^{n}.
\end{equation*}%
In fact, since $\mathcal{Z}_0$ contains only coinvariant elements, it
follows that $C_{\ast }(\mathcal{A},M)$ is a complex in $_{\mathcal{Z}_{0}}%
\mathfrak{M}^{\mathcal{C}_{\mathcal{H}}}$. Therefore, $\mathrm{HH}_{n}(%
\mathcal{A},M)$ is an object in the same category, for every $n$. In view of
\S \ref{nn:cotensor} it follows that $\mathrm{HH}_{\ast }(\mathcal{A}%
,M)\square _{\mathcal{C}_{\mathcal{H}}}V$ is a left $\mathcal{Z}_{0}$%
-module, for any injective left $\mathcal{C}_{\mathcal{H}}$-comodule $V$.
Therefore
\begin{equation*}
H_{V}:{}_{\mathcal{A}}\mathfrak{M}_{\mathcal{A}}^{\mathcal{H}}\rightarrow
{}_{\mathcal{Z}_{0}}\mathfrak{M,\quad }H_{V}(M):=M_{\mathcal{A}}\square _{%
\mathcal{C}_{\mathcal{H}}}V.
\end{equation*}%
is a well defined functor, as by the foregoing remarks $M_{\mathcal{A}}=%
\mathrm{HH}{}_{0}(\mathcal{A},M)$ is a right $\mathcal{C}_{\mathcal{H}}$%
-module and $H_{V}(M)$ is a $\mathcal{Z}_{0}$-module.
\end{noname}

\begin{proposition}
\label{pr:H_V}Let ${\mathcal{B}\subseteq \mathcal{A}}$ be a faithfully flat $%
\mathcal{H}$-Galois extension, where $\mathcal{H}$ is a Hopf algebra with
bijective antipode. If $V$ is an injective $\mathcal{C}_{\mathcal{H}}$%
-comodule, then there is an isomorphism of $\delta $-functors%
\begin{equation}
\mathrm{L}_{\ast}H_{V}\cong \mathrm{HH}{}_{\ast}(\mathcal{A},-) \square _{%
\mathcal{C}_{\mathcal{H}}}V.  \label{ec:LH_V}
\end{equation}
\end{proposition}

\begin{proof}
First, let us show that $T_{\ast }:=\mathrm{HH}{}_{\ast }(\mathcal{A}%
,-)\square _{\mathcal{C}_{\mathcal{H}}}V$ is a homological and effaceable $%
\delta $-functor to the category of left $\mathcal{Z}_{0}$-modules. For a
short exact sequence as in (\ref{eq:S}),
\begin{equation*}
0\longrightarrow C_{\ast }(\mathcal{A},M^{\prime })\longrightarrow C_{\ast }(%
\mathcal{A},M)\longrightarrow C_{\ast }(\mathcal{A},M^{\prime \prime
})\longrightarrow 0
\end{equation*}%
is an exact sequence of complexes in $_{\mathcal{Z}_{0}}\mathfrak{M}^{%
\mathcal{C}_{\mathcal{H}}}.$ Therefore, the corresponding long exact
sequence in homology lives in the same category. In particular, its
connecting maps $\delta _{\ast }$ are morphisms of $\mathcal{Z}_{0}$-modules
and $\mathcal{C}_{\mathcal{H}}$-comodules, so are $\mathcal{Z}_{0}$-linear.
Since $V$ is injective, the functor $(-)\square _{\mathcal{C}_{\mathcal{H}%
}}V $ is exact. Thus $T_{\ast }$ is a homological functor with connecting
maps $\delta _{\ast}\square _{\mathcal{C}_{\mathcal{H}}}V$.

By Proposition \ref{pr:progenerator}, the Hopf bimodule $U:=\mathcal{A}%
\otimes \mathcal{A}$ is a generator. We also have
\begin{equation*}
\mathrm{HH}{}_{n}(\mathcal{A},U)\square _{\mathcal{C}_{\mathcal{H}}}V\cong
\mathrm{Tor}{}_{n}^{\mathcal{A}^{e}}(\mathcal{A},\mathcal{A}^{e})\square _{%
\mathcal{C}_{\mathcal{H}}}V=0.
\end{equation*}%
In conclusion $T_{\ast }$ is effaceable, as Hochschild homology and the
cotensor product commute with direct sums. The isomorphisms in (\ref{ec:LH_V}%
) are obtained by lifting the identity $\mathrm{L}_{0}H_{V}=T_{0}$, as in
the proof of Proposition \ref{pr:G_V}.
\end{proof}

\begin{noname}
Since $\mathcal{H}$ is a Hopf algebra, the category of right $\mathcal{H}$%
-comodules is monoidal, with respect to the tensor product of vector spaces,
on which we put the diagonal coaction. More precisely, if $V$ and $W$ are
right $\mathcal{H}$-comodules then $\mathcal{H}$ coacts on $V\otimes W$ via
the map
\begin{equation*}
\rho _{V\otimes W}(v\otimes w)=\textstyle\sum v_{\left\langle 0\right\rangle
}\otimes w_{\left\langle 0\right\rangle }\otimes v_{\left\langle
1\right\rangle }w_{\left\langle 1\right\rangle },\quad
\end{equation*}%
where $v\in V,\ w\in W.$ Let us denote the right adjoint coaction of $%
\mathcal{H}$ on itself by $\mathcal{H}^{\mathrm{ad}}.$ Recall that the map $%
\rho _{\mathcal{H}}:\mathcal{H}\rightarrow \mathcal{H}\otimes \mathcal{H}$
that define this coaction is given by%
\begin{equation*}
\rho _{\mathcal{H}}(h)=\textstyle\sum h_{(2)}\otimes Sh_{(1)}h_{(3)}.
\end{equation*}%
Hence $\mathcal{A}\otimes \mathcal{H}^{\mathrm{ad}}$ is a right $\mathcal{H}$%
-comodule. Consequently, it is a $\mathcal{C}_{\mathcal{H}}$-comodule via
\begin{equation}  \label{ro:AxH}
\rho _{\mathcal{A}\otimes \mathcal{H}}(a\otimes h)=\textstyle\sum
a_{\left\langle 0\right\rangle }\otimes h_{(2)}\otimes \pi _{\mathcal{H}%
}\left( a_{\left\langle 1\right\rangle }Sh_{(1)}h_{(3)}\right) .\quad
\end{equation}
\end{noname}

\begin{noname}
For a Hopf algebra $\mathcal{H}$ let
\begin{equation*}
\mathcal{H}^{+}:=\ker\varepsilon_\mathcal{H}\qquad\text{and}\qquad\mathcal{R}%
_{\mathcal{H}}^{+}:=\mathcal{H}^{+}\cap \mathcal{R}_{\mathcal{H}}.
\end{equation*}
If $X$ is an $\mathcal{R}_{\mathcal{H}}\otimes \mathcal{Z}_{0}$-module, then
$\mathcal{R}_{\mathcal{H}}^{+}X$ is a $\mathcal{Z}_{0}$-submodule of $X$, as
the corresponding actions of $\mathcal{R}_{\mathcal{H}}$ and $\mathcal{Z}%
_{0} $ on $X$ commute. For the same reason, $\mathbb{K}\otimes _{\mathcal{R}%
_{\mathcal{H}}}X$ is a $\mathcal{Z}_{0}$-module. Obviously, with respect to
these module structures, the canonical isomorphism
\begin{equation*}
\mathbb{K}\otimes _{\mathcal{R}_{\mathcal{H}}}X\cong X/\left( \mathcal{R}_{%
\mathcal{H}}^{+}X\right)
\end{equation*}%
is $\mathcal{Z}_{0}$-linear. For example, $\mathcal{Z}_{0}$ acts on $\mathbb{%
K}\otimes _{\mathcal{R}_{\mathcal{H}}}G_{V}(M)$ such that, for
\begin{equation}
\zeta =\textstyle\sum_{i=1}^{n}\left[ m_{i}\right] _{\mathcal{B}}\otimes
v_{i}  \label{ec:csi}
\end{equation}%
in $G_{V}(M)$ and $z$ in $\mathcal{Z}_{0}$ we have
\begin{equation*}
z\cdot \left( 1\otimes _{\mathcal{R}_{\mathcal{H}}}\zeta \right) =\textstyle%
\sum_{i=1}^{n}1\otimes _{\mathcal{R}_{\mathcal{H}}}\left[ zm_{i}\right] _{%
\mathcal{B}}\otimes v_{i}.
\end{equation*}%
In view of this observation, we shall regard $\mathbb{K}\otimes _{\mathcal{R}%
_{\mathcal{H}}}G_{V}$ as a functor from the category of Hopf bimodules to
the category of left $\mathcal{Z}_{0}$-modules.
\end{noname}

\begin{lemma}
\label{pr:Tor=comodule}Let $\mathcal{R}$ and $\mathcal{C}$ denote an
algebra and a coalgebra, respectively, over a field $\mathbb{K}$. If
$(M,\rho_M )\in {}_{\mathcal{R}}\mathfrak{M}^{\mathcal{C}}$ then
there is a unique morphism of $\delta $-functors
\begin{equation*}
\rho _{\ast }(-):\mathrm{Tor}_{\ast }^{\mathcal{R}}(-,M)\rightarrow \mathrm{%
Tor}_{\ast }^{\mathcal{R}}(-,M)\otimes \mathcal{C}
\end{equation*}%
that lifts $\rho _{0}(-):=(-)\otimes _{\mathcal{R}}\rho_M $ and defines a $%
\mathcal{C}$-coaction on $\mathrm{Tor}_{\ast }^{\mathcal{R}}(-,M).$
\end{lemma}

\begin{proof}
Obviously, $\mathrm{Tor}_{\ast }^{\mathcal{R}}(-,M)$ and $\mathrm{Tor}_{\ast
}^{\mathcal{R}}(-,M)\otimes \mathcal{C}$ are homological and effaceable $%
\delta $-functors, which are defined on the category of right $\mathcal{R}$%
-modules. In degree zero, $\rho _{0}(-):=(-)\otimes_\mathcal{R} \rho_M $
defines a $\mathcal{C}$-comodule structure on $\mathrm{Tor}_{0}^{\mathcal{R}%
}(-,M)=(-)\otimes _{\mathcal{R}}\mathcal{C}$. By the universal property,
there is a morphism of $\delta $-functors
\begin{equation*}
\rho _{\ast }(-):\mathrm{Tor}_{n}^{\mathcal{R}}(-,M)\rightarrow \mathrm{Tor}%
_{n}^{\mathcal{R}}(-,M)\otimes \mathcal{C}
\end{equation*}%
that lifts $\rho _{0}.$ We want to prove that $\rho _{\ast }(-)$ defines a
coaction on $\mathrm{Tor}_{\ast }^{\mathcal{R}}(-,M).$ We have already
remarked that this property holds in degree zero, so
\begin{equation}
\big( \rho _{0}(-)\otimes \mathcal{C}\big) \circ \rho _{0}(-)=\left(
-\otimes \Delta _{\mathcal{C}}\right) \circ \rho _{0}(-).
\label{eq:C_comodule}
\end{equation}%
We need a similar identity for $\rho _{n}(-)$. Note that $\mathrm{Tor}_{\ast
}^{\mathcal{R}}(-,M)\otimes \mathcal{C} \otimes \mathcal{C}$ is also a
homological and effaceable $\delta $-functor. Clearly, $\big( \rho _{\ast
}(-)\otimes \mathcal{C}\big) \circ \rho _{\ast }(-)$ and $( -\otimes \Delta
_{\mathcal{C}}) \circ \rho _{\ast }(-)$ lift the natural transformations in
the left and respectively right hand sides of (\ref{eq:C_comodule}). Hence,
by the uniqueness of the lifting, these morphisms of $\delta $-functors are
equal (see the universal property in \S \ref{delta_functor}).
\end{proof}

\begin{remark}
In view of the previous lemma, $\rho _{\ast }(-)$ is a homomorphism of $%
\delta $-functors. Hence, for an exact sequence of right $\mathcal{R}$%
-modules as in \eqref{eq:short_exact} and an object $M$ in $_{\mathcal{R}}%
\mathfrak{M}^{\mathcal{C}},$ the connecting maps
\begin{equation*}
\delta _{\ast}:\mathrm{Tor}_{\ast}^{\mathcal{R}}(X^{\prime\prime},M)%
\rightarrow \mathrm{Tor}_{\ast-1}^{\mathcal{R}}(X^{\prime},M)
\end{equation*}%
are morphisms of $\mathcal{C}$-comodules.
\end{remark}

\begin{lemma}
\label{le: fi}Let ${\mathcal{B}\subseteq \mathcal{A}}$ be a faithfully flat $%
\mathcal{H}$-Galois extension over a Hopf algebra such that $S_\mathcal{H}^2=%
\mathrm{Id}_\mathcal{H}$. Let $V$ be an injective $\mathcal{C}_{\mathcal{H}}$%
-comodule and set $U:=\mathcal{A}\otimes \mathcal{A}.$

\begin{enumerate}
\item Let $G_{V}$ and $H_{V}$ be the functors defined respectively in \emph{%
\S \ref{nn:FunctorG}} and \emph{\S \ref{nn:FunctorH}}. For every Hopf
bimodule $M$ and $\zeta$ as in (\ref{ec:csi}), the formula
\begin{equation*}
\phi(M)(1\otimes \zeta)=\textstyle\sum_{i=1}^{n}\left[ m_{i}\right] _{%
\mathcal{A}}\otimes v_{i}
\end{equation*}
defines a natural transformation $\phi(-) :\mathbb{K}\otimes _{\mathcal{R}%
_{H}}G_{V}(-)\longrightarrow H_{V}(-)$.

\item The algebra $\mathcal{R}_\mathcal{H}$ acts on $\mathcal{A}\otimes
\mathcal{H}$ via the multiplication in $\mathcal{H}$ so that $\mathcal{A}%
\otimes \mathcal{H}$ is an object in ${}_{\mathcal{R}_{\mathcal{H}}}%
\mathfrak{M}^{\mathcal{C}_{\mathcal{H}}}$ with respect to the coaction (\ref%
{ro:AxH}).

\item The $\mathcal{C_H}$-coaction on $\mathcal{A} \otimes \mathcal{H}$
induces a comodule structure on $\mathrm{Tor}_{\ast}^{\mathcal{R}_{\mathcal{H%
}}}(\mathbb{K},\mathcal{A}\otimes \mathcal{H})$.

\item We assume, in addition, that $\mathcal{R}_{\mathcal{H}}^{+}\mathcal{H}%
= \mathcal{H}^{+}$ and $\mathrm{Tor}_{n }^{\mathcal{R}_{\mathcal{H}}}(%
\mathbb{K},\mathcal{H})=0,$ for every $n>0$. Then $\phi (U)$ is an
isomorphism and $\mathrm{\mathrm{Tor}}_{n}^{\mathcal{R}_{\mathcal{H}}}\left(
\mathbb{K},G_{V}(U)\right) =0$, for $n>0$.
\end{enumerate}
\end{lemma}

\begin{proof}
(1) Let $M$ be a Hopf bimodule and let $p(M):M_{\mathcal{B}}\rightarrow M_{%
\mathcal{A}}$ denote the canonical projection. Clearly, $p(M)$ is a natural
morphism of $\mathcal{C}_{\mathcal{H}}$-modules. To simplify the notation,
set $p:=p(M).$ By the foregoing, $f:=p\;\square _{\mathcal{C}_{\mathcal{H}%
}}V $ is well-defined. Furthermore, by \cite[Proposition 2.6]{JS},
\begin{equation*}
\textstyle\sum \kappa ^{1}(r)\kappa ^{2}(r)=\varepsilon (r).
\end{equation*}%
Let $\zeta $ be an element in $M_{\mathcal{B}}\square _{\mathcal{C}_{%
\mathcal{H}}}V$ satisfying relation (\ref{ec:csi}). Since $[am]_{\mathcal{A}%
}=[ma]_{\mathcal{A}},$ for every $r\in \mathcal{R}_{\mathcal{H}}$, a
straightforward computation yields
\begin{equation*}
f(r\cdot \zeta )=\textstyle\sum_{i=1}^{n}[\kappa ^{2}(r)m_{i}\kappa
^{1}(r)]_{\mathcal{A}}\otimes v_{i}=\varepsilon (r)f(\zeta ).
\end{equation*}%
Thus, there exists a natural map $\phi (M):\mathbb{K}\otimes _{\mathcal{R}%
_{H}}G_{V}(M)\rightarrow H_{V}(M)$ of $\mathcal{Z}_{0}$-modules, which is
uniquely defined such that
\begin{equation*}
\phi (M)(1_{\mathbb{K}}\otimes _{\mathcal{R}_{\mathcal{H}}}\zeta )=f(\zeta ).
\end{equation*}

(2) We regard $\mathcal{A}\otimes \mathcal{H}$ as a left $\mathcal{R}_{%
\mathcal{H}}$-module via the multiplication in $\mathcal{H}.$ Let us prove
that $\rho _{\mathcal{A}\otimes \mathcal{H}}$ is a morphism of $\mathcal{R}_{%
\mathcal{H}}$-modules. We pick up $a\in \mathcal{A},$ $h\in \mathcal{H}$ and
$r\in \mathcal{R}_{\mathcal{H}}.$ Thus%
\begin{align*}
\rho _{\mathcal{A}\otimes \mathcal{H}}(a\otimes rh)& =\textstyle\sum
a_{\left\langle 0\right\rangle }\otimes r_{(2)}h_{(2)}\otimes \pi _{\mathcal{%
H}}\left( a_{\left\langle 1\right\rangle }S_{\mathcal{H}%
}(r_{(1)}h_{(1)})r_{(3)}h_{(3)}\right) \\
& =\textstyle\textstyle\sum a_{\left\langle 0\right\rangle }\otimes
r_{(2)}h_{(2)}\otimes \pi _{\mathcal{H}}\left( a_{\left\langle
1\right\rangle }S_{\mathcal{H}}h_{(1)}S_{\mathcal{H}}r_{(1)}r_{(3)}h_{(3)}%
\right) \\
& =\textstyle\sum a_{\left\langle 0\right\rangle }\otimes
r_{(3)}h_{(2)}\otimes \pi _{\mathcal{H}}\left( a_{\left\langle
1\right\rangle }S_{\mathcal{H}}h_{(1)}S_{\mathcal{H}}r_{(2)}r_{(1)}h_{(3)}%
\right) \\
& =\textstyle\sum a_{\left\langle 0\right\rangle }\otimes rh_{(2)}\otimes
\pi _{\mathcal{H}}\left( a_{\left\langle 1\right\rangle }S_{\mathcal{H}%
}h_{(1)}h_{(3)}\right) \\
& =r\cdot \rho _{\mathcal{A}\otimes \mathcal{H}}(a\otimes h).
\end{align*}%
Note that the third equality follows by (\ref{eq:R_H}), while in the fourth
one we used (\ref{ec:involutie}).

(3) This part is a direct application of Lemma \ref{pr:Tor=comodule}.

(4) Let $\lambda :=\beta \circ \eta ,\ $where $\beta \mathcal{\ }$is the
canonical map in the definition of Hopf-Galois extensions and $\eta $ is the
following $\mathbb{K}$-linear isomorphism
\begin{equation*}
\eta :(\mathcal{A}\otimes \mathcal{A})_{\mathcal{B}}\rightarrow \mathcal{A}%
\otimes _{\mathcal{B}}\mathcal{A},\quad \eta \left( \lbrack a\otimes x]_{%
\mathcal{B}}\right) \ =x\otimes _{\mathcal{B}}a.
\end{equation*}%
By the definition of $\beta $ and $\eta $ we can easily show that
\begin{equation*}
\lambda \left( \lbrack a\otimes x]_{\mathcal{B}}\right) =\textstyle\sum
xa_{\left\langle 0\right\rangle }\otimes a_{\left\langle 1\right\rangle }.
\end{equation*}%
We claim that $\lambda $ is an isomorphism in $_{\mathcal{R}_{\mathcal{H}}}%
\mathfrak{M}^{\mathcal{C}_{\mathcal{H}}}.$ Obviously, $\lambda $ is
bijective as $\beta $ and $\eta $ are so. If $r\in \mathcal{R}_{\mathcal{H}}$
and $a,$ $x\in \mathcal{A}$ then%
\begin{align*}
\lambda (r\cdot \lbrack a\otimes x]_{\mathcal{B}})& =\lambda \left( %
\textstyle\sum [\kappa ^{2}(r)a\otimes x\kappa ^{1}(r)]_{B}\right) \\
& =\textstyle\sum x\kappa ^{1}(r)\kappa ^{2}(r)_{\left\langle 0\right\rangle
}a_{\left\langle 0\right\rangle }\otimes \kappa ^{2}(r)_{\left\langle
1\right\rangle }a_{\left\langle 1\right\rangle } \\
& =\textstyle\sum xa_{\left\langle 0\right\rangle }\otimes ra_{\left\langle
1\right\rangle },
\end{align*}%
where for the last equality we used (\ref{eq:k}). Thus $\lambda $ is a
morphism of $\mathcal{R}_{\mathcal{H}}$-modules. Let $\rho $ denote the
coaction of $\mathcal{C}_{\mathcal{H}}$ on $(\mathcal{A}\otimes \mathcal{A}%
)_{\mathcal{B}}.$ Hence, by the definition of $\rho $ and the fact that $\pi
_{\mathcal{H}}$ is a trace map, we get%
\begin{align*}
(\lambda \otimes \mathcal{H})\circ \rho ([a\otimes x]_{\mathcal{B}}) &=%
\textstyle\sum \lambda \left( \lbrack a_{\left\langle 0\right\rangle
}\otimes x_{\left\langle 0\right\rangle }]_{\mathcal{B}}\right) \otimes \pi
_{\mathcal{H}}(a_{\left\langle 1\right\rangle }x_{\left\langle
1\right\rangle }) \\
&=\textstyle\sum x_{\left\langle 0\right\rangle }a_{\left\langle
0\right\rangle }\otimes a_{\left\langle 1\right\rangle }\otimes \pi _{%
\mathcal{H}}(x_{\left\langle 1\right\rangle }a_{\left\langle 2\right\rangle
}).
\end{align*}%
On the other hand, by (\ref{ro:AxH}), it follows%
\begin{eqnarray*}
\rho _{\mathcal{A}\otimes \mathcal{H}}\circ \lambda ([a\otimes x]_{\mathcal{B%
}}) &=&\rho _{\mathcal{A}\otimes \mathcal{H}}\left( \textstyle\sum
xa_{\left\langle 0\right\rangle }\otimes a_{\left\langle 1\right\rangle
}\right) \\
&=&\textstyle\sum (xa_{\left\langle 0\right\rangle })_{\left\langle
0\right\rangle }\otimes (a_{\left\langle 1\right\rangle })_{(2)}\otimes \pi
_{\mathcal{H}}\left( (xa_{\left\langle 0\right\rangle })_{\left\langle
1\right\rangle }S_{\mathcal{H}}(a_{\left\langle 1\right\rangle
_{(1)}})a_{\left\langle 1\right\rangle _{(3)}}\right) \\
&=&\textstyle\sum x_{\left\langle 0\right\rangle }a_{\left\langle
0\right\rangle }\otimes a_{\left\langle 3\right\rangle }\otimes \pi _{%
\mathcal{H}}(x_{\left\langle 1\right\rangle }a_{\left\langle 1\right\rangle
}S_{\mathcal{H}}a_{\left\langle 2\right\rangle }a_{\left\langle
4\right\rangle }) \\
&=&\textstyle\sum x_{\left\langle 0\right\rangle }a_{\left\langle
0\right\rangle }\otimes a_{\left\langle 1\right\rangle }\otimes \pi _{%
\mathcal{H}}(x_{\left\langle 1\right\rangle }a_{\left\langle 2\right\rangle
}).
\end{eqnarray*}%
Summarizing, the computation above shows us that $\lambda $ is a morphism of
$\mathcal{C}_{\mathcal{H}}$-comodules too. Furthermore, for a right $%
\mathcal{R}_{\mathcal{H}}$-module $N$, we get
\begin{equation*}
N\otimes _{\mathcal{R}_{\mathcal{H}}}G_{V}(U)\cong N\otimes _{\mathcal{R}_{%
\mathcal{H}}}\left[ (\mathcal{A}\otimes \mathcal{H})\square _{\mathcal{C}_{%
\mathcal{H}}}V\right] \cong \left[ N\otimes _{\mathcal{R}_{\mathcal{H}}}(%
\mathcal{A}\otimes \mathcal{H})\right] \square _{\mathcal{C}_{\mathcal{H}}}V,
\end{equation*}%
where the first isomorphism is defined by $N\otimes _{\mathcal{R}_{\mathcal{H%
}}}\left( \lambda \square _{\mathcal{C}_{\mathcal{H}}}V\right) $ and the
second one comes from the commutation of the tensor product and the cotensor
product. We have obtained a natural isomorphism
\begin{equation*}
\nu (N):N\otimes _{\mathcal{R}_{\mathcal{H}}}G_{V}(U)\rightarrow \left[
\mathcal{A}\otimes (N\otimes _{\mathcal{R}_{\mathcal{H}}}\mathcal{H})\right]
\square _{\mathcal{C}_{\mathcal{H}}}V
\end{equation*}
given, for $z:=\sum_{i=1}^{n}\left[ a^{i}\otimes b^{i}\right] _{\mathcal{B}%
}\otimes v^{i}$ in $G_{V}(U)$ and $x\in N$, by
\begin{equation*}
\nu (N)(x\otimes _{\mathcal{R}_{\mathcal{H}}}z)=\textstyle\sum_{i=1}^{k}\sum
b^{i}a_{\left\langle 0\right\rangle }^{i}\otimes (x\otimes _{\mathcal{R}_{%
\mathcal{H}}}a_{\left\langle 1\right\rangle }^{i})\otimes v^{i}.
\end{equation*}%
Let us prove that $\phi (U)$ is an isomorphism. Since $\mathcal{H}/\mathcal{R%
}_{\mathcal{H}}^{+}\mathcal{H}\cong \mathbb{K},$
\begin{equation*}
\mathbb{K}\otimes _{\mathcal{R}_{\mathcal{H}}}\mathcal{H}\cong \mathcal{R}_{%
\mathcal{H}}/\mathcal{R}_{\mathcal{H}}^{+}\otimes _{\mathcal{R}_{\mathcal{H}%
}}\mathcal{H}\cong \mathcal{H}/\mathcal{R}_{\mathcal{H}}^{+}\mathcal{H}\cong
\mathbb{K}.
\end{equation*}%
Note that this isomorphism maps $1\otimes _{\mathcal{R}_{\mathcal{H}}}h$ to $%
\varepsilon (h).$ Let $\gamma $ be the composition of the isomorphism $\left[
\mathcal{A}\otimes (\mathbb{K}\otimes _{\mathcal{R}_{\mathcal{H}}}\mathcal{H}%
)\right] \square _{\mathcal{C}_{\mathcal{H}}}V\cong $ $\mathcal{A}\square _{%
\mathcal{C}_{\mathcal{H}}}V$ and $\nu (\mathbb{K}).$ Then
\begin{equation*}
\gamma (1\otimes _{\mathcal{R}_{\mathcal{H}}}z)=\textstyle\sum_{i=1}^{k}\sum
b^{i}a^{i}\otimes v^{i},
\end{equation*}
where $z\in G_{V}(U)$ is given by the same formula as above. Furthermore,
the multiplication in $\mathcal{A}$ induces an isomorphism of $\mathcal{C}_{%
\mathcal{H}}$-comodules
\begin{equation*}
\mu:U_{\mathcal{A}}\longrightarrow \mathcal{A}, \qquad\mu ([a^{\prime
}\otimes a^{\prime \prime }]_{\mathcal{A}})=a^{\prime \prime }a^{\prime }.
\end{equation*}
It is easy to see that $\phi (U)=(\mu\; \square _{\mathcal{C}_{\mathcal{H}%
}}V)\circ \gamma$, so $\phi (U) $ is an isomorphism.

Let $P_{\ast }$ be a resolution of $\mathbb{K}$ in $\mathfrak{M}_{\mathcal{R}%
_{\mathcal{H}}}$. The natural transformation $\nu $ yields isomorphisms
\begin{equation*}
P_{\ast }\otimes _{\mathcal{R}_{\mathcal{H}}}G_{V}(U)\cong \left[ \mathcal{A}%
\otimes (P_{\ast }\otimes _{\mathcal{R}_{\mathcal{H}}}\mathcal{H})\right]
\square _{\mathcal{C}_{\mathcal{H}}}V.
\end{equation*}%
Since $\mathcal{A}\otimes (-)$ and $(-)\square _{\mathcal{C}_{\mathcal{H}}}V$
are exact functors it follows%
\begin{equation*}
\mathrm{Tor}_{n}^{\mathcal{R}_{\mathcal{H}}}(\mathbb{K},G_{V}(U))\simeq
\lbrack \mathcal{A}\otimes \mathrm{Tor}_{n}^{\mathcal{R}_{\mathcal{H}}}(%
\mathbb{K},\mathcal{H})]\square _{\mathcal{C}_{\mathcal{H}}}V.
\end{equation*}%
Hence the lemma is completely proven as $\mathrm{Tor}_{n}^{\mathcal{R}_{%
\mathcal{H}}}(\mathbb{K},\mathcal{H})=0$, for every $n>0.$
\end{proof}

\begin{definition}
We say that a Hopf algebra $\mathcal{H}$ \emph{has enough cocommutative}
\emph{elements} if $\mathcal{R}_{\mathcal{H}}^{+}\mathcal{H}=\mathcal{H}%
^{+}. $
\end{definition}

\begin{theorem}
\label{te:sir}Let $\mathcal{H}$ be a Hopf algebra such that $S_{\mathcal{H}%
}^{2}=\mathrm{Id}_{\mathcal{H}}$. We assume that $\mathcal{H}$ has enough
cocommutative elements and $\mathrm{Tor}_{\ast }^{\mathcal{R}_{\mathcal{H}}}(%
\mathbb{K},\mathcal{H})=0.$ If ${\mathcal{B}\subseteq \mathcal{A}}$ is a
faithfully flat $\mathcal{H}$-Galois extension and $V$ is an injective left $%
\mathcal{C}_{\mathcal{H}}$-comodule then, for every Hopf bimodule $M$, there
is a spectral sequence in the category $_{\mathcal{Z}_{0}}\mathfrak{M}$
\begin{equation}
\mathrm{Tor}_{p}^{\mathcal{R}_{\mathcal{H}}}(\mathbb{K},\mathrm{HH}_{q}(%
\mathcal{B},M\square _{\mathcal{C}_{\mathcal{H}}}V))\Longrightarrow \mathrm{%
HH}_{p+q}(\mathcal{A},M)\square _{\mathcal{C}_{\mathcal{H}}}V.
\label{eq:sir}
\end{equation}
\end{theorem}

\begin{proof}
We know that $U:=\mathcal{A}\otimes \mathcal{A}$ is a generator in $_{%
\mathcal{A}}\mathfrak{M}_{\mathcal{A}}^{\mathcal{H}}.$ In view of Lemma \ref%
{le: fi}, one can apply Proposition \ref{pr:sir spectral} to the
following categories:
\begin{equation*}
\mathfrak{A}:={}_{\mathcal{A}}\mathfrak{M}_{\mathcal{A}}^{\mathcal{H}%
},\qquad \mathfrak{B}:=_{\mathcal{R}_{\mathcal{H}}\otimes \mathcal{Z}_{0}}%
\mathfrak{M},\qquad \mathfrak{C}:={}_{\mathcal{Z}_{0}}\mathfrak{M}.
\end{equation*}%
The functors $F,$ $G_{V}$ and $H$ are given by
\begin{equation*}
F:=\mathbb{K}\otimes _{\mathcal{R}_{\mathcal{H}}}(-),\qquad G_{V}:=(-)_{%
\mathcal{B}}\square _{\mathcal{C}_{\mathcal{H}}}V,\qquad H_{V}:=(-)_{%
\mathcal{A}}\square _{\mathcal{C}_{\mathcal{H}}}V
\end{equation*}%
and the natural transformation $\phi :F\circ G_{V}\rightarrow H_{V}$ is
defined in Lemma \ref{le: fi} (1). To compute the left derived functors of $%
G_{V}$ and $H_{V}$ we use Propositions \ref{pr:G_V} and \ref{pr:H_V}. Since
any projective $\mathcal{R}_{\mathcal{H}}\otimes \mathcal{Z}_{0}$-module is
also projective as an $\mathcal{R}_{\mathcal{H}}$-module, it follows that $%
\mathrm{L}_{n}F\cong \mathrm{Tor}_{n}^{\mathcal{R}_{\mathcal{H}}}(\mathbb{K},%
\mathbb{-)}.$
\end{proof}

\begin{proposition}
\label{pr:R_H-semisimple}Let $\mathcal{H}$ be a finite-dimensional
Hopf
algebra over a field $\mathbb{K}$ of characteristic zero such that $S_{%
\mathcal{H}}^{2}=\mathrm{Id}_{\mathcal{H}}$. Then $\mathcal{R}_{\mathcal{H}}$
is semisimple and $\mathcal{C}_{\mathcal{H}}$ is cosemisimple.
\end{proposition}

\begin{proof}
We first prove that $\mathcal{R}_{\mathcal{H}}$ is semisimple in the case
when $\mathbb{K}$ is algebraically closed. By Larson-Radford Theorem \cite[%
Theorem 7.4.6]{DNR}, it follows that $\mathcal{H}$ is semisimple and
cosemisimple. We claim that, in this particular case, $\mathcal{R}_{\mathcal{%
H}^{\ast}}$ equals the $\mathbb{K}$-subalgebra $C_{\mathbb{K}}(\mathcal{H)%
}$ of $\mathcal{H}^{\ast}$, which is generated by the set of characters of $%
\mathcal{H}.$ For the definition and properties of characters of a
semisimple Hopf algebra, the reader is referred to \cite[Section 7.5]{DNR}.
Recall that an element $\alpha \in \mathcal{H}^{\ast }$ is said to be a
trace map on $\mathcal{H}$ if and only if $\alpha \mathcal{\ }$vanishes on
the space of commutators $[\mathcal{H},\mathcal{H}].$ Let us show that $%
\mathcal{R}_{\mathcal{H}^{\ast }}$ equals the space of all trace maps on $%
\mathcal{H}.\ $By the definition of comultiplication of $\mathcal{H}^{\ast
}, $%
\begin{equation*}
\Delta (\alpha )=\sum \alpha _{(1)}\otimes \alpha _{(2)}
\end{equation*}%
if and only if $\alpha (xy)=\sum \alpha _{(1)}(x)\alpha _{(2)}(y),$ for
all $x,y\in \mathcal{H}.$ On the other hand, $\alpha \in \mathcal{R}_{%
\mathcal{H}^{\ast }}$ if and only if $\Delta (\alpha )=\sum \alpha
_{(2)}\otimes \alpha _{(1)}.$ Therefore, for $\alpha \in \mathcal{R}_{%
\mathcal{H}^{\ast }}$, we get%
\begin{equation*}
\alpha (xy)=\sum \alpha _{(2)}(x)\alpha _{(1)}(y)=\sum \alpha
_{(1)}(y)\alpha _{(2)}(x)=\alpha (yx),
\end{equation*}%
so $\alpha $ is a trace map. The other implication can be proved similarly.
We can now show that $\mathcal{R}_{\mathcal{H}^{\ast }}=C_{\mathbb{K}}(%
\mathcal{H)}.$ By definition, a character is a trace map, so $C_{\mathbb{K}}(%
\mathcal{H)}$ is a subspace of $\mathcal{R}_{\mathcal{H}^{\ast }}.$
Therefore, it is enough to show that $\dim \mathcal{R}_{\mathcal{H}^{\ast
}}\leq \dim C_{ \mathbb{K}}(\mathcal{H)}$. As the base field is
algebraically closed, $\mathcal{H\cong }\prod_{i=1}^{n}M_{d_{i}}(\mathbb{K}%
). $ For every $i=1,\dots ,n,$ let $V_{i}$ be a simple left
$\mathcal{H}$-module associated to the block $M_{d_{i}}(\mathbb{K})$
and let $\chi _{i}$
denote the irreducible character corresponding to $V_{i}.$ By \cite[%
Proposition 7.5.7]{DNR}, $\chi _{1},\dots ,\chi _{n}$ are linearly
independent over $\mathbb{K}$, as elements in $\mathcal{H}^{\ast }.$ On the
other hand, using the canonical basis $\{E_{ip_{i}q_{i}}\mid i=1,\dots
,n,~p_{i},q_{i}=1,\dots ,d_{i}\}$ on $\prod_{i=1}^{n}M_{d_{i}}(\mathbb{K})$,
one can show that $\alpha $ is a trace map if and only if there are $%
a_{1},\dots ,a_{n}$ in $\mathbb{K}$ such that
\begin{equation*}
\alpha (E_{ip_{i}q_{i}})=\left\{
\begin{array}{cc}
0, & \text{if }p_{i}\neq q_{i}, \\
a_{i}, & \text{if }p_{i}=q_{i}.%
\end{array}%
\right.
\end{equation*}%
Hence, $\dim \mathcal{R}_{\mathcal{H}^{\ast }}=n\leq \dim C_{\mathbb{K}}(%
\mathcal{H)}$. To deduce that $\mathcal{R}_{\mathcal{H}^{\ast }}$ is
semisimple, we now use \cite[Theorem 7.5.12]{DNR} and the fact that $C_{%
\mathbb{K}}(\mathcal{H})=\mathbb{K}\otimes _{\mathbb{Q}}C_{\mathbb{Q}}(%
\mathcal{H)}$. We have already remarked that $\mathcal{H}$ is
cosemisimple too. Thus, $\mathcal{R}_{\mathcal{H}}\cong $
$\mathcal{R}_{\mathcal{H}^{\ast \ast }}$ is also semisimple.

We now assume that $\mathbb{K}$ is an arbitrary field of characteristic
zero. Let $\overline{\mathbb{K}}$ be an algebraic closure of $\mathbb{K}$
and let $\overline{\mathcal{H}}:=\overline{\mathbb{K}}\otimes _{\mathbb{K}}%
\mathcal{H}.$ We claim that $\overline{\mathbb{K}}\otimes _{\mathbb{K}}%
\mathcal{R}_{\mathcal{H}}=\mathcal{R}_{\overline{\mathcal{H}}}.$ Indeed, let
$\{\alpha _{i}\mid i\in I\}$ be a basis of $\overline{\mathbb{K}}$ as a $%
\mathbb{K}$-vector space and $z=\sum_{i\in I}\alpha ^{i}\otimes h^{i}\in
\overline{\mathcal{H}}$. By the definition of the comultiplication of $%
\overline{\mathcal{H}},$ $z$ belongs to $\mathcal{R}_{\overline{\mathcal{H}}%
} $ if and only if
\begin{equation*}
\textstyle\sum_{i\in I}(\alpha ^{i}\otimes h_{(1)}^{i})\otimes _{\overline{%
\mathbb{K}}}(1\otimes h_{(2)}^{i})=\sum_{i\in I}(1\otimes
h_{(2)}^{i})\otimes _{\overline{\mathbb{K}}}(\alpha ^{i}\otimes h_{(1)}^{i}).
\end{equation*}%
Thus $\sum_{i\in I}\alpha ^{i}\otimes h_{(1)}^{i}\otimes
h_{(2)}^{i}=\sum_{i\in I}\alpha ^{i}\otimes h_{(2)}^{i}\otimes h_{(1)}^{i}.$
Since the elements $\alpha _{i}$ are linearly independent over $\mathbb{K}$
it results that $z\in \mathcal{R}_{\overline{\mathcal{H}}}$ if and only if
each $h_{i}$ is an element in $\mathcal{R}_{\mathcal{H}}.$ Consequently, the
claimed equality is proven.

Obviously, $S_{\overline{\mathcal{H}}}^{2}=\mathrm{Id}_{\overline{\mathcal{H}%
}}.$ Since $\overline{\mathcal{H}}$ is a Hopf algebra over an
algebraically closed field, it follows that
$\mathcal{R}_{\overline{\mathcal{H}}}$ is semisimple. Let $J$ be the
Jacobson radical of $\mathcal{R}_{\mathcal{H}}$ which is a
finite-dimensional algebra. Thus $J$ is a nilpotent ideal. Clearly
$\overline{\mathbb{K}}\otimes _{\mathbb{K}}J$ is a nilpotent ideal
in $\overline{\mathbb{K}}\otimes
_{\mathbb{K}}\mathcal{R}_{\mathcal{H}}\cong
\mathcal{R}_{\overline{\mathcal{H}}},$ so it is contained in the
Jacobson radical of $\mathcal{R}_{\overline{\mathcal{H}}}.$ We
deduce that $\overline{
\mathbb{K}}\otimes _{\mathbb{K}}J=0.$ Thus $J=0,$ so $\mathcal{R}_{\mathcal{%
H }}$ is semisimple, being finite-dimensional.

It remains to prove that $\mathcal{C}_{\mathcal{H}}$ is cosemisimple. The
dual algebra $\mathcal{C}_{\mathcal{H}}^{\ast }$ is isomorphic to the
subalgebra of trace maps on $\mathcal{H}$, i.e. $\mathcal{C}_{\mathcal{H}%
}^{\ast }\cong \mathcal{R}_{ \mathcal{H}^{\ast }}.$ As $\mathcal{H}$ is
cosemisimple, we have already seen that $\mathcal{C}_{\mathcal{H}}^{\ast
}\cong\mathcal{R}_{ \mathcal{H}^{\ast }}$ is semisimple. Hence $\mathcal{C}_{%
\mathcal{H}}$ is cosemisimple.
\end{proof}

\begin{theorem}
Let ${\mathcal{B}\subseteq \mathcal{A}}$ be an $\mathcal{H}$-Galois
extension, where $\mathcal{H}$ is a Hopf algebra
of finite dimension over a field of characteristic zero. If $\mathcal{H}$ has enough
cocommutative elements and
$S_{\mathcal{H}}^{2}=\mathrm{Id}_{\mathcal{H}}$ then, for every Hopf
bimodule $M$ and every left $\mathcal{C}_{\mathcal{H}}$-comodule
$V$, there is  isomorphisms of ${\mathcal{Z}_{0}}$-modules
\begin{equation*}
\mathbb{K}\otimes _{\mathcal{R}_{\mathcal{H}}}\mathrm{HH}_{n}(\mathcal{B}%
,M\square _{C_{\mathcal{H}}}V)\cong \mathrm{HH}_{n}(\mathcal{A},M)\square
_{C_{\mathcal{H}}}V.
\end{equation*}
\end{theorem}

\begin{proof}
By the proof of the previous proposition, $\mathcal{H}$ is cosemisimple.
Thus $\mathcal{A}$ is injective as an $\mathcal{H}$-comodule, so the
extension ${\mathcal{B}\subseteq \mathcal{A}}$ is faithfully flat, cf. \cite[%
Theorem 4.10]{SS}. In view of the same proposition $V$ is injective, as $%
\mathcal{C}_{\mathcal{H}}$ is cosemisimple, and $\mathbb{K}$ is projective
as a right $\mathcal{R}_{\mathcal{H}}$-module. Therefore, under the
assumptions of the theorem, the spectral sequence (\ref{eq:sir}) exists and
collapses. The edge maps of this spectral sequence yields the required
isomorphisms.
\end{proof}

For another application of Theorem \ref{te:sir}, let us take the Hopf
algebra $\mathcal{H}$ to be cocommutative. In this case $\mathcal{R}_{%
\mathcal{H}}=\mathcal{H},$ so the assumptions on $\mathcal{H}$ are trivially
satisfied. We obtain the spectral sequence from the following corollary.
Note that a related result can be found in \cite[Theorem 3.1]{St2}, where
the extension $\mathcal{B}\subseteq \mathcal{A}$ is not necessarily
faithfully flat but $V$ is just a subcoalgebra of $\mathcal{C}_{\mathcal{H}%
} $ which is injective in $^{\mathcal{C}_{\mathcal{H}}}\mathfrak{M.}$

\begin{corollary}
Let $\mathcal{B}\subseteq \mathcal{A}$ be a faithfully flat $\mathcal{H}$%
-Galois extension, with $\mathcal{H}$ a cocommutative Hopf algebra. If $V$
is an injective right $\mathcal{C}_{\mathcal{H}}$-comodule and $M$ is a Hopf
bimodule then there exists a spectral sequence in the category $_{\mathcal{Z%
}_{0}}\mathfrak{M}$%
\begin{equation*}
\mathrm{Tor}_{p}^{\mathcal{H}}(\mathbb{K},\mathrm{HH}_{q}(\mathcal{B}%
,M\square _{\mathcal{C}_{\mathcal{H}}}V))\Longrightarrow \mathrm{HH}_{{p+q}}(%
\mathcal{A},M)\square _{\mathcal{C}_{\mathcal{H}}}V.
\end{equation*}
\end{corollary}

\begin{remark}
By taking $V:=\mathcal{C}_\mathcal{H}$ in the above corollary we obtain
(only in the case of cocommutative Hopf algebras) the spectral sequence \cite%
[Theorem 4.5]{St1}.
\end{remark}

Let now consider the case when the Hopf algebra $\mathcal{H}$ is the group
algebra $\mathbb{K}G$ of an arbitrary group $G$. By \cite[Theorem 8.1.7]{Mo}%
, an extension $\mathcal{B}\subseteq \mathcal{A}$ is $\mathbb{K}G$-Galois if
and only if $\mathcal{A}$ is $G$-strongly graded and $\mathcal{B}$ is its
homogeneous component of degree one, i.e. $\mathcal{A}$ is a direct sum of
linear subspaces $\mathcal{A}=\bigoplus_{g\in G}\mathcal{A}_{g}$ such that $%
\mathcal{A}_{1}=\mathcal{B}$ and
\begin{equation*}
\mathcal{A}_{g}\mathcal{A}_{h}=\mathcal{A}_{gh},
\end{equation*}%
for any $g, h$ in $G$. A strongly graded algebra, i.e. a $\mathbb{K}G$%
-Galois extension $\mathcal{B}\subseteq \mathcal{A}$, is always faithfully
flat, as $\mathbb{K}G$ is cosemisimple. Furthermore, the coalgebra $\mathcal{C%
}_{\mathbb{K}G}$ is cosemisimple and pointed, cf. \cite[Exemple 1.2 (a)]{St2}%
. A direct application of the preceding corollary, for $V:=\mathcal{C}_{%
\mathbb{K}G}$, yields \cite[Theorem 2.5 (a)]{Lo1}.

Furthermore, a $\mathbb{K}G$-comodule is a vector space $M$ together with a
decomposition as a direct sum of subspaces
\begin{equation*}
M=\textstyle\bigoplus_{g\in G}M_{g}.
\end{equation*}
Hence, an $\mathcal{A}$-bimodule $M$ is a Hopf bimodule if and only if the above decomposition satisfies, for any $g$ and $h$ in
$G$, the following relations
\begin{equation*}
\mathcal{A}_{h}M_{g}\subseteq M_{hg}\qquad \text{and\qquad }M_{g}\mathcal{A}%
_{h}\subseteq M_{gh}.
\end{equation*}%
Thus, to give a Hopf bimodule is equivalent to give a $G$-graded $\mathcal{A}$%
-bimodule.

The coalgebra $\mathcal{C}_{\mathbb{K}G}$ is cosemisimple and pointed, cf.
\cite[Exemple 1.2 (a)]{St2}. Recall that on $\mathcal{C}_{\mathbb{K}G}$
there is a canonical basis $\{e_{\sigma }\mid \sigma \in T(G)\}$, where $%
T(G) $ is the set of conjugacy classes in $G$ and each $e_{\sigma }$ is a
group-like element. A left (or right) $\mathcal{C}_{\mathbb{K}G}$ comodule
structure $\rho_V:V\longrightarrow V\otimes \mathbb{K}G$ on a given vector
space $V$ is uniquely defined by a decomposition of $V$ as a direct sum
\begin{equation*}
V=\textstyle\bigoplus_{\sigma\in T(G)}V_{\sigma}.
\end{equation*}
Note that the subspace $V_\sigma$ is given by
\begin{equation*}
V_\sigma=\{v\in V\mid\rho_{V}(v)=e_{\sigma }\otimes v\}.
\end{equation*}
We shall say that $V_\sigma$ is the homogeneous component of $V$ of degree $%
\sigma$.

We now fix a conjugacy class $\sigma$ in $G$ and we put $V:=\mathbb{K}%
e_\sigma$. Since $e_\sigma$ is a group-like element, $V$ is a left $\mathcal{%
C}_{\mathbb{K}G}$-subcomodule of $\mathcal{C}_{\mathbb{K}G}$ and $V_\sigma=V$%
. Let $M$ be a $G$-graded bimodule and $M_\sigma:=\bigoplus_{g\in\sigma}M_g$%
. Thus
\begin{equation*}
M\square _{\mathcal{C}_{\mathcal{H}}}V=M_\sigma.
\end{equation*}
Therefore, if $G_{V} $ and $H_V$ are the functors in the proof of Theorem %
\ref{te:sir} then
\begin{equation*}
G_{V}(M)\cong \mathcal{B}\otimes _{\mathcal{B}^{e}}M_{\sigma },\qquad \qquad
H_V(M):=\left(M_{\mathcal{A}}\right)_{\sigma }.
\end{equation*}%
Let us notice that $M_\mathcal{A}$ is a right $\mathcal{C}_{\mathbb{K}G}$%
-comodule, so it makes sense to speak about $(M_\mathcal{A})_\sigma$. More
generally, the coaction of $\mathcal{C}_{\mathbb{K}G}$ on the Hochschild
homology of $\mathcal{A}$ with coefficients in $M$ induces a decomposition
of $\mathrm{HH}_\ast(\mathcal{A},M)$ as a direct sum of its homogeneous
components $\mathrm{HH}_\ast(\mathcal{A},M)_\sigma$, for $\sigma$ arbitrary
in $T(G)$.

As an application of Theorem \ref{te:sir} we now get the following result,
that was also proved in \cite[Theorem 2.5 (b)]{Lo1} by a different method.

\begin{corollary}
If ${\mathcal{B}\subseteq \mathcal{A}}$ is a strongly $G$-graded algebra
then, for any graded $\mathcal{A}$-bimodule $M$ and $\sigma\in T(G),$ there
is a natural spectral sequence in $_{\mathcal{Z}_{0}}\mathfrak{M}$
\begin{equation}  \label{sir-graded}
\mathrm{H}_{p}(G,\mathrm{\mathrm{HH}}_{q}(\mathcal{B},M_{\sigma}))%
\Longrightarrow \mathrm{\mathrm{HH}}_{p+q}(\mathcal{A},M)_{\sigma}.
\end{equation}
\end{corollary}

\begin{proof}
Group homology $\mathrm{H}_{\ast }(G,X)$ and $\mathrm{Tor}_{\ast }^{\mathbb{K%
}G}(\mathbb{K},X)$ are equal for any $\mathbb{K}G$-module $X$, cf. \cite[%
Theorem 3.6.2]{We}.
\end{proof}

\begin{remark}
We keep the notation from the previous corollary. Let us pick up an element $%
g$ in $\sigma$ and denote the centralizer of $g$ in $G$ by $C_G(g)$. One can
show that
\begin{equation*}
\mathrm{HH}_{q}(\mathcal{B},M_{\sigma}))\cong\mathbb{K}G\otimes_{\mathbb{K}%
C_G(g)}\mathrm{HH}_{q}(\mathcal{B},M_{g}).
\end{equation*}
Thus, by Shapiro's Lemma, the terms in the second page of the spectral
sequence (\ref{sir-graded}) are isomorphic to
\begin{equation*}
E^2_{p,q}=\mathrm{H}_{p}(C_G(g),\mathrm{\mathrm{HH}}_{q}(\mathcal{B},M_{g})).
\end{equation*}
For details the reader is referred to \cite[p. 504]{Lo1}.
\end{remark}

Now we are going to investigate the case when $\mathcal{H}$ is
commutative
but not necessarily cocommutative. Thus $\mathcal{C}_{\mathcal{H}}=\mathcal{%
H.}$ Our first aim is to show that we can drop the assumptions on $\mathcal{R%
}_{\mathcal{H}}$ in Theorem \ref{te:sir}. To this end, we need the
following.

\begin{lemma}
\label{KG^ast} \label{le:KG*}Let $G$ be a finite group. If $\mathcal{H}:=(%
\mathbb{K}G)^{\ast }$ then $\mathcal{H}$ has enough cocommutative
elements and $\mathcal{R}_{\mathcal{H}}$ is semisimple.
\end{lemma}

\begin{proof}
Let $\{p_{x}\mid x\in G\}$ be the dual basis of the canonical basis on $%
\mathbb{K}G.$ By definition, the coalgebra structure on $\mathcal{H}$ is
given by
\begin{equation*}
\Delta (p_{x})\textstyle=\sum_{g\in G}p_{xg^{-1}}\otimes p_{g}\qquad \text{%
and\qquad }\varepsilon (p_{x})=\delta _{x,1}.
\end{equation*}%
Thus, an element $z=\sum_{x\in G}a_{x}p_{x}$ belongs to $\mathcal{R}_{%
\mathcal{H}}$ if and only if
\begin{equation*}
\textstyle\sum_{x,g\in G}a_{x}p_{xg^{-1}}\otimes p_{g}=\sum_{x,g\in
G}a_{x}p_{g}\otimes p_{xg^{-1}}.
\end{equation*}%
Since $\{p_{x}\otimes p_{y}\mid x,y\in G\}$ is a basis on $\mathcal{H}%
\otimes \mathcal{H,}$ we deduce that $a_{xg}=a_{gx}$. Therefore, if $\sigma
\in T(G)$ then there is $a_{\sigma }$ in $\mathbb{K}$ such that $%
a_{g}=a_{\sigma },$ for all $g\in\sigma$. It follows%
\begin{equation*}
\textstyle z=\sum_{\sigma \in T(G)}a_{\sigma }p_{\sigma },
\end{equation*}%
where $p_{\sigma }:=\sum_{x\in \sigma }p_{x}$. In conclusion, $\mathcal{R}_{%
\mathcal{H}}$ is the $\mathbb{K}$-linear subspace generated by all $%
p_{\sigma },$ with $\sigma \in T(G).$ On the other hand, $p_{x}p_{y}=\delta
_{x,y}p_{x},$ for arbitrary $x,y\in G.$ Thus, if $x\in G$ and $x\neq 1$,
then
\begin{equation*}
p_{x}=p_{\sigma}p_{x},
\end{equation*}
where $\sigma$ denotes the conjugacy class of $x.$ Hence the relation $%
\mathcal{R}_{\mathcal{H}}^{+}\mathcal{H=H}^{+}$ is proven. To conclude the
proof we remark that
\begin{equation*}
\textstyle p_{\sigma }p_{\tau }=\delta _{\sigma ,\tau }p_{\sigma }\qquad
\text{and\qquad }\sum_{\sigma \in T(G)}p_{\sigma }=1_{\mathcal{R}_{\mathcal{H%
}}},
\end{equation*}%
so $\mathcal{R}_{\mathcal{H}}$ is semisimple. In fact, the above relations
shows us that $\mathcal{R}_{\mathcal{H}}\simeq \mathbb{K}^{\# T(G)}.$
\end{proof}

\begin{proposition}
\label{pr:ss} Let $\mathcal{H}$ be a commutative Hopf algebra. If $\mathcal{H%
}$ is semisimple then $\mathcal{H}$ has enough cocommutative elements and $%
\mathcal{R}_{ \mathcal{H}}$ is a semisimple $\mathbb{K}$-algebra.
\end{proposition}

\begin{proof}
Let us assume first that $\mathcal{H}$ is semisimple. Then
$\mathcal{H}$ is finite-dimensional by \cite[Remark 3.8(b)]{St1}. Let
$\overline{\mathbb{K}}$
be an algebraic closure of $\mathbb{K}$ and let $\overline{\mathcal{H}}:=%
\overline{\mathbb{K}}\otimes _{\mathbb{K}} \mathcal{H}.$ Hence $\overline{%
\mathcal{H}}$ is a finite-dimensional commutative Hopf algebra over $%
\overline{\mathbb{K}}.$ Since $\mathcal{H}$ is semisimple it follows that $%
\overline{\mathcal{H}}$ is semisimple too. Thus, the dual Hopf algebra $%
\overline{\mathcal{H}}^{\ast }$ is cosemisimple and cocommutative. Since $%
\overline{\mathbb{K}}$ is algebraically closed, there is a finite group $G$
such that $\overline{\mathcal{H}}=(\overline{ \mathbb{K}}G)^{\ast }.$ By the
previous lemma, $\mathcal{R}_{\overline{ \mathcal{H}}}^{+}\overline{\mathcal{%
H}}=\overline{\mathcal{H}}^{+}$ and $\mathcal{R}_{\overline{\mathcal{H}}}$
is semisimple. As
\begin{equation*}
\overline{\mathbb{K}}\otimes _{\mathbb{K}}(\mathcal{H}^{+}/\mathcal{R}_{
\mathcal{H}}^{+}\mathcal{H})\cong (\overline{\mathbb{K}}\otimes _{\mathbb{K}%
} \mathcal{H}^{+})/(\overline{\mathbb{K}}\otimes _{\mathbb{K}}\mathcal{R}_{
\mathcal{H}}^{+}\mathcal{H})\cong \overline{\mathcal{H}}^{+}/\mathcal{R}_{
\overline{\mathcal{H}}}^{+}\overline{\mathcal{H}}=0
\end{equation*}%
we get $\mathcal{R}_{\mathcal{H}}^{+}\mathcal{H=H}^{+}.$ To prove that $%
\mathcal{R}_{\mathcal{H}}$ is semisimple, we can proceed as in the proof of
Proposition \ref{pr:R_H-semisimple}.

If $\mathcal{H}$ is finite-dimensional over a field of
characteristic zero then it is semisimple (and cosemisimple). Hence
we can apply the first part of the proposition.
\end{proof}

\begin{remark}
For a commutative Hopf algebra $\mathcal{H}$ over a field
$\mathbb{K}$, we have
$S_{\mathcal{H}}^{2}=\mathrm{Id}_{\mathcal{H}}$. If $\mathcal{H}$ is
finite-dimensional then the trace of $S_{ \mathcal{H}}^{2}$ equals
$(\mathrm{dim}\,\mathcal{H})1_{\mathbb{K}}$. Therefore, by
\cite[Theorem 7.4.1]{DNR}, $\mathcal{H}$ is semisimple and
cosemisimple if and only if $\dim\mathcal{H}$ is not zero in
$\mathbb{K}$. In this case, $\mathcal{H}$ has enough cocommutative
elements.
\end{remark}

\begin{theorem}
\label{te:iso}Let ${\mathcal{B}\subseteq \mathcal{A}}$ be an $\mathcal{H}$%
-Galois extension, where $\mathcal{H}$ is a
commutative Hopf algebra of finite dimension over a field $\mathbb{K}$ such that
$\dim\mathcal{H}$ is not zero in $\mathbb{K}$. If $V$ is a left
$\mathcal{H}
$-comodule and $M$ is a Hopf bimodule then there is an isomorphisms of ${%
\mathcal{Z}_{0}}$-modules
\begin{equation}
\mathbb{K}\otimes _{\mathcal{R}_{\mathcal{H}}}\mathrm{HH}_{n}(\mathcal{B}%
,M\square _{\mathcal{H}}V)\cong \mathrm{HH}_{n}(\mathcal{A},M)\square _{%
\mathcal{H}}V.  \label{ec:isom}
\end{equation}
\end{theorem}

\begin{proof}
In view of the above remark, $\mathcal{H}$ has enough cocommutative
elements and $\mathcal{H}$ is semisimple and cosemisimple. Thus $V$
is injective and ${\mathcal{B}\subseteq \mathcal{A}}$ is a
faithfully flat
extension. Since $\mathcal{R}_{\mathcal{H}}$ is semisimple, $\mathrm{Tor}%
_{p}^{\mathcal{R}_{\mathcal{H}}}(\mathbb{K},\mathcal{H})=0$ for
$p>0$, so we can apply Theorem \ref{te:sir}. Furthermore, for $p>0$
\begin{equation*}
\mathrm{Tor}_{p}^{\mathcal{R}_{\mathcal{H}}}(\mathbb{K},\mathrm{HH}_{q}(%
\mathcal{B},M\square _{\mathcal{H}}V))=0,
\end{equation*}%
as any $\mathcal{R}_{\mathcal{H}}$-module is projective. It follows that the
spectral sequence in Theorem \ref{te:sir} collapses, its edge maps giving
the isomorphism in (\ref{ec:isom}). Obviously these maps are ${\mathcal{Z}%
_{0}}$-linear, as the spectral sequence lives in $_{\mathcal{Z}_{0}}%
\mathfrak{M}$ by construction.
\end{proof}

\begin{noname}
\label{nn:G-Galois}Recall that if $G$ is a finite group of algebra
automorphisms of $\mathcal{A}$ and $\mathcal{B}=\mathcal{A}^{G}$ then $%
\mathcal{A}$ is a $(\mathbb{K}G)^{\ast }$-comodule algebra and $\mathcal{B}:=%
\mathcal{A}^{\mathrm{co}(\mathbb{K}G)^{\ast }}.$ Note that the corresponding
coaction $\rho :\mathcal{A}\longrightarrow \mathcal{A}\otimes (\mathbb{K}%
G)^{\ast }\mathcal{\ }$satisfies the relation%
\begin{equation*}
\textstyle\rho (a)=\sum_{x\in G}x(a)\otimes p_{x},
\end{equation*}%
where $\{p_{x}\mid x\in G\}$ is the dual basis of $\{x\mid x\in G\}\subseteq
\mathbb{K}G.$ It is not difficult to see that $\mathcal{B}\subseteq \mathcal{%
A}$ is $(\mathbb{K}G)^{\ast }$-Galois if and only if there are elements $%
a_{1}^{\prime },\ldots ,a_{n}^{\prime }$ and $a_{1}^{\prime \prime },\ldots
,a_{n}^{\prime \prime }$ in $\mathcal{A}$ such that
\begin{equation*}
\textstyle\sum_{i=1}^{n}a_{i}^{\prime }g(a_{i}^{\prime \prime })=\delta
_{g,1},
\end{equation*}%
for all $g\in G.$ Thus $(\mathbb{K}G)^{\ast }$-Galois extensions
generalize Galois extensions of commutative rings. For the
definition of Galois extension of commutative rings, the reader is
referred to \cite[Chapter III]{MI}. More particularly, a finite field extension is $(\mathbb{K}G)^{\ast }$%
-Galois if and only if it is separable and normal. In this case, the Galois
group of the extension is $G,$ cf. \cite[Example 6.4.3 (1)]{DNR}. For this
reason, in this paper, $(\mathbb{K}G)^{\ast }$-Galois extensions will be
called (classical) $G$-Galois extensions.

Note that an object in $_{\mathcal{A}}\mathfrak{M}_{\mathcal{A}}^{(\mathbb{K}
G)^{\ast }}$ is an $\mathcal{A}$-bimodule $M$ together with a $G$-action on $%
M$ such that, for $g\in G,~a\in \mathcal{A\ }$and $m\in M$
\begin{equation*}
g\cdot (am)=g(a)\left[ g\cdot m\right] \qquad \text{and\qquad }g\cdot (ma)= %
\left[ g\cdot m\right] g(a).
\end{equation*}%
We shall say that such an $M$ is a $(G,\mathcal{A})$-bimodule.
\end{noname}

\begin{corollary}
\label{co:iso}Let $\mathcal{B}\subseteq \mathcal{A}$ be a
$G$-Galois extension over a field $\mathbb{K}$ such that the order
of $G$ is not zero in
$\mathbb{K}$. If $M$ is a $(G,\mathcal{A})$%
-bimodule then
\begin{equation*}
\mathrm{HH}_{n}(\mathcal{A},M)^{G}\cong p_{1}\cdot \mathrm{HH}_{n}(\mathcal{B%
},M^{G})
\end{equation*}%
as ${\mathcal{Z}_{0}}$-modules, where $\{p_{x}\mid x\in G\}$ is the dual
basis of the canonical basis on $\mathbb{K}G$.
\end{corollary}

\begin{proof}
By Theorem \ref{te:iso},
\begin{equation*}
\mathbb{K}\otimes _{\mathcal{R}_{(\mathbb{K}G)^{\ast }}}\mathrm{HH}_{n}(%
\mathcal{B},M\square _{(\mathbb{K}G)^{\ast }}\mathbb{K})\cong \mathrm{HH}%
_{n}(\mathcal{A},M)\square _{(\mathbb{K}G)^{\ast }}\mathbb{K}.
\end{equation*}%
On the other hand, by the proof of Lemma \ref{KG^ast}, we get $\mathcal{R}_{(%
{\mathbb{K}G})^{\ast }}\cong \mathbb{K}^{\#T(G)}$, as
\begin{equation*}
\mathcal{S}:=\{p_{\sigma }\mid \sigma \in T(G)\}
\end{equation*}%
is a basis on $({\mathbb{K}G})^{\ast }$ and a complete set of orthogonal
idempotents. Therefore, for a $({\mathbb{K}G})^{\ast }$-module $W$, we have $%
W=\bigoplus_{\sigma \in T(G)}p_{\sigma }\cdot W$. Clearly,
\begin{equation*}
\mathbb{K}\otimes _{\mathcal{R}_{(\mathbb{K}G)^{\ast }}}W\cong \left(
\mathcal{R}_{(\mathbb{K}G)^{\ast }}/\mathcal{R}_{(\mathbb{K}G)^{\ast
}}^{+}\right) \otimes _{\mathcal{R}_{(\mathbb{K}G)^{\ast }}}W\cong W/\left(
\mathcal{R}_{(\mathbb{K}G)^{\ast }}^{+}W\right) \cong p_{1}\cdot W.
\end{equation*}%
Note that for the last isomorphism we used that $\mathcal{R}_{(\mathbb{K}%
G)^{\ast }}^{+}$ is spanned by $\mathcal{S}\setminus \{p_{1}\}.$

We conclude the proof in view of the foregoing remarks and of the
isomorphisms
\begin{equation*}
X\square _{(\mathbb{K}G)^{\ast }}V\cong X^{\mathrm{co}{(\mathbb{K}G)^{\ast }}%
}\cong X^{G}.
\end{equation*}%
In the above identifications, for a right $(\mathbb{K}G)^{\ast }$-comodule $%
X $, the $G$-invariants are taken with respect to the left $G$-action on $X$
that corresponds to the $(\mathbb{K}G)^{\ast }$-comodule structure on $X$
via the isomorphism of categories $\mathfrak{M}^{(\mathbb{K}G)^{\ast }}\cong
{}_{\mathbb{K}G}\mathfrak{M.}$
\end{proof}

\section{Centrally Hopf-Galois extensions}

Throughout this section we fix a commutative Hopf algebra $\mathcal{H}$. In
the case when $\mathcal{H}$ is a finite-dimensional Hopf algebra and ${%
\mathcal{B}\subseteq \mathcal{A}}$ is an $\mathcal{H}$-comodule algebra we
shall prove that $\mathcal{Z}$, the center of ${\mathcal{A}}$, is an $%
\mathcal{H}$-subcomodule. For a given Hopf bimodule $M$, our main purpose is
to show that, under some assumptions on $\mathcal{H}$ and $\mathcal{Z}^{%
\mathrm{co}\mathcal{H}}\subseteq \mathcal{Z}$, the homology groups $\mathrm{%
HH}_{\ast }(\mathcal{A},M)^{\mathrm{co}\mathcal{H}}$ and $\mathrm{HH}_{\ast
}(\mathcal{B},M^{\mathrm{co}\mathcal{H}})$ are isomorphic. A similar result
will be proved for cyclic homology.

\begin{proposition}
\label{pr:centrally_Galois}Let ${\mathcal{B}\subseteq }$ $\mathcal{A}$ be an
$\mathcal{H}$-comodule algebra. Let $\mathcal{Z}$ denote the center of $%
\mathcal{A}$ and set $\mathcal{Z}^{\prime }:=$ $\mathcal{Z}\bigcap \mathcal{B%
}.$

\begin{enumerate}
\item If $\mathcal{H}$ is commutative and finitely generated as an algebra
then $\mathcal{Z}$ is an $\mathcal{H}$-subcomodule of $\mathcal{A}.$

\item If $\mathcal{Z}$ is an $\mathcal{H}$-subcomodule of $\mathcal{A}$and $%
\mathcal{Z}^{\prime }\subseteq \mathcal{Z}$ is an
$\mathcal{H}$-Galois extension then $\mathcal{H}$ is commutative.
Let us assume, in addition, that $\mathcal{Z}^{\prime }\subseteq
\mathcal{Z}$ is a faithfully flat extension. Then
$\mathcal{A}^{\mathrm{co}\mathcal{H}}\subseteq \mathcal{A}$ is a
faithfully flat $\mathcal{H}$-Galois extension.
\end{enumerate}
\end{proposition}

\begin{proof}
(1) As $\mathcal{A}$ is an $\mathcal{H}$-comodule, $(\mathcal{A},\cdot )$ is
a left $\mathcal{H}^{\ast }$-module, where for $\alpha $ in $\mathcal{H}%
^{\ast }$ and $a$ in $\mathcal{A}$
\begin{equation}
\textstyle\alpha \cdot a={\sum }\alpha (a_{\left\langle 1\right\rangle
})a_{\left\langle 0\right\rangle }.  \label{ec:dual_action}
\end{equation}%
To prove that $\mathcal{Z}$ is an $\mathcal{H}$-subcomodule we must check
that $\mathcal{Z}$ is an $\mathcal{H}^{\ast }$-submodule. Let $\mathcal{H}%
^{\circ }$ denote the finite dual of $\mathcal{H}\ $(for the definition of
the finite dual of an algebra see \cite[Section 1.5]{DNR}). It is well-known
that $\mathcal{H}^{\circ }$ is an $S_{\mathcal{H}}$-invariant subalgebra of $%
\mathcal{H}^{\ast },$ so $\ $it has a canonical structure of Hopf
algebra. The comultiplication of $\mathcal{H}^{\circ }$ is uniquely
defined such that $\Delta (\alpha ):=\sum_{i=1}^{n}\alpha
_{i}^{\prime }\otimes \alpha _{i}^{\prime \prime }$ if and only if
\begin{equation*}
\alpha (xy)=\sum_{i=1}^{n}\alpha _{i}^{\prime }(x)\alpha _{i}^{\prime \prime
}(y),
\end{equation*}%
for all $x,y\in \mathcal{H}.$ Clearly, $\mathcal{A}$ is an $\mathcal{H}%
^{\circ }$-module. In fact $\mathcal{A}$ is an $\mathcal{H}^{\circ }$-module
algebra, that is
\begin{equation}
\textstyle\alpha \cdot (a^{\prime }a^{\prime \prime })={\sum }(\alpha
_{(1)}\cdot a^{\prime })(\alpha _{(2)}\cdot a^{\prime \prime }),
\label{mod_alg}
\end{equation}%
for $\alpha \in \mathcal{H}^{\circ }$ and $a^{\prime },$ $a^{\prime \prime
}\in \mathcal{A}.$ We now want to show that $\mathcal{Z}$ is an $\mathcal{H}%
^{\circ }$-submodule. For $\alpha \in \mathcal{H}^{\circ }$ and
$a\in \mathcal{Z}$, we get
\begin{align*}
(\alpha \cdot a)x& ={\textstyle\sum }(\alpha _{(1)}\cdot a)[\alpha
_{(2)}\cdot (S_{\mathcal{H}^{\circ }}\alpha _{(3)}\cdot x)] \\
& ={\textstyle\sum }\alpha _{(1)}\cdot \lbrack a(S_{\mathcal{H}^{\circ
}}\alpha _{(2)}\cdot x)] \\
& ={\textstyle\sum }\alpha _{(1)}\cdot \lbrack
(S_{\mathcal{H}^{\circ }}\alpha _{(2)}\cdot x)a],
\end{align*}%
where for the second equality we used (\ref{mod_alg}). By \cite[Corollary
2.3.17(ii)]{Ab}, $\mathcal{H}^{\circ }$ is cocommutative. Thus
\begin{align*}
{\textstyle\sum }\alpha _{(1)}\cdot \lbrack (S_{\mathcal{H}^{\circ }}\alpha
_{(2)}\cdot x)a]& ={\textstyle\sum }[(\alpha _{(1)}S_{\mathcal{H}^{\circ
}}\alpha _{(3)})\cdot x](\alpha _{(2)}\cdot a) \\
& ={\textstyle\sum }[(\alpha _{(1)}S_{\mathcal{H}^{\circ }}\alpha
_{(2)})\cdot x](\alpha _{(3)}\cdot a) \\
& =x(\alpha \cdot a).
\end{align*}%
By the foregoing computation, we conclude that $\alpha \cdot a\in \mathcal{Z}%
.$ Since $\mathcal{H}$ is finitely generated as an algebra it follows that $%
\mathcal{H}^{\circ }$ is dense in $\mathcal{H}^{\ast },$ with respect to the
finite topology, cf. \cite[Theorems 2.2.17 and 2.3.19 ]{Ab}. This means
that, for every $\alpha \in \mathcal{H}^{\ast }$ and every finite set $%
X\subseteq \mathcal{H}$ there is $\beta \in \mathcal{H}^{\circ }$
such that
$\alpha =\beta $ on $X.$ We can now prove that $\mathcal{Z}$ is an $\mathcal{%
H}^{\ast }$-submodule of $\mathcal{A}.$ Let $\alpha \in \mathcal{H}^{\ast }$
and $a\in \mathcal{A}.$ If $\rho (a)=\sum_{i=1}^{n}a_{i}\otimes h_{i},$ then
there is $\beta \in $ $\mathcal{H}^{\circ }$ such that $\alpha (h_{i})=\beta
(h_{i}),$ for every $i=1,\dots ,n.$ Thus%
\begin{equation*}
\alpha \cdot a=\sum_{i=1}^{n}\alpha (h_{i})a_{i}=\sum_{i=1}^{n}\beta
(h_{i})a_{i}=\beta \cdot a.
\end{equation*}%
It follows that $\alpha \cdot a\in \mathcal{Z},$ as $\beta \in \mathcal{H}%
^{\circ }$ and $a\in \mathcal{Z}.$

(2) Since $\mathcal{Z}$ is a subcomodule of $\mathcal{A}$, it follows that $%
\mathcal{Z}^{\mathrm{co}\mathcal{H}}=\mathcal{Z}^{\prime }$. The canonical
map $\beta _{\mathcal{Z}}:\mathcal{Z}\otimes _{\mathcal{Z}^{\prime }}%
\mathcal{Z\rightarrow Z}\otimes \mathcal{H}$, that corresponds to the $%
\mathcal{H}$-comodule algebra $\mathcal{Z}^{\prime }\subseteq \mathcal{Z}$,
is bijective by assumption. As $\mathcal{Z}$ is a commutative algebra, $%
\beta _{\mathcal{Z}}$ is a morphism of algebras and $\mathcal{Z}\otimes _{%
\mathcal{Z}^{\prime }}\mathcal{Z}$ is commutative. We conclude that $%
\mathcal{H}$ is commutative by remarking that $\mathcal{H}$ is a subalgebra
of $\mathcal{Z}\otimes \mathcal{H}$, which is commutative.

We now assume that $\mathcal{Z}^{\prime }\subseteq \mathcal{Z}$ is a
faithfully flat $\mathcal{H}$-Galois extension. For each $h\in
\mathcal{H}$ there are $a_{1}^{\prime },\dots ,a_{r}^{\prime }$ and
$a_{1}^{\prime \prime },\dots ,a_{r}^{\prime \prime }$ in
$\mathcal{Z}$ such that
\begin{equation}
\textstyle\beta _{\mathcal{Z}}(\sum_{i=1}^{r}a_{i}^{\prime }\otimes _{%
\mathcal{Z}^{\mathrm{co}\mathcal{H}}}a_{i}^{\prime \prime })=1\otimes h.
\label{eq:beta}
\end{equation}
Obviously, $\beta _{\mathcal{A}}(\sum_{i=1}^{r}a_{i}^{\prime }\otimes _{%
\mathcal{B}}a_{i}^{\prime \prime })=\beta _{\mathcal{Z}}(%
\sum_{i=1}^{r}a_{i}^{\prime }\otimes _{\mathcal{Z}^{\mathrm{co}\mathcal{H}%
}}a_{i}^{\prime \prime })=1\otimes h$ and $\beta _{\mathcal{A}}$ is
a morphism of left $\mathcal{A}$-modules. Thus $\beta
_{\mathcal{A}}$ is surjective too. Since $\mathcal{Z}^{\prime
}\subseteq \mathcal{Z}$ is faithfully flat it follows that
$\mathcal{Z}$ is injective as an $\mathcal{H}
$-comodule. By \cite[Lemma 4.1.]{SS} there is an $\mathcal{H}$-comodule map $%
\phi :\mathcal{H}\rightarrow \mathcal{Z}$ such that $\phi (1)=1.$ We may
regard $\phi $ as an $\mathcal{H}$-colinear map from $\mathcal{H}$ to $%
\mathcal{A},$ so $\mathcal{A}$ is injective as an $\mathcal{H}$-comodule.
Hence $\mathcal{B}\subseteq \mathcal{A}$ is a faithfully flat $\mathcal{H}$%
-Galois extension, cf. \cite[Theorem 4.10]{SS}.
\end{proof}

\begin{definition}
Let $\mathcal{H}$ be a commutative Hopf algebra. We say that an $\mathcal{H}$%
-comodule algebra ${\mathcal{B}\subseteq \mathcal{A}}$ is a \emph{centrally }%
$\mathcal{H}$-\emph{Galois extension} if the center $\mathcal{Z}$ of $%
\mathcal{A}$ is a subcomodule and $\mathcal{Z}^{\prime }\subseteq \mathcal{Z}
$ is a faithfully flat $\mathcal{H}$-Galois extension, where $\mathcal{Z}%
^{\prime }:=\mathcal{Z}^{\mathrm{co}\mathcal{H}}$.
\end{definition}

\begin{remark}
In the case when $\mathcal{H}$ is cosemisimple and finitely generated as an
algebra, an $\mathcal{H}$-comodule algebra $\mathcal{A}$ is centrally $%
\mathcal{H}$-Galois if and only if $\mathcal{Z}^{\prime }\subseteq
\mathcal{Z}$ is $\mathcal{H}$-Galois.
\end{remark}

\begin{noname}
Throughout the remaining part of this section we fix a commutative Hopf
algebra $\mathcal{H}$ and a centrally $\mathcal{H}$-Galois extension ${%
\mathcal{B}\subseteq }$ $\mathcal{A}$. We also fix a Hopf bimodule $M$ and a
left $\mathcal{H}$-comodule $V$.

We have seen that $\mathcal{B}\subseteq \mathcal{A}$ is $\mathcal{H}$%
-Galois, so $\mathrm{HH}{}_{n}(\mathcal{B},M)$ is a left $\mathcal{H}$%
-module. Our aim now is to give an equivalent description of this action. We
fix $h\in \mathcal{H}$ and we pick up $a_{1}^{\prime },\dots ,a_{r}^{\prime
} $ and $a_{1}^{\prime \prime },\dots ,a_{r}^{\prime \prime }$ in $\mathcal{Z%
}$ such that (\ref{eq:beta}) holds true. We now define $\lambda
_{n}^{h}(M):C_{n}(\mathcal{B},M)\rightarrow C_{n}(\mathcal{B},M)$ by%
\begin{equation*}
\textstyle\lambda _{n}^{h}(M)(m\otimes b^{1}\otimes \mathfrak{\cdots }%
\otimes b^{n})=\sum_{i=1}^{r}a_{i}^{\prime \prime }ma_{i}^{\prime }\otimes
b^{1}\otimes \mathfrak{\cdots }\otimes b^{n}.
\end{equation*}%
It is easy to see that $\lambda _{\ast }^{h}(M)$ is a morphism of complexes,
as $a_{i}^{\prime }$ and $a_{i}^{\prime \prime }$ are in the center of $%
\mathcal{A}$ for all $i=1,\dots ,r$. Let $\bar{\lambda}_{n}^{h}(M)$ be the
endomorphism of $\mathrm{HH}_{n}(\mathcal{B},M)$ induced by $\lambda
_{n}^{h}(M)$. Clearly, both $\lambda _{n}^{h}$ and $\bar{\lambda}_{n}^{h}$
are natural transformations.
\end{noname}

\begin{proposition}
\label{pr: trivial}Let $h\in \mathcal{H}.$ For a Hopf bimodule $M$ and $%
\omega \in \mathrm{HH}_{n}(\mathcal{B},M)$
\begin{equation*}
h\cdot \omega =\bar{\lambda}_{n}^{h}(M)(\omega ).
\end{equation*}%
If in addition $M$ is a symmetric $\mathcal{Z}$-bimodule then the
above action is trivial. In this case, for an injective left
$\mathcal{H}$-comodule $V$, the
action of $\mathcal{R}_{\mathcal{H}}$ on $\mathrm{HH}_{n}(\mathcal{B}%
,M\square _{\mathcal{H}}V)$ is trivial too.
\end{proposition}

\begin{proof}
Let $\mu _{\ast }^{h}$ be the natural transformations that lift the $%
\mathcal{H}$-action on $M_{\mathcal{B}},$ as in the proof of Proposition \ref%
{pr:HH_n(B,M)}. Thus, for $\omega $ in
$\mathrm{HH}_{n}(\mathcal{B},M)$
\begin{equation*}
\mu _{n}^{h}(M)(\omega )=h\cdot \omega .
\end{equation*}%
We shall prove by induction on $n$ that $\bar{\lambda}_{n}^{h}(M)=\mu
_{n}^{h}(M)$. For $n=0$ that is obvious, by construction of $\bar{\lambda}%
_{0}^{h}$ and the definition of the $\mathcal{H}^{{}}$-module structure in (%
\ref{eq:R-module}). Let us assume that $\bar{\lambda}_{n}^{h}(K)=\mu
_{n}^{h}(K),$ for any Hopf bimodule $K.$ Since ${\mathcal{B}}{\subseteq }%
\mathcal{A}$ is a faithfully flat $\mathcal{H}$-Galois extension, $U:=%
\mathcal{A}\otimes \mathcal{A}$ is a projective generator in $_{\mathcal{A}}%
\mathcal{M}_{\mathcal{A}}^{\mathcal{H}^{{}}}$. Thus, there is an exact
sequence
\begin{equation*}
0\longrightarrow K_{0}\longrightarrow L\longrightarrow M\longrightarrow 0
\end{equation*}%
in $_{\mathcal{A}}\mathcal{M}_{\mathcal{A}}^{\mathcal{H}^{{}}}$ such that $%
L\cong U^{(I)},$ where $I$ is a certain set. On the other hand, $\mathcal{A}$
is projective as a left and right ${\mathcal{B}}$-module, so $U$ is
projective as a ${\mathcal{B}}$-bimodule. Hence $\mathrm{HH}_{n}(\mathcal{B}%
,L)=0$, for $n>0.$ Consequently, $\delta _{n+1}:\mathrm{HH}_{n+1}(\mathcal{B}%
,M)\rightarrow \mathrm{HH}_{n}(\mathcal{B},K_{0})$ is injective. On the
other hand, by construction, $\mu _{\ast }^{h}$ is a morphism of $\delta $%
-functors. Thus%
\begin{equation}
\delta _{n+1}\circ \mu _{n+1}^{h}(M)=\mu _{n}^{h}(K_{0})\circ \delta _{n+1}.
\label{ec:delta}
\end{equation}%
Since the long exact sequence in homology is natural and $\lambda _{\ast
}^{h}$ is a natural morphism of complexes
\begin{equation}
\delta _{n+1}\circ \bar{\lambda}_{n+1}^{h}(M)=\bar{\lambda}%
_{n}^{h}(K_{0})\circ \delta _{n+1}.  \label{ec:Lambda}
\end{equation}%
Using relations (\ref{ec:delta}) and (\ref{ec:Lambda}), the induction
hypothesis and the fact that $\delta _{n+1}$ is injective one gets $\mu
_{n+1}^{h}(M)=\bar{\lambda}_{n+1}^{h}(M).$

Let us assume that $M$ is symmetric as a $\mathcal{Z}$-bimodule, i.e. $%
z\cdot m=m\cdot z,$ for any $z\in \mathcal{Z}$ and $m\in M.$ Thus%
\begin{equation*}
\textstyle\lambda _{n}^{h}(M)(m\otimes b^{1}\otimes \mathfrak{\cdots }%
\otimes b^{n})=\sum_{i=1}^{k}ma_{i}^{\prime }a_{i}^{\prime \prime }\otimes
b^{1}\otimes \mathfrak{\cdots }\otimes b^{n}=\varepsilon (h)m\otimes
b^{1}\otimes \mathfrak{\cdots }\otimes b^{n},
\end{equation*}%
where for the second equality we used \cite[Relation (5)]{JS}. Thus $\bar{%
\lambda}_{n}^{h}(M)(\omega )=\varepsilon (h)\omega .$ By the first part of
the proposition we deduce that the action of $\mathcal{H}$ on $\mathrm{HH}%
_{\ast }(\mathcal{B},M)$ is trivial. Finally, if $V$ is an injective left $\mathcal{H}%
^{{}}$-comodule, then there is an isomorphism
\begin{equation*}
\mathrm{HH}_{\ast }(\mathcal{B},M\square _{\mathcal{H}}V)\cong \mathrm{HH}%
_{\ast }(\mathcal{B},M)\square _{\mathcal{H}}V
\end{equation*}%
of $\mathcal{R}_{\mathcal{H}}$-modules. Note that $\mathrm{HH}_{\ast }(%
\mathcal{B},M)\square _{\mathcal{H}}V$ is a $\mathcal{R}_{\mathcal{H}}$%
-submodule of $\mathrm{HH}_{\ast }(\mathcal{B},M)\otimes V$. Hence,
the
action of $\mathcal{R}_{\mathcal{H}^{{}}}$ on $\mathrm{HH}_{\ast }(\mathcal{B%
},M)\square _{\mathcal{H}}V$ is also trivial.
\end{proof}

\begin{theorem}
\label{izo1}Let ${\mathcal{B}\subseteq \mathcal{A}}$ be a centrally $%
\mathcal{H}$-Galois extension, where $\mathcal{H}$ is a
finite-dimensional Hopf algebra over a field $\mathbb{K}$ such
that $\dim \mathcal{H}$ is not
zero in $\mathbb{K}$. Let $M$ be a Hopf bimodule which is symmetric as a $%
\mathcal{Z}$-bimodule. If $V$ is a left $\mathcal{H}$-comodule then there
are isomorphisms of $\mathcal{Z}^{\prime }$-modules
\begin{equation}
\mathrm{HH}{}_{\ast }(\mathcal{A},M)\square _{\mathcal{H}}V\simeq \mathrm{HH}%
{}_{\ast }(\mathcal{B},M\square _{\mathcal{H}}V).  \label{eq:izo1}
\end{equation}
\end{theorem}

\begin{proof}
We have already noticed that $\mathcal{B}\subseteq \mathcal{A}$ is a
faithfully flat $\mathcal{H}$-Galois extension and that $\mathrm{H}_{\ast }(\mathcal{B%
},M\square _{\mathcal{H}}V)$ is a trivial $\mathcal{R}_{\mathcal{H}}$%
-module. The isomorphism of left $\mathcal{Z}^{\prime }$-modules(\ref%
{eq:izo1}) follows by applying Theorem \ref{te:iso}.
\end{proof}

\begin{corollary}
\label{iso2} Keeping the notation and the assumptions from the preceding
theorem, there are isomorphisms of $\mathcal{Z}^{\prime }$-modules
\begin{equation*}
\mathrm{HH}{}_{\ast }(\mathcal{A},M)^{\mathrm{co}{}\mathcal{H}}\simeq
\mathrm{HH}{}_{\ast }(\mathcal{B},M^{\mathrm{co}{}\mathcal{H}}).
\end{equation*}
\end{corollary}

\begin{proof}
Take $V=\mathbb{K}$ in Theorem \ref{izo1} and note that $(-)^{\mathrm{co}{}%
\mathcal{H}}\cong (-)\square _{\mathcal{H}}\mathbb{K}.$
\end{proof}

\begin{remark}
A faithfully flat $\mathcal{H}$-Galois extension of commutative
algebras is centrally Hopf-Galois. Thus the isomorphisms in the
preceding corollary exist for such an extension, provided that
$\mathcal{H}$ is finite-dimensional and $\dim \mathcal{H}$ is not
zero in $\mathbb{K}.$
\end{remark}

\begin{corollary}
\label{iso3}Let $G$ be a finite group of automorphisms of an algebra ${%
\mathcal{A}}$ over a field $\mathbb{K}$ such that$\ $the order of
$G$ is not zero in $\mathbb{K}$. Let $\mathcal{Z}$ denote the center
of $\mathcal{A}
$ and let $M$ be an $(G,\mathcal{A})$-Hopf bimodule. If $\mathcal{Z}%
^{G}\subseteq \mathcal{Z}$ is a $G$-Galois extension then there are
isomorphisms of $\mathcal{Z}^{G}$-modules
\begin{equation*}
\mathrm{HH}{}_{\ast }(\mathcal{A},M)^{G}\simeq \mathrm{HH}{}_{\ast }(%
\mathcal{A}^{G},M^{G}).
\end{equation*}
\end{corollary}

\begin{proof}
Apply Corollary \ref{iso2} for $\mathcal{H}:=(\mathbb{K}G)^{\ast }.$ Note
that, for a $(\mathbb{K}G)^{\ast }$-comodule $X$, we have $X^{\mathrm{co}(%
\mathbb{K}G)^{\ast }}=X^{G}$, cf. \S \ref{nn:G-Galois}.
\end{proof}

\begin{remark}
Note that the proof of the previous corollary works only if the
order of $G$ is not invertible in $\mathbb{K},$ as $\mathbb{K}$
must be injective as $(\mathbb{K}G)^{\ast }$-comodule in order to
apply Theorem \ref{izo1}. On the other hand  the isomorphisms in
\cite[\S 6]{Lo2}  hold true without any assumption on the characteristic of $%
\mathbb{K}.$
\end{remark}

\begin{theorem}
\label{izo2}Let ${\mathcal{B}\subseteq \mathcal{A}}$ be a centrally Galois
extension over a finite-dimensional Hopf algebra $\mathcal{H}$ such that $%
\dim \mathcal{H}$ is not zero in $\mathbb{K}$. Let $M$ be a Hopf bimodule
which is symmetric as a $\mathcal{Z}$-bimodule. Then, there are isomorphisms
of $\mathcal{Z}$-modules and $\mathcal{H}$-comodules
\begin{equation*}
\mathrm{HH}{}_{\ast }(\mathcal{A},M)\simeq \mathcal{Z}\otimes _{\mathcal{Z}%
^{\prime }}\mathrm{HH}{}_{\ast }(\mathcal{B},M^{\mathrm{co}\mathcal{H}}).
\end{equation*}
\end{theorem}

\begin{proof}
Since $\mathcal{Z}^{\prime }\subseteq \mathcal{Z}$ is a faithfully flat $%
\mathcal{H}^{{}}$-Galois extension the categories $_{\mathcal{Z}}\mathfrak{M}%
^{\mathcal{H}}$ and ${}_{\mathcal{Z}^{\prime }}\mathfrak{M}$ are equivalent,
cf. \cite[Theorems 4.9 and 4.10]{SS}. More precisely,
\begin{equation*}
\mathcal{Z}\otimes _{\mathcal{Z}^{\prime }}(-):{}_{\mathcal{Z}^{\prime }}%
\mathfrak{M}\longrightarrow {}_{\mathcal{Z}}\mathfrak{M}^{\mathcal{H}}
\end{equation*}%
is an equivalence of categories, whose inverse is the functor $X\mapsto X^{%
\mathrm{co}H^{{}}}$. Thus,
\begin{equation*}
\mathrm{H}_{\ast }(\mathcal{A},M)\cong \mathcal{Z}\otimes _{\mathcal{Z}%
^{\prime }}\mathrm{H}_{\ast }(\mathcal{A},M)^{\mathrm{co}{}\mathcal{H}}.
\end{equation*}%
We conclude de proof in view of Corollary \ref{iso2}.
\end{proof}

\begin{noname}
Recall that $\mathcal{B}\subseteq \mathcal{A}$ is a centrally $\mathcal{H}$%
-Galois extension. In particular, $\mathcal{H}$ is commutative. Therefore,
the map $t_{n}:C_{n}(\mathcal{A},\mathcal{A})\longrightarrow C_{n}(\mathcal{A%
},\mathcal{A})$ given by
\begin{equation*}
t_{n}(a^{0}\otimes a^{1}\otimes \cdots \otimes a^{n})=a^{n}\otimes
a^{0}\otimes \cdots \otimes a^{n-1}.
\end{equation*}%
is a morphism of right $\mathcal{H}$-comodules, where $\mathcal{H}$ coacts
on $\mathcal{A}^{\otimes n+1}$ as in (\ref{nn:FunctorH}). Cyclic homology of
$\mathcal{A}$, denoted by $\mathrm{HC}_{\ast }(\mathcal{A}),$ is defined as
the homology of the total complex of the bicomplex $\mathrm{CC}_{\ast \ast }(%
\mathcal{A})$, see \cite[Definition 9.6.6]{We}. As the operator $t_{n}$ is $%
\mathcal{H}$-colinear for every $n,$ it follows that
$\mathrm{CC}_{\ast
\ast }(\mathcal{A})$ is a bicomplex in the category of right $\mathcal{H}$%
-comodules. Thus $\mathrm{HC}_{n}(\mathcal{A})$ is an $\mathcal{H}$-comodule
too.

We can now prove the following.
\end{noname}

\begin{theorem}
\label{izo3}Let ${\mathcal{B}\subseteq \mathcal{A}}$ be a centrally $%
\mathcal{H}$-Galois extension, where $\mathcal{H}$ is a
finite-dimensional Hopf algebra such that $\dim \mathcal{H}$ is
not zero in $\mathbb{K}$. Then
\begin{equation}
\mathrm{HC}_{n}(\mathcal{A})^{\mathrm{co}{}\mathcal{H}}\cong \mathrm{HC}_{n}(%
\mathcal{B}).  \label{iso cyc}
\end{equation}
\end{theorem}

\begin{proof}
The case $n=0$ is obvious, in view of Corollary \ref{iso2} and of the fact
that cyclic homology and Hochschild homology are equal in degree zero. For
each right $\mathcal{H}$-comodule $X$ the natural transformation
\begin{equation*}
\nu (X):X^{\mathrm{co}\mathcal{H}}\longrightarrow X\square _{\mathcal{H}}%
\mathbb{K},\qquad \nu (X)(x):=x\otimes 1
\end{equation*}%
is an isomorphism. Since $\mathcal{H}$ is commutative and the characteristic
of $\mathbb{K}$ does not divide the dimension of $\mathcal{H},$ we deduce
that $\mathcal{H}$ is cosemisimple. Hence $\mathbb{K}$ is an injective
comodule. Thus the functor that maps a right $\mathcal{H}$-comodule $X$ to $%
X\square _{\mathcal{H}}\mathbb{K}$ is exact. Consequently, by applying the
functor $(-)^{\mathrm{co}\mathcal{H}}$ to Connes' exact sequence \cite[%
Proposition 9.6.11]{We}, we get the exact sequence on the top of the
following diagram.
\begin{equation*}
\xymatrix{ \widetilde{\mathrm{HC}}_{n}(\mathcal{A})
\ar[r]^{\widetilde{B}} &
\widetilde{\mathrm{HH}}_{n+1}(\mathcal{A}) \ar[r]^{\widetilde{I}}
& \widetilde{\mathrm{HC}}_{n+1}(\mathcal{A})
\ar[r]^{\widetilde{S}} &
\widetilde{\mathrm{HC}}_{n-1}(\mathcal{A}) \ar[r]^{\widetilde{B}}
&\widetilde{\mathrm{HH}}_{n}(\mathcal{A}) \\
\mathrm{HC}_{n}(\mathcal{B}) \ar[r]_{B} \ar[u]^{\cong}&
\mathrm{HH}_{n+1}(\mathcal{B}) \ar[r]_{I} \ar[u]^{\cong}&
\mathrm{HC}_{n+1}(\mathcal{B})\ar[r]_{S} \ar[u]&
\mathrm{HC}_{n-1}(\mathcal{B}) \ar[r]_{B} \ar[u]_{\cong}
&\mathrm{HH}_{n}(\mathcal{B}) \ar[u]_{\cong}}
\end{equation*}%
Here, $\widetilde{\mathrm{HC}}_{\ast }(\mathcal{A})$ and $\widetilde{\mathrm{%
HH}}_{\ast }(\mathcal{A})$ denote $\mathrm{HC}_{\ast }(\mathcal{A})^{\mathrm{%
co}\mathcal{H}}$ and $\mathrm{HH}_{\ast }(\mathcal{A})^{\mathrm{co}\mathcal{H%
}}$, respectively. Note that by induction hypothesis, the first and
the
fourth vertical arrows are isomorphisms. Furthermore, by taking $M=\mathcal{A%
}$ in Corollary \ref{iso2}, we get that the second and the fifth
vertical maps are isomorphisms. Thus, by $5$-Lemma \cite[p.13]{We}
the vertical map in the middle is also an isomorphism.
\end{proof}

We conclude this paper showing that, under some extra assumptions, Ore
extensions provide non-trivial examples of centrally Hopf-Galois extensions.
To define an Ore extension of a $\mathbb{K}$-algebra $\mathcal{A}$, we need
an algebra automorphism $\sigma :\mathcal{A}\rightarrow \mathcal{A}$ and a $%
\sigma$-derivation $\delta :\mathcal{A}\rightarrow \mathcal{A}$. Recall that
$\delta$ is a $\sigma$-derivation if, for $a$ and $b$ in $\mathcal{A}$,
\begin{equation*}
\delta (ab)=\sigma(a)\delta (b)+\delta (a)b.
\end{equation*}
For $\sigma$ and $\delta$ as above one defines a new algebra $\mathcal{A}%
[X,\sigma ,\delta ]$, the Ore extension of $\mathcal{A}$. As a left $%
\mathcal{A}$-module, $\mathcal{A}[X,\sigma ,\delta ]$ is free with basis $%
\{1,X,X^{2},\dots \}$ and its multiplication is the unique left $\mathcal{A}$%
-linear morphism such that $X^{n}X^{m}=X^{n+m}$ and
\begin{equation}  \label{Ore}
Xa=\sigma (a)X+\delta (a).
\end{equation}%
For simplicity, we shall denote the Ore extension $\mathcal{A}[X,\sigma
,\delta ]$ by $\mathcal{T}$.

We now assume, in addition, that $\mathcal{A}$ is an $\mathcal{H}$-comodule
algebra and that $\sigma$ and $\delta$ are morphisms of comodules. Set $%
\mathcal{B}:=\mathcal{A}^{\mathrm{co}\mathcal{H}}$. Since $\sigma $ and $%
\delta $ are morphisms of $\mathcal{H}$-comodules they map $\mathcal{B}$
into $\mathcal{B}.$ We still denote the restrictions of these maps to $%
\mathcal{B}$ by $\sigma $ and $\delta .$ Clearly, $\delta $ can be regarded
as a $\sigma $-derivation of $\mathcal{B}$, so we can construct the Ore
extension $\mathcal{S}:=\mathcal{B}[X,\sigma ,\delta ]$.

\begin{lemma}
The comodule structure map $\rho_{\mathcal{A}}:\mathcal{A}\longrightarrow%
\mathcal{A}\otimes\mathcal{H}$ can be extended in a unique way to an $%
\mathcal{H}$-coaction $\rho _{\mathcal{T}}$ on $\mathcal{T}$ such that, for $%
a\in \mathcal{A}$ and $n\in \mathbb{N},$
\begin{equation}
\textstyle \rho _{\mathcal{T}}(aX^{n})=\sum a_{\left\langle 0\right\rangle
}X^{n}\otimes a_{\left\langle 1\right\rangle }.  \label{ec:ro}
\end{equation}%
With respect to this coaction the subalgebra of coinvariant elements in $%
\mathcal{T}$ is $\mathcal{S}$.
\end{lemma}

\begin{proof}
For $n\in \mathbb{N}$ and $0\leq k\leq n$ let $f_{k}^{(n)}$ be the
non-commutative polynomial in $\sigma $ and $\delta $ with coefficients in
the prime subfield of $\mathbb{K}$ such that
\begin{equation}  \label{recurenta}
\textstyle X^{n}a=\sum_{k=0}^{n}f_{k}^{(n)}(a)X^{k}.
\end{equation}%
Let us put $f_{-1}^{(n)}=f_{n+1}^{(n)}=0$. Thus, by multiplying to the left
both sides of (\ref{recurenta}) by $X$ and using (\ref{Ore}), for $0\leq
k\leq n+1$, we get
\begin{equation*}
f_{k}^{(n+1)}=\sigma f_{k-1}^{(n)}+\delta f_{k}^{(n)},
\end{equation*}%
For $a^{0},\dots ,a^{n}$ in $\mathcal{A}$ we now define
\begin{equation*}
\textstyle \rho _{\mathcal{T}}(\sum_{i=0}^{n}a^{i}X^{i})=\sum_{i=0}^{n}a_{%
\left\langle 0\right\rangle }^{i}X^{i}\otimes a_{\left\langle 1\right\rangle
}^{i}.
\end{equation*}%
Clearly, $\rho _{\mathcal{T}}$ defines a coaction of $\mathcal{H}$ on $%
\mathcal{T}$ and verifies the identity (\ref{ec:ro}). We have to prove that $%
\rho _{\mathcal{T}}$ is a morphism of algebras, i.e. $\rho _{\mathcal{T}%
}(fg)=\rho _{\mathcal{T}}(f)\rho _{\mathcal{T}}(g)$ for any $f,g\in \mathcal{%
T}$. In fact, it is enough to prove this equality for $f=X^{n}$ and
$g=a,$ with $a\in \mathcal{A}$ and $n\in \mathbb{N}^{\ast }$. Since
$f_{k}^{(n)}$ are non-commutative polynomials in $\sigma $ and
$\delta $ and these maps are morphism of $\mathcal{H}$-modules, it
follows that $f_{k}^{(n)}$ are also $\mathcal{H}$-colinear. Hence,
\begin{align*}
\rho _{\mathcal{T}}(X^{n}a)&=\rho _{\mathcal{T}}(\textstyle%
\sum_{k=0}^{n}f_{k}^{(n)}(a)X^{k}) \\
&=\textstyle\sum_{k=0}^{n}f_{k}^{(n)}(a)_{\left\langle 0\right\rangle
}X^{k}\otimes f_{k}(a)_{\left\langle 1\right\rangle } \\
&=\textstyle\sum_{k=0}^{n}f_{k}^{(n)}(a_{\left\langle 0\right\rangle
})X^{k}\otimes a_{\left\langle 1\right\rangle }.
\end{align*}%
On the other hand,
\begin{equation*}
\textstyle \rho _{\mathcal{T}}(X^{n})\rho _{\mathcal{A}}(a)=\sum
X^{n}a_{\left\langle 0\right\rangle }\otimes a_{\left\langle 1\right\rangle
}=\sum_{k=0}^{n}f_{k}^{(n)}(a_{\left\langle 0\right\rangle })X^{k}\otimes
a_{\left\langle 1\right\rangle }=\rho _{\mathcal{T}}(X^{n}a).
\end{equation*}%
Obviously $\rho _{\mathcal{T}}\ $is unital. Thus $\mathcal{T}$ is an $%
\mathcal{H}$-comodule algebra. It remains to prove that $\mathcal{T}^{%
\mathrm{co}\mathcal{H}}=\mathcal{S}.$ For this, we fix a basis $\{h_{j}\mid
j\in J\}$ on $\mathcal{H}.$ We may assume that there is $j_{0}\in J$ such
that $h_{j_{0}}=1.$ Let us take $f=\sum_{i=0}^{n}a^{i}X^{i}\ $in $\mathcal{T}
$ and write $\rho (a^{i})=\sum_{j\in J}a_{j}^{i}\otimes h_{j}.$ Therefore,
\begin{equation*}
\textstyle \rho (f)=\sum_{i=0}^{n}\sum_{j\in J}a_{j}^{i}X^{i}\otimes h_{j}.
\end{equation*}%
It follows that $f\in \mathcal{T}^{\mathrm{co}\mathcal{H}}$ if and only if $%
\sum_{i=0}^{n}a_{j}^{i}X^{i}=\delta _{j,j_{0}}\sum_{i=0}^{n}a^{i}X^{i}.$
Thus, $f$ is $\mathcal{H}$-coinvariant if and only if
\begin{equation*}
\textstyle \rho (a^{i})=\sum_{j\in J}\delta _{j,j_{0}}a^{i}\otimes
h_{j}=a^{i}\otimes 1,
\end{equation*}%
for all $i=0,\dots ,n.$ We deduce that $f\in \mathcal{T}^{\mathrm{co}%
\mathcal{H}}$ if and only if $f\in \mathcal{S}$.
\end{proof}

\begin{lemma}
Let $f=\sum_{i=0}^{n}{a^{i}X^{i}}$ be an element in $\mathcal{T}$ and $%
a^{-1}=a^{n+1}=0.$ Then $f$ is in $\mathcal{Z}(\mathcal{T})$, the center of $%
\mathcal{T}$, if and only if %
\begin{eqnarray}
&\sum_{k=i}^{n}{a}^{k}{f_{i}^{(k)}(a)}=aa^{i},&\quad \text{for } i=0,\dots ,n%
\text{ and } a\in \mathcal{A},  \label{B1} \\
&\sigma (a^{i})+\delta (a^{i+1})=a^{i},&\quad \text{for } i=-1,0,\dots
,n+1.\   \label{B2}
\end{eqnarray}
\end{lemma}

\begin{proof}
As an algebra, $\mathcal{T}$ is generated by $\mathcal{A}$ and $X$. Hence, $%
f\ $is central if and only if $Xf=fX$ and $af=fa,$ for all $a\in \mathcal{A}$%
. We get
\begin{equation*}
\textstyle fa=\sum_{k=0}^{n}a^{k}X^{k}a=\sum_{k=0}^{n}%
\sum_{i=0}^{k}a^{k}f_{i}^{(k)}(a)X^{i}=\sum_{i=0}^{n}\left(
\sum_{k=i}^{n}a^{k}f_{i}^{(k)}(a)\right) X^{i}.
\end{equation*}%
We deduce that $fa=af$ and (\ref{B1}) are equivalent. On the other hand, $%
fX=Xf$ is equivalent to
\begin{equation*}
\textstyle \sum_{i=0}^{n}\sigma (a^{i})X^{i+1}+\sum_{i=0}^{n}\delta
(a^{i})X^{i}=\sum_{i=0}^{n}a^{i}X^{i+1}.
\end{equation*}%
In conclusion, $fX=Xf$ and (\ref{B2}) are equivalent.
\end{proof}

\begin{corollary}
Let $\mathcal{A}^{\sigma }=\{a\mid \sigma (a)=a\}$ and
$\mathcal{A}^{\delta }=\{a\mid \delta (a)=0\}.$ If $\mathcal{Z}$ is
the center of $\mathcal{A}$ then
\begin{equation*}
\mathcal{Z}(\mathcal{T})\cap \mathcal{A}=\mathcal{A}^{\sigma }\cap \mathcal{A%
}^{\delta }\cap \mathcal{Z}.
\end{equation*}
\end{corollary}

\begin{proof}
We regard $\mathcal{A}\ $as a subalgebra of $\mathcal{T}$. Thus,
$a^{0}\in \mathcal{A}$ is in the center of $\mathcal{T}$ if and only
if for any $a\in \mathcal{A}$ we have
\begin{equation*}
\delta (a^{0})=0,\quad \sigma (a^{0})=a^{0},\quad aa^{0}=a^{0}f_{0}^{(0)}(a).
\end{equation*}%
Since $f_{0}^{(0)}=\mathrm{Id}_{A}$, we get $\mathcal{Z}(\mathcal{T})\cap
\mathcal{A}=\mathcal{Z}\cap \mathcal{A}^{\sigma }\cap \mathcal{A}^{\delta }$.
\end{proof}

\begin{theorem}
\label{te:example}Let $\mathcal{B}\subseteq \mathcal{A}$ be an $\mathcal{H}$%
-comodule algebra, where $\mathcal{H}$ is a commutative finite-dimensional
Hopf algebra over a field of characteristic zero. Let $\sigma :\mathcal{A}%
\rightarrow \mathcal{A}$ be an algebra map and  $\delta :\mathcal{A}%
\rightarrow \mathcal{A}$ be a $\sigma$-derivation. Assume that both $\sigma$
and $\delta$ are morphisms of $\mathcal{H}$-comodules. Let $\mathcal{T}:=%
\mathcal{A}[X,\sigma ,\delta ]$ and $\mathcal{S}:=\mathcal{B}[X,\sigma
,\delta ]$.

\begin{enumerate}
\item The center $\mathcal{Z}$ of $\mathcal{A}$ and
$\mathcal{A}^{\sigma }\cap \mathcal{A}^{\delta }\cap \mathcal{Z}$
are $\mathcal{H}$-comodule
subalgebras of $\mathcal{A}.$ The algebra of coinvariant elements in $%
\mathcal{A}^{\sigma }\cap \mathcal{A}^{\delta }\cap \mathcal{Z}$ is $%
\mathcal{B}^{\sigma }\cap \mathcal{B}^{\delta }\cap \mathcal{Z}.$

\item If $\mathcal{B}^{\sigma }\cap \mathcal{B}^{\delta }\cap \mathcal{Z}%
\subseteq \mathcal{A}^{\sigma }\cap \mathcal{A}^{\delta }\cap \mathcal{Z}$
is an $\mathcal{H}$-Galois extension then the extension $\mathcal{S}%
\subseteq \mathcal{T}$ is a centrally $\mathcal{H}$-Galois extension.
\end{enumerate}
\end{theorem}

\begin{proof}
(1) By Proposition \ref{pr:centrally_Galois} (2), $\mathcal{Z}$ is an $%
\mathcal{H}$-subcomodule of $\mathcal{A}$ as $\mathcal{H}$ is
finite-dimensional and commutative. Since $\sigma$ and $\delta $ are
morphisms
of $\mathcal{H}$-comodules it follows that $\mathcal{A}^{\sigma }=\mathrm{ker%
}(\sigma -\mathrm{Id}_{\mathcal{A}})$ and $\mathcal{A}^{\delta }=\mathrm{ker}%
\delta $ are $\mathcal{H}$-subcomodules of $\mathcal{A}$. We deduce that $%
\mathcal{A}^{\sigma }\cap \mathcal{A}^{\delta }\cap \mathcal{Z}$ is an $%
\mathcal{H}$-comodule algebra. Its subalgebra of coinvariant elements is
\begin{equation*}
\lbrack \mathcal{A}^{\sigma }\cap \mathcal{A}^{\delta }\cap \mathcal{Z}]^{%
\mathrm{co}\mathcal{H}}=\mathcal{A}^{\sigma }\cap \mathcal{A}^{\delta }\cap
\mathcal{Z}\cap \mathcal{B}=\mathcal{B}^{\sigma }\cap \mathcal{B}^{\delta
}\cap \mathcal{Z}.
\end{equation*}

(2) Again by Proposition \ref{pr:centrally_Galois} (2), the center $\mathcal{%
Z}(\mathcal{T})$ of $\mathcal{T}$ is an $\mathcal{H}$-subcomodule of $%
\mathcal{T}$. Since $\mathcal{H}$ is commutative and finite-dimensional over
a field of characteristic zero, we deduce that $\mathcal{H} $ is
cosemisimple. Hence, $\mathcal{Z}(\mathcal{T})$ is injective as an $\mathcal{H%
}$-comodule. In view of \cite[Theorem 4.10]{SS}, to prove that $\mathcal{Z}(%
\mathcal{T})\cap \mathcal{S}\subseteq \mathcal{Z}(\mathcal{T})$ is $\mathcal{%
H}$-Galois and faithfully flat, we have to show that the canonical map
\begin{equation*}
\beta _{\mathcal{Z}(\mathcal{T})}:\mathcal{Z}(\mathcal{T})\otimes _{\mathcal{%
Z}(\mathcal{T})\cap \mathcal{S}}\mathcal{Z}(\mathcal{T})\longrightarrow
\mathcal{Z}(\mathcal{T})\otimes \mathcal{H}
\end{equation*}%
is surjective. Proceeding as in the proof of Proposition \ref%
{pr:centrally_Galois} (3) it is enough to show that $1\otimes h$ is in the
image of $\beta _{\mathcal{Z}(\mathcal{T})},$ for every $h\in \mathcal{H}.$
Let $\mathcal{Z}^{\prime }:=\mathcal{A}^{\sigma }\cap \mathcal{A}^{\delta
}\cap \mathcal{Z}.$ By assumption, the canonical map
\begin{equation*}
\beta _{\mathcal{Z}^{\prime }}:\mathcal{Z}^{\prime }\otimes _{\mathcal{Z}%
^{\prime }\cap \mathcal{B}}\mathcal{Z}^{\prime }\longrightarrow \mathcal{Z}%
^{\prime }\otimes \mathcal{H}
\end{equation*}%
is bijective. Thus, there are $a_{1}^{\prime },\dots ,a_{r}^{\prime }$ and $%
a_{1}^{\prime \prime },\dots, a_{r}^{\prime \prime }$ in $\mathcal{Z}%
^{\prime }$ such that
\begin{equation*}
\textstyle \beta _{\mathcal{Z}^{\prime }}(\sum_{i=1}^{r}a_{i}^{\prime
}\otimes _{\mathcal{Z}^{\prime }\cap \mathcal{B}}a_{i}^{\prime \prime
})=1\otimes h.
\end{equation*}%
By the previous corollary, $\mathcal{Z}^{\prime }$ is an $\mathcal{H}$%
-submodule of $\mathcal{Z}(\mathcal{T})$. Therefore,
\begin{equation*}
\textstyle \beta _{\mathcal{Z}(\mathcal{T})}(\sum_{i=1}^{n}a_{i}^{\prime
}\otimes _{\mathcal{Z}(\mathcal{T})\cap \mathcal{S}}a_{i}^{\prime \prime
})=\beta _{\mathcal{Z}^{\prime }}(\sum_{i=1}^{n}a_{i}^{\prime }\otimes _{%
\mathcal{Z}^{\prime }\cap \mathcal{B}}a_{i}^{\prime \prime })=1\otimes h.
\end{equation*}%
Hence, the theorem is completely proven.
\end{proof}

A more concrete example can be obtained as follows. Let $\mathbb{K}\subseteq
\mathcal{K}\subseteq \mathcal{A}$ be field extensions such that $\mathcal{K}%
\subseteq \mathcal{A}$ is finite, separable and normal of Galois group $G. $
We assume that $G=NH,$ where $H$ and $N$ are subgroups in $G$ such that $%
N\bigcap H=\left\{ 1\right\} $ and $N$ is generated by a central element $%
\sigma $ in $G.$ We set $\mathcal{B}:=\mathcal{A}^{H}.$ We wish to prove
that this setting fulfils the conditions in the preceding theorem, to get
the following.

\begin{corollary}
With the above notation, $\mathcal{B}\left[ X,\sigma ,0\right] \subseteq
\mathcal{A}\left[ X,\sigma ,0\right] $ is a centrally $\left( \mathbb{K}%
H\right) ^{\ast }$-Galois extensions.
\end{corollary}

\begin{proof}
In order to apply Theorem \ref{te:example}, we have to check that
$\sigma $
is a morphism of $\left( \mathbb{K}H\right) ^{\ast }$-comodules and that $%
\mathcal{B}^{\sigma }\subseteq \mathcal{A}^{\sigma }$ is a $\left( \mathbb{K}%
H\right) ^{\ast }$-Galois extension. The former condition is
equivalent to the fact that $\sigma $ is a morphism of
$\mathbb{K}H$-modules, which in our case means that $\sigma
h=h\sigma $ for any $h$ in $H.$ Trivially this equality is satisfied
as, by assumption, $\sigma $ is central in $G$.
Furthermore,%
\begin{equation*}
\mathcal{B}^{\sigma }=\left( \mathcal{A}^{H}\right) ^{\sigma }=\left(
\mathcal{A}^{H}\right) ^{N}=\mathcal{A}^{HN}=\mathcal{K}.
\end{equation*}%
A similar computation yields us $\left( \mathcal{A}^{\sigma }\right) ^{H}=%
\mathcal{A}^{NH}=\mathcal{B}^{\sigma }.$ On the other hand, since
$N\bigcap H=\left\{ 1\right\}$ one can embed $H$ into the group of
field automorphisms of $\mathcal{A}^{\sigma }$ via the restriction
map $u\mapsto u|_{\mathcal{A}^{\sigma }}.$ By Artin's Lemma,
$\mathcal{B}^{\sigma
}\subseteq \mathcal{A}^{\sigma }$ is separable and normal of Galois group $%
H. $ We have noticed in \S \ref{nn:G-Galois} that a finite field extension
is $\left( \mathbb{K}H\right) ^{\ast }$-Galois if and only if it is
separable and normal of Galois group $H.$ In conclusion, the second
requirement is also satisfied.
\end{proof}
\noindent\textbf{Acknowledgements.} The authors thank the referee
for his valuable comments and suggestions.

\bigskip

\noindent
\begin{minipage}{148mm}\sc\footnotesize Universit\' e de
Haute Alsace, Laboratoire de Math\'{e}matiques Informatique et
Applications, 4, Rue des Fr\`eres Lumi\`ere,
68093 Mulhouse Cedex, France\\
{\it E--mail address}: {\tt
Abdenacer.Makhlouf@uha.fr}\bigskip

University of Bucharest, Faculty of Mathematics, Str.
Academiei 14,
RO-70109, Bucharest, Romania\\
{\it E--mail address}: {\tt drgstf@gmail.com}
\end{minipage}

\end{document}